\documentclass[a4paper]{article}
\usepackage{cite}
\usepackage[all]{xy}
\usepackage[latin1]{inputenc}        
\usepackage[dvips]{graphics,graphicx}
\usepackage{amsfonts,amssymb,amsmath,color,mathrsfs, amstext}
\usepackage{amsbsy, amsopn, amscd, amsxtra, amsthm,authblk}
\usepackage{enumerate,algorithmicx,algorithm,enumitem}
\usepackage{algpseudocode}
\usepackage{upref}
\usepackage{geometry,wrapfig}
\usepackage{amsthm,amsmath,amssymb}
\usepackage{relsize}
\usepackage{mathrsfs,mathtools}
\usepackage{bm,bbm}
\usepackage{ulem}
\usepackage{exscale}
\usepackage{tabularx}
\usepackage{threeparttable,multirow,dcolumn,booktabs}
\usepackage{algorithmicx,algorithm}
\usepackage{tabularx}
\usepackage{caption}
\usepackage{subfigure}
\usepackage[displaymath]{lineno}
\usepackage{float}
\usepackage{bm}
\usepackage{caption}
\usepackage{yhmath}
\usepackage{booktabs}
\usepackage{multirow}
\usepackage{makecell}
\usepackage[colorlinks,
            linkcolor=red,
            anchorcolor=red,
            citecolor=red
            ]{hyperref}

\usepackage{cleveref}

\usepackage[table,dvipsnames]{xcolor}

\definecolor{cgray}{gray}{0.80}
\newcolumntype{a}{>{\columncolor{cgray}}c}

\def\R{\mathbbm{R}}
\def\N{\mathbbm{N}}

\newcommand{\calD}{\mathcal{D}}
\newcommand{\cC}{\mathcal C}

\numberwithin{equation}{section}
\usepackage{epstopdf}
\usepackage{pdfpages}

\graphicspath{{./Figs/}}

\begin{document}

\title{TGPT-PINN: Nonlinear model reduction with transformed GPT-PINNs}

\author{
Yanlai Chen\footnote{Department of Mathematics, University of Massachusetts Dartmouth, North Dartmouth, MA 02747. Email: {\tt{yanlai.chen@umassd.edu}}. Y. Chen is partially supported by National Science Foundation grant DMS-2208277 and by the UMass Dartmouth MUST program, N00014-22-1-2012, established by Dr. Ramprasad Balasubramanian via sponsorship from the Office of Naval Research.}, \quad  
Yajie Ji\footnote{School of Mathematical Sciences, Shanghai Jiao Tong University, Shanghai 200240, China. Email: {\tt jiyajie595@sjtu.edu.cn}. },\quad
Akil Narayan\footnote{Scientific Computing and Imaging Institute and Department of Mathematics, University of Utah, Salt Lake City, UT 84112. Email: {\tt{akil@sci.utah.edu}}. A. Narayan is partially supported by NSF DMS-1848508 and AFOSR FA9550-23-1-0749.}, \quad
Zhenli Xu\footnote{School of Mathematical Sciences, CMA-Shanghai and MOE-LSC, Shanghai Jiao Tong University, Shanghai 200240, China. Email: {\tt xuzl@sjtu.edu.cn}. Z. Xu acknowledges the support
from the NSFC (grant No. 12325113) and the HPC center of Shanghai Jiao Tong University.}
}

\date{}
\maketitle

\begin{abstract}
  We introduce the \textit{Transformed} Generative Pre-Trained Physics-Informed Neural Networks (TGPT-PINN) for accomplishing nonlinear model order reduction (MOR) of transport-dominated partial differential equations in an MOR-integrating PINNs framework. Building on the recent development of the GPT-PINN that is a network-of-networks design achieving snapshot-based model reduction, we design and test a novel paradigm for nonlinear model reduction that can effectively tackle problems with parameter-dependent discontinuities. Through incorporation of a shock-capturing loss function component as well as a parameter-dependent transform layer, the  TGPT-PINN overcomes the limitations of linear model reduction in the transport-dominated regime. We demonstrate this new capability for nonlinear model reduction in the PINNs framework by several nontrivial parametric partial differential equations.

\end{abstract}

\section{Introduction}

Reduced order models (ROMs) are mainstays in computational science, playing an integral role in design optimization, uncertainty quantification, and in the construction of digital twins. ROMs can accelerate the evaluation of one-shot computational models in time-dependent settings, and are also capable of substantial acceleration in multi-query modeling for, e.g., parametric problems. We investigate ROMs for partial differential equation (PDE) models in this paper.  Although nonlinear model reduction has gained attention in recent years, most of the established ROM theory and algorithms focus on linear model reduction techniques that are known to be effective for diffusion-dominated PDE problems. However, such linear reduction approaches frequently fail when applied to PDEs that involve transport-dominated phenomena. This failure arises from a fundamental mathematical limitation: the slow decay of the Kolmogorov $n$ width \cite{pinkus_n-widths_1985} for transport problems with discontinuities \cite{OhlbergerRave2016}. Many recent advances in ROMs have therefore sought nonlinear reduction strategies that are not bound by the limitations of linear reduction.

Recent years have witnessed a wealth of nonlinear reduction methods being developed. The first set of approaches involve transforms of the trial basis \cite{Welper2017, cagniart2019model, NairBalajewicz2019, KrahSrokaReiss2021, Taddei2020, ReissSchulzeSesterhenn2018, RimPeherstorferMandli2023}. While these transforms improve the approximation property of the surrogate space, they tend to rely on additional knowledge about the problem and the transforms are typically not computable in an automatic way. 
The second set adopts neural networks as a component of its algorithm 
\cite{LeeCarlberg2020,KimChoiWidemannZohdi2020,KimChoiWidemannZohdi2022,dahmen2022compositional,PapapiccoDemoGirfoglioStabileRozza2022}. While being more general, these methods tend to be purely data driven, attempting to capture low-dimensional structure in a latent space representation, with various degrees of success.
The third is the optimal transport, specifically the Wasserstein metric, based approaches 
\cite{IolloLombardi2014, ehrlacher2020nonlinear, BattistiBlickhanEncheryEhrlacherLombardiMula2023}. While still nascent and more limited in applicability, it provides a unique perspective grounded in a mathematically sound theory. Another class of approaches use adaptive techniques where linear subspaces are progressively updated, for example in time-dependent problems \cite{PeherstorferWillcox2015,Peherstorfer2020,RimPeherstorferMandli2023}. Some existing work also use a specific types of nonlinear ans\"{a}tze (such as judicious combinations of polynomial-type functions and neural network functions) to inject nonlinear dynamics into a reduced model \cite{GeelenWrightWillcox2023,BarnettFarhatMaday2023}.

This paper makes advances that sit at the intersection of physics-informed machine learning and nonlinear reduced order modeling. We introduce the Transformed Generative Pre-Trained Physics-Informed Neural Networks (TGPT-PINN), a framework that extends Physics-Informed Neural Networks (PINNs) and reduced basis methods (RBM) to the nonlinear model reduction regime while maintaining the type of network structure and the unsupervised nature of its learning. By introducing a novel transform layer in a GPT-PINN, we enable the resolution of parameter-dependent discontinuity locations. Through the use of PINNs in snapshot-based model reduction, we eliminate the need for intrusive analysis and manipulation of existing PDE solvers; in particular, the ``online'' phase of the multi-query contexts leverages a network of pre-trained PINN-type networks, and hence all that must be explicitly manipulated is a PDE residual term in strong form and a transform layer capturing the discontinuity. Moreover, the determination of the network parameters in the transform layer are automatically achieved via back-propagation simultaneously with weights of the pre-trained PINNs serving as the super neurons of the reduced network.

We demonstrate the efficacy of the TGPT-PINN on several nontrivial examples. Its efficacy on explicit parametric functions is first investigated to reveal the expressive strength of the procedures for nonlinear model reduction. Next, the TGPT-PINN is applied to solutions of parametric PDEs, where explicit solutions are not available and are computed instead through a PINN-type formulation. We observe in this case as well that the TGPT-PINN is very effective at accomplishing model reduction, with very little \textit{a priori} knowledge about a PDE, or even whether or not the solution has discontinuities. 

As a novel effort in developing nonlinear reduction strategies, the main advance of our proposed method is an extension of the newly developed GPT-PINN \cite{chen2024gpt} to the nonlinear reduction regime by integrating shifts and transforms into the neural network/GPT-PINN framework. 
This is achieved by adding a \textit{transform layer} to the developed GPT-PINN. 
The benefit of this integration is two-fold. First, training inherits the greedy alrogithm-based model reduction capability of the GPT-PINN automatically. Second, the addition of the transform layer allows simultaneous training of the shift-and-transform network parameters of the transform layer and the mode coefficients of the GPT-PINN layer.
A particularly important practical improvement is that the involved neural network is also hyper-reduced, i.e. the architecture has orders of magnitude fewer trainable parameters than those in \cite{LeeCarlberg2020,KimChoiWidemannZohdi2020,KimChoiWidemannZohdi2022,dahmen2022compositional}. 

The newly proposed approach is tested in two steps. First, to isolate any influence due to choice of a numerical (``truth'') solver, we test the case when the parametric function set is analytically given. We show that, for the more challenging case when the Kolmogorov width decays slowly due to a moving kink or discontinuity, the TGPT-PINN is able to capture the parameter dependence exactly by one or a handful of neurons (while linear approaches, such as the Empirical Interpolation Method (EIM) \cite{Barrault_Nguyen_Maday_Patera}, fail to be effective). 
Next, we investigate when $u(x,\mu)$ is the numerical solution to a parametric PDE, which involves usage of a numerical solver. We show that, on examples when linear reduction does not work, the TGPT-PINN \textit{achieves machine accuracy with one neuron} (e.g., the transport equation with a discontinuous initial condition whose speed is a parameter). On examples when linear reduction does work well, the TGPT-PINN works by employing only one (nonlinear reaction) or three (nonlinear reaction diffusion) neurons while achieving better accuracy than the GPT-PINN with $10$ neurons. This significant improvement is a manifestation of the novel design of the TGPT-PINN, i.e., the addition of a trainable transform layer to a network of pre-trained networks.

The rest of this paper is organized as follows. Background materials including a motivating example, PINN and GPT-PINN are briefly presented in Section \ref{sec:bgd}. The main algorithm is given in Section \ref{sec:tgpt-pinn} with numerical results on three classess of analytically-given functions and three types of equations detailed in Section \ref{sec:numerics}. We draw conclusions in Section \ref{sec:conclusion}.

\section{Background}

\label{sec:bgd}

This section is devoted to a brief introduction of the two ingredients leading to the TGPT-PINN, namely snapshot transformations and the GPT-PINN for neural network based model order reduction for parametric PDEs.

\subsection{The Kolmogorov $n$-width and barrier}

Let $U$ be a Banach space, and let $K$ be a set (or ``manifold'') of functions that we wish to approximate using elements of an $n$-dimensional subspace. The theoretically optimal performance of such \textit{linear} reduced order models is quantified by the Kolmogorov $n$-width \cite{pinkus_n-widths_1985}, which measures how well elements in $K \subset U$ can be best approximated from a dimension-$n$ linear subspace $U_n \subset U$,
\begin{equation*}
  d_n(K) = \inf_{\dim U_n = n} \sup_{u \in K} \inf_{\phi \in U_n} \|u - \phi\|_U.
\end{equation*}
The class $K$ can be a set of functions that are encoded with a parameterization; in our setting a more salient example is to consider $K$ as the solution manifold of a (parametric) partial differential equation, $K = \{u_\mu \mid \mu \in \mathcal{D}\} \subset U$, where $\mu \mapsto u_{\mu}$ is the $\mu$-parameterized PDE solution operator, and $\mathcal{D}$ is frequently a subset of $\R^k$ for some $k \in \N$. 

The $n$-width $d_n(K)$ provides an accuracy lower bound for all reduced order models that rely on linear reduction, i.e., by building approximations from a(ny) dimension-$n$ subspace $U_n$. For example, the $n$-width is the best possible error when building $U_n$ from snapshots or proper orthogonal decomposition (POD) modes. 
For many parametric PDEs that are elliptic or parabolic in nature, $d_n(K)$ is known to decay exponentially in $n$ \cite{CohenDeVore2015}. For such cases, constructive (i.e., computationally feasible) greedy methods have been proven to generate reduced spaces that asymptotically match the Kolmogorov $n$-width 
\cite{BinevCohenDahmenDeVorePetrovaWojtaszczyk,BuffaMadayPateraPrudhommeTurinici2011}. This fast decay of the $n$-width along with the development of effective greedy algorithms are the driving forces behind linear model reduction and its descendants, such as RBM and GPT-PINN, respectively.

However, when facing solution manifolds generated by transport-dominated problems, the Kolmogorov $n$-width decreases much slower. Ohlberger and Rave \cite{OhlbergerRave2016} proved that $d_n(K) \sim n^{-1/2}$ where $K$ is the solution manifold of a (relatively simple) transport dominated problem. This extremely slow decay for transport-dominated problems is the fundamental reason that linear model order reduction methods are ineffective in this case, and motivates the need to develop and exercise \textit{nonlinear} model reduction approaches.

\subsection{A motivating example for transformed snapshots}
\label{sec:example}
\begin{figure}[htbp]
	\centering
 \includegraphics[width=0.42\textwidth]{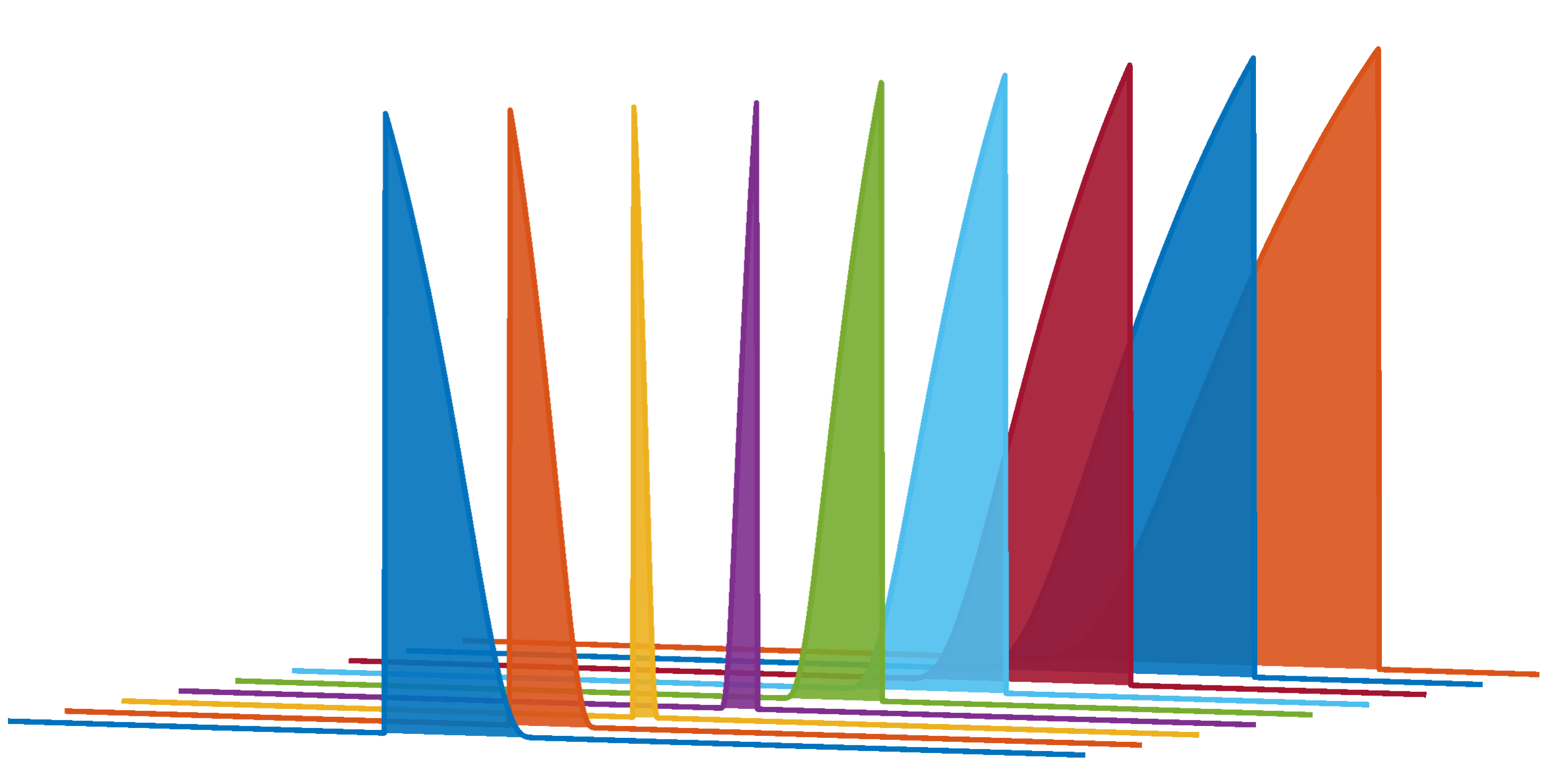}
\includegraphics[width=0.48\textwidth]{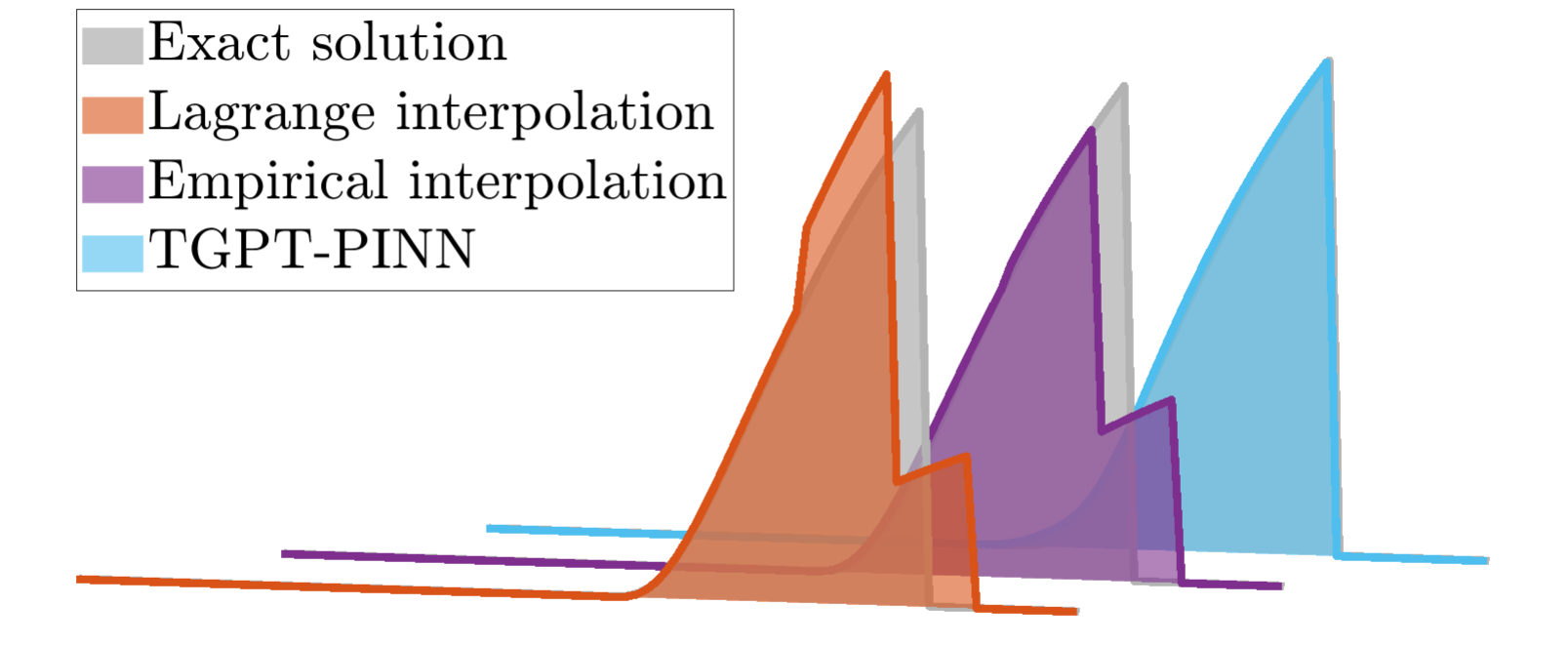}
\caption{Shown on the left are $9$ snapshots of $u(x,\mu)$ with various $\mu$ values. On the right are the function at $\mu = 1.0$ and its three kinds of interpolations,  the polynomial interpolation from snapshots with $\mu= 0.5, 0.85$ and $1.2$, the EIM approximation, and the TGPT-PINN interpolation.}
\label{fig:movingjump_eim}
\end{figure}

Consider the following example from \cite{Welper2017},
\begin{align*}
u(x, \mu):=\psi\left(\frac{x}{0.4+\mu}-1\right), \quad \psi(x):= \begin{cases}\exp \left(-\frac{1}{1-x^2}\right) & -1 \leq x<-\frac{1}{2}, \\ 0 & \text { else}.\end{cases}
\end{align*}
The function $\psi(x)$ is a discontinuous function at $x=-1/2$, which makes $u(x,\mu)$ discontinuous at $x=(\mu+0.4)/2$. We are therefore in a situation that the location of the discontinuity depends on the parameter. Interpolatory type techniques in the parameter $\mu$, including classical polynomial interpolations and the EIM which is a standard tool in model order reduction for the nonlinear and nonaffine setting, have limited effectiveness as demonstrated in Figure \ref{fig:movingjump_eim}. However, if we change the interpolation ans\"{a}tze 
\begin{equation}
    u(x,\mu) \approx \sum_{i=0}^n \ell_i(\mu) u(x, \mu^i) \quad \longrightarrow \quad u(x,\mu) \approx \sum_{i=0}^n \ell_i(\mu) u(\phi(\mu, \mu^i)(x), \mu^i),
    \label{eq:nonmoransatz}
\end{equation}
where $\ell_i$ are any basis (e.g., cardinal Lagrange interpolants associated to the $\mu^i$), we can design $\phi$ to eliminate the staircasing behavior. Here, $\phi(\mu, \eta)$ represents a transform from the ``source'' parameter $\eta$ to the ``target'' parameter $\mu$ so that the jump locations match. For this example, the choice
\begin{align*}
  \phi(\mu,\eta)(x) = \frac{0.4 + \eta}{0.4 + \mu} x,
\end{align*}
would enable exact reconstruction with $n=0$, i.e., $u(\phi(\mu, \eta)(x), \eta) = u(x, \mu)$. This shows that sufficiently complicated transformation maps $\phi$ can restore expressive power in reduced order models. However, we make a further observation that the even simpler map,
\[
\phi(\mu, \eta)(x) = x + \frac{\eta - \mu}{2},
\]
which is linear in $(\eta, \mu)$, is enough to ensure that the discontinuity locations are exactly captured, and hence even relatively simple maps $\phi$ can provide enormous benefit. The proposed TGPT-PINN procedure attempts to identify the map $\phi$ automatically. As a result, the TGPT-PINN approximation of this function is indistinguishable from the exact solution while the standard approaches produce noticeable staircases leading to large interpolation errors, see Figure \ref{fig:movingjump_eim}.

\subsection{PINN and GPT-PINN}
\begin{figure}[htbp]
    \centering
    \includegraphics[width=0.99\textwidth]{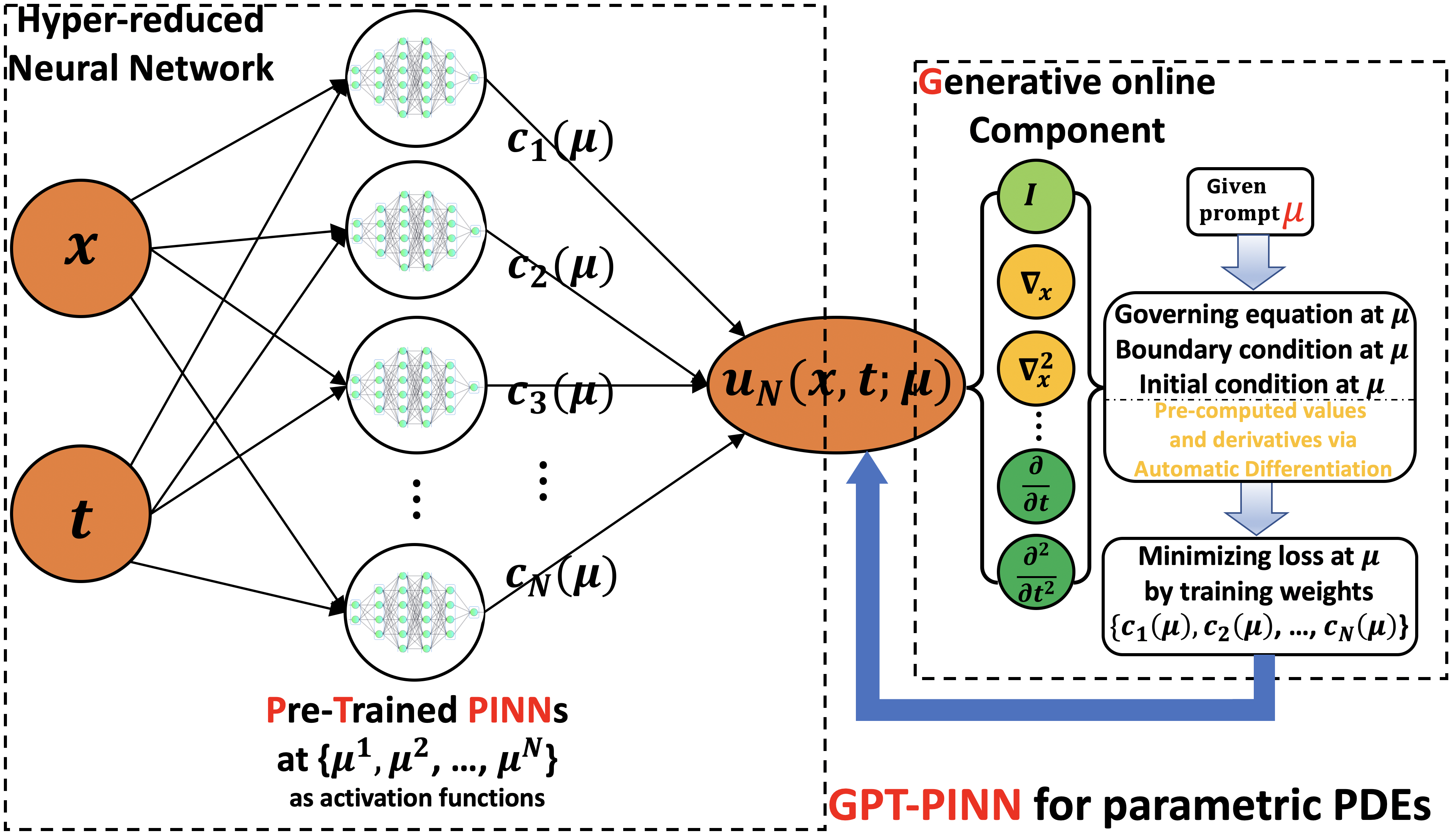}
   \caption{The GPT-PINN architecture \cite{chen2024gpt}. A hyper-reduced network adaptively embeds pre-trained PINNs at the nodes of its sole hidden layer. It then allows a quick online generation of a surrogate solution at any given parameter value.}
    \label{fig:gptpinn_diagram}
\end{figure}
Physics-informed neural networks (PINNs), popularized by Raissi {\it et al. }\cite{raissi2019physics}, have emerged as an increasingly popular alternative to traditional numerical methods for PDEs in recent years. They adopt Deep Neural Networks (DNNs) to approximate the solutions of PDEs while incorporating a strong physics-based PDE prior encoded into the loss function to constrain the output of the DNN. 
In comparison to traditional numerical solvers, the advantages of PINNs include that they are able to numerically solve the PDE without discretizing the spatiotemporal domain, and that they can leverage automatic differentiation \cite{baydin2017automatic,Paszke2017AutomaticDI} to algorithmically minimize the residual in the loss. 
However, PINNs have some weaknesses. For example, training a vanilla PINN is usually significantly slower than employing a classic numerical method to solve the corresponding PDE. To address this shortcoming, the Generative Pre-Trained PINNs (GPT-PINNs) were developed \cite{chen2024gpt} as a meta-learning approach for parametric systems to 
reduce the architecture and corresponding number of trained parameters in a PINN, thereby reducing the computational time required to solve parametric PDEs with PINNs. 
The GPT-PINN requires an initial (``offline'') investment cost dedicated to learning the parametric dependence during the offline stage, which is guided by a mathematically reliable greedy algorithm. With this initial investment, The GPT-PINN is capable of providing significant computational savings in the multi-query and real-time settings thanks to the fact that their marginal cost is of orders of magnitude lower than that of an individual PINN solve \cite{chen2024gpt}. 

The GPT-PINN architecture, depicted in Figure \ref{fig:gptpinn_diagram}, is a network-of-networks, where activation functions in the  sole hidden layer of the outer network are chosen in a customized way. The inner networks, defining the activation functions, are full pre-trained PINNs instantiated by the PDE solutions at a  set of  adaptively-selected parameter values $\{\mu^1, \mu^2, \cdots, \mu^n\}$ chosen by a greedy algorithm. The corresponding outer-/meta-network is hyper-reduced in comparison to the inner networks, having only one hidden layer. 
Moreover, this meta layer adaptively ``learns'' the parametric dependence of the system and can ``grow'' its hidden layer one super neuron/network at a time. The GPT-PINN is capable of generating approximate solutions for the parametric system across the entire parameter domain accurately and efficiently, with a cost independent of the size of the full PINN.

\section{The TGPT-PINN algorithm}

\label{sec:tgpt-pinn}

\begin{figure}[thbp]
\centering
\includegraphics[width=\textwidth]{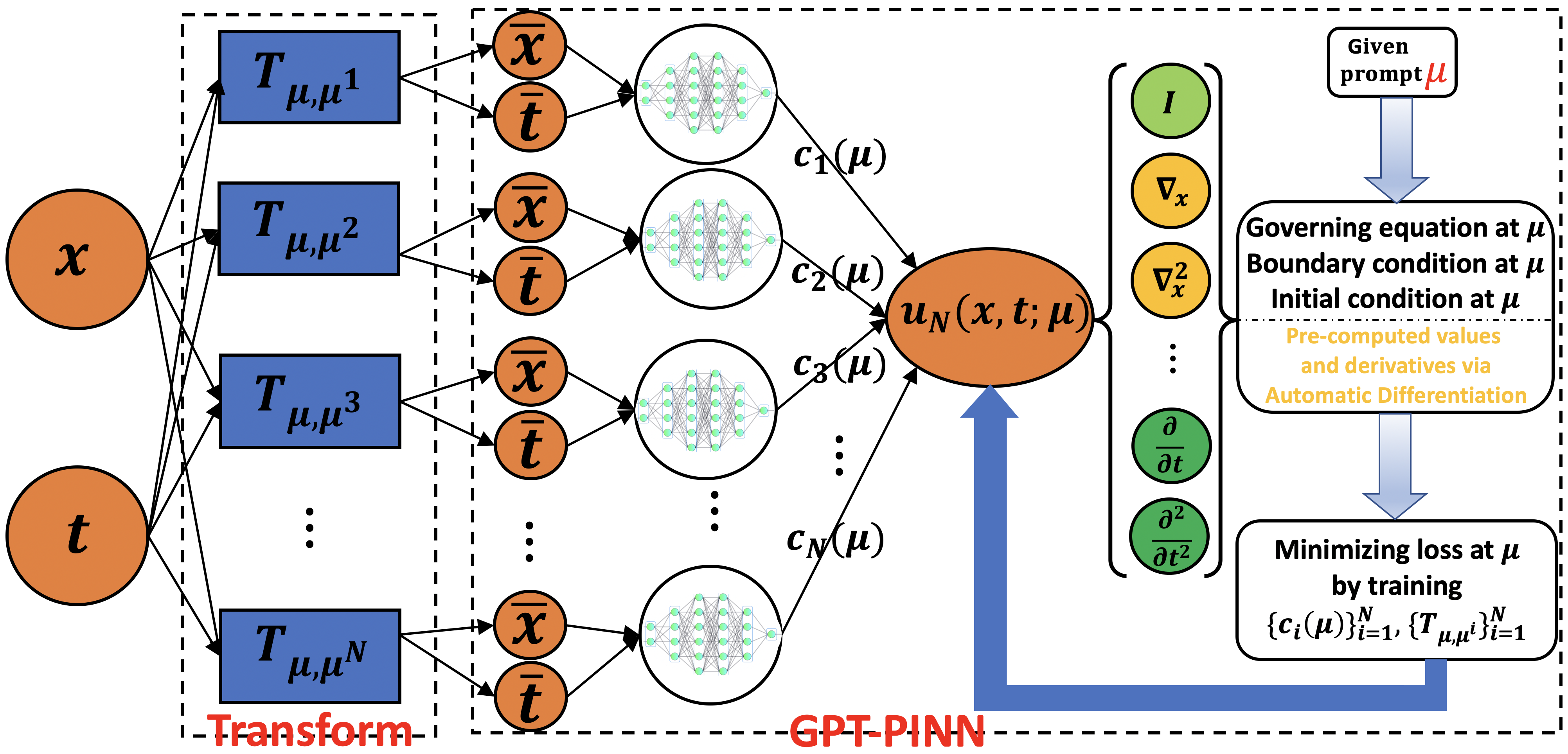}
  \caption{The TGPT-PINN design schematic. For any given parameter value $\mu$, a $\mu$-dependent loss is constructed and the coefficients $c_j(\mu)$ and the weights and biases in $T_{\mu,\mu^i}$ are trained.}
\label{fig:tgpt-pinn}
\end{figure}

The insight provided by Section \ref{sec:example} is that, to achieve nonlinear model order reduction, one way is to compose the snapshot with a \textit{parameter-dependent transformation}. This is the basic idea of the TGPT-PINN whose schematic design is shown in Figure \ref{fig:tgpt-pinn}. 
It seeks to approximate $u(\bm{x},t;\mu$) as follows,
\begin{equation}
u(\bm{x},t;\mu) \approx \sum_{i=1}^N c_i(\mu) u(T_{\mu,\mu^i}(\bm{x},t); \mu^i).
\label{eq:tgpt_ansatz}
\end{equation}
All the quantities, $N$, $c_i(\mu)$, $\mu^i$, $T_{\mu,\mu^i}$ and $u(\bm{x},t; \mu^i)$, are obtained through training in a two-step procedure. The first step is an ``offline'' step that trains $\mu^i$, $N$ and $u(\bm{x},t; \mu^i)$; this step can be expensive as we require well-resolved solutions $u$ to PDEs and we must sweep over the parameter domain $\mathcal{D} \subset \R^k$ to identify parameter values $\mu^i$. The second step is an ``online'' step that, given an arbitrary $\mu$, computes $T_{\mu,\mu^i}$ and $c_i(\mu)$ for $i\in[N]=\{1,\cdots, N\}$. This second step is much more computationally efficient, training a neural network with only $N (d^2 + 3 d + 3)$ degrees of freedom, i.e., linear in the number of snapsots $N$ and \textit{substantially} fewer degrees of freedom than required for even a single snapshot $u$. More explicitly, the TGPT-PINN is the following two-step procedure:
\begin{align*}
  \textrm{$i=1$, $\mu^1$, ``Offline'':} & \hskip 20pt \textrm{PDE/residual} \xrightarrow{\textrm{Minimize } \mathcal{L}_{\mathrm{PINN}}} \hskip 5pt u(\cdot,\cdot;\mu^i) \hskip 5pt \xrightarrow{\textrm{Sufficient accuracy}} \hskip 5pt N \coloneqq i \\
  &\hskip 40pt  \uparrow \hskip 120pt \downarrow \textrm{{\scriptsize Insufficient accuracy}} \\
  &\hskip 40pt \mu^i \hskip 10pt \xleftarrow[\textrm{Search parameter space}]{} \hskip 10pt i \mathrel{+}= 1 \\\\
  \textrm{``Online'':} & \hskip 20pt \mu \xrightarrow{\textrm{Minimize } \mathcal{L}^{\mathrm{TGPT}}_{\mathrm{PINN}}} c_i(\mu),\;\; T_{\mu,\mu^i}(\cdot,\cdot), \hskip 5pt i \in [N] \xrightarrow{\eqref{eq:tgpt_ansatz}} u(\cdot,\cdot;\mu)
\end{align*}
The ``online'' step is visualized in Figure \ref{fig:tgpt-pinn}. The main strength of this approach is that function composition is natural and straightforward in neural networks - simply adding a layer or block, and the network parameters defining the transformation can be trained together with the mode coefficients $\{c_i(\mu)\}_{i=1}^N$ in the GPT-PINN block of the network. The loss functions $\mathcal{L}_{\mathrm{PINN}}$ and $\mathcal{L}^{\mathrm{TGPT}}_{\mathrm{PINN}}$ are described in the coming sections. All minimization steps are performed using somewhat standard neural network optimization and back-propagation procedures.

\subsection{PDE formulation}
We define our problem to be the following time-dependent PDE  on the spatial domain $\Omega \subset \mathbb{R}^d$ with boundary $\partial \Omega$, and parametric domain $\mu \in \mathcal{D}$:
\begin{equation}
\begin{aligned}
  & \frac{\partial}{\partial t} u(\bm{x}, t;\mu)+\mathcal{F}[u(\bm{x}, t;\mu)]=0, \quad \bm{x} \in \Omega, \quad t \in[0, T], \quad \mu \in \mathcal{D}\\
  & \mathcal{G}(u)(\bm{x}, t;\mu)=0, \quad \bm{x} \in \partial \Omega, \quad t \in[0, T], \quad \mu \in \mathcal{D}\\
  & u(\bm{x}, 0;\mu)=u_0(\bm{x};\mu), \quad \bm{x} \in \Omega, \quad \mu \in \mathcal{D}.
\end{aligned}
\label{eq:pde}
\end{equation}
Here $\mathcal{F}$ is a differential operator and $\mathcal{G}$ denotes a boundary operator \cite{chen2024gpt}. Our goal is to compute the numerical solution to this PDE, which we shall accomplish using the TGPT-PINN methodology.

The remainder of this section is devoted to describing the particulars of the TGPT-PINN. In Section \ref{sec:shock_apd} we will introduce an approach that modifies the ``offline'' portion of the TGPT-PINN to enhance its ability to resolve discontinuities in the solution for transport-dominated problems.

\subsection{Offline: The PINN solver}
\label{sec:shock_apd}
We discuss the ``offline'' procedure of the TGPT-PINN that computes $u(\cdot,\cdot;\mu)$ for a given, fixed $\mu$. In particular, we will use $\mu = \mu^i$ in $u$ for the TGPT-PINN. Because $\mu$ is fixed in this section, we frequently omit writing the $\mu$-dependence in what follows. The idealized loss function that we use corresponds to a standard PINN loss function,
\begin{equation}
  \mathcal{L}(u)=\mathcal{L}_{\rm int}(u) + \varepsilon_i \int_{\Omega} \left\|u(\bm{x}, 0)-u_0(\bm{x})\right\|_2^2 d \bm{x}+\varepsilon_b \int_{\partial \Omega \times [0,T]}\|\mathcal{G}(u)(\bm{x}, t)\|_2^2 d \bm{x}\, d t,
\label{eq:continousloss}
\end{equation}
which is constructed by summing losses corresponding to the residual of the PDE, the initial values and the boundary conditions, and the domain-interior loss is given by,
\[
\mathcal{L}_{\rm int}(u) \coloneqq 
 \int_{\Omega \times (0,T]}\left\|\lambda(\bm{x},t)\,\left(\frac{\partial}{\partial t} u(\bm{x}, t)+\mathcal{F}(u)(\bm{x}, t)\right)\right\|_2^2 d \bm{x}\, d t.
\]
We choose the weighting constants as $\varepsilon_i=1$ and $\varepsilon_b=1$.  In practice, PINNs-type losses can suffer from the exploding or vanishing gradient problem when solutions are discontinuous. The factor $\lambda$ above is introduced to ameliorate this effect, and is designed to attenuate the loss when the gradient of the solution is large. We therefore introduce the following ``shock-capturing weighting'' factor in the calculation of the residual of the PDE, which was originally proposed in Liu {\it et al.} \cite{liu2024discontinuity},
\begin{equation*}
\lambda(\bm{x},t)=\frac{1}{\varepsilon_\lambda\,|\nabla \cdot u|+1}.
\end{equation*}
We choose $\varepsilon_\lambda=0.1$. The PINNs approach chooses a discretization for $u$ and minimizes the loss function through optimization. Like a standard PINN, we use a fully connected neural network $\psi^{\mu}_{\mathsf{NN}}(\bm{x},t)$ to approximate $u$, where the weights and biases of $\psi^{\mu}_{\mathsf{NN}}$ are optimized. The loss $\mathcal{L}$ is approximated through a sampling/quadrature-type procedure. In particular, the actual discretized loss that is optimized is,
\begin{equation}
\begin{aligned}
\mathcal{L}_{\mathrm{PINN}}&\left(\psi_{\mathrm{\mathsf{NN}}}^\nu\right)=  \frac{1}{\left|\mathcal{C}_o\right|} \sum_{(\bm{x}, t) \in \mathcal{C}_o}\left\|\lambda\cdot\left(\frac{\partial}{\partial t}\left(\psi_{\mathrm{\mathsf{NN}}}^\nu\right)(\bm{x}, t)+\mathcal{F}\left(\psi_{\mathrm{\mathsf{NN}}}^\nu\right)(\bm{x}, t)\right)\right\|_2^2 \\
& +\varepsilon_b \cdot \frac{1}{\left|\mathcal{C}_{\partial}\right|} \sum_{(\bm{x}, t) \in \mathcal{C}_{\partial}}\left\|\mathcal{G}\left(\psi_{\mathrm{\mathsf{NN}}}^\nu\right)(\bm{x},t)\right\|_2^2 +\varepsilon_i \cdot \frac{1}{\left|\mathcal{C}_i\right|} \sum_{\bm{x} \in \mathcal{C}_i}\left\|\psi_{\mathrm{\mathsf{NN}}}^\nu(\bm{x}, 0)-u_0(\bm{x})\right\|_2^2.
\end{aligned}
\label{eq:discreteloss}
\end{equation}
Above, we sample collocation/training points in certain fashion (described below) from the PDE spatiotemporal domain $\cC_o \subset \Omega \times (0,T)$, spatiotemporal boundary $\cC_\partial \subset \partial \Omega \times [0,T]$, and spatial interior $\cC_i \subset \Omega$, and use them to form an approximation of the continuous loss \eqref{eq:continousloss}. In our examples, the boundary operator $\mathcal{G}$ is taken to be a spatially periodic one for simplicity.

The training points must be carefully chosen in the presence of discontinuities. In general it is well known that the distribution of training points can impact the network performance. In particular, dense concentration of training points near the spatiotemporal vicinity of a discontinuity can effectively enhance the network performance in approximating the PDE solutions with jumps. For the transport equation, the location of the discontinuity is often determined by the parameters. Therefore, to achieve better network performance, in practice it is necessary for the training/sampling points to be constructed in a parameter-dependent fashion. 

Finally, we remark that the above construction optimizes for ``meta-neurons'', i.e. the networks within the GPT-PINN. For some of our test cases, in particular those in Section \ref{ssec:example-function}, we use analytically available formulas for $u(\cdot,\cdot;\mu^i)$, and in those examples the entire PINN apparatus is simply replaced with this explicit function evaluation.

\subsection{TGPT-PINN}

As shown by Figure \ref{fig:tgpt-pinn}, the hidden layers of the TGPT-PINN include initially a transform layer and, similar to GPT-PINN \cite{chen2024gpt}, are followed by a meta-PINN layer with pre-trained PINNs as activation functions. The novel contribution is the transform layer that allows the TGPT-PINN to achieve nonlinear model order reduction. This transform layer allows us to further explore the potential of the full PINNs used as the activation function, thus significantly enhancing the expressivity of the (pre-trained) full PINNs, even with only one pre-trained PINN in a single neuron in the hidden layer. This boost of expressivity from the parameter-dependent transform layer allow us to construct a hyper-reduced neural network of pre-trained full networks featuring parameter-dependent discontinuous functions. In particular, we achieve sizes of the outer network that are much smaller than what the Kolmogorov width of the problem implies is optimal for linear reduction of transport-dominated problems. 

\subsubsection{Design of the transform layer}

The \textit{source-to-target transform layer} $T_{\mu,\eta}$ in \eqref{eq:tgpt_ansatz} of the proposed TGPT-PINN is abstractly,
\[
T_{\mu,\mu^i}(\bm{x},t): \Omega \times [0, T] \longrightarrow \Omega \times [0, T],
\]
and in principle it should be surjective. In this paper, we set it as
\begin{equation}
T_{\mu, \eta}(\bm{x},t) \coloneqq 
{\rm Mod}_{\Omega, T}
\left(
W_{\mu,\eta} 
\left(
\begin{tabular}{c}
     $\bm{x}$\\
     $t$ 
\end{tabular}
\right)
+ b_{\mu,\eta}
\right)
, \quad \eta = \mu^1, \dots, \mu^N.
\label{eq:transform}
\end{equation}
Here, $W_{\mu,\eta} \in \R^{(d+1) \times (d+1)}$, and $b_{\mu,\eta} \in \R^{d+1}$ making $T_{\mu, \eta}$ a simple linear transformation that depends on $\mu$ and $\eta$. In particular, the dependence on $\mu$ is enforced through ``online'' training. Moreover, to make sure that the range of $T_{\mu, \eta}$ is exactly $\Omega \times [0,T]$, we apply ${\rm Mod}_{\Omega, T}(\cdot)$ which is an element-wise modulo map, ensuring that each component of $T$ outputs on the appropriate slice of $\Omega \times [0,T]$. Our notation of $(\mu,\eta)$ subscripts on $W_{\mu,\eta}$ is meant to reveal that $W$ depends on both $\eta$ (i.e., indexing the particular snapshot) as well as $\mu$ (i.e., the given value of the parameter in the online phase affects the value of the weights trained through minimization).

We use the particular form of $T$ above for all the experiments in this paper. However, the general framework of the TGPT-PINN is not limited by the relatively simple choice of $T$ above. For example, $T$ could itself be a deep neural network. Our simple choice is motivted by the discussion in Section \ref{sec:example}, demonstrating that even simple linear-type maps can be effective for transport-based problems. In particular, $T$ can be trained to stretch or shrink the $(\bm{x},t)$ variables, and in particular can be trained to expand spectral content in the frequency domain. For this reason, as shown by our numerical examples, the proposed architecture also works well for much more complicated problems that don't just feature transport with constant speed.

\subsubsection{Online phase training}

The transformation \eqref{eq:transform} and the TGPT-PINN ans\"{a}tze \eqref{eq:tgpt_ansatz} mean that a given TGPT-PINN with $n$ PINNs (pre-trained at $\{\mu^1, \dots, \mu^n\}$) has the following $n (d^2+3d+3)$ network parameters to train
\begin{equation}
    \Theta(\mu) \coloneqq \left\{ \{W_{\mu, \mu^i}\}_{i=1}^n, \{b_{\mu, \mu^i}\}_{i=1}^n, \{c_i(\mu)\}_{i=1}^n \right\}.
    \label{eq:tgpt_paras}
\end{equation}
Similar to the GPT-PINN, this number of network parameters is \textit{independent} of the architecture parameters used to train the individual PINNs, and depends strictly linearly on $n$, the number of snapshots. With the network version of the TGPT-PINN ans\"{a}tze \eqref{eq:tgpt_ansatz} denoted by
\begin{equation}
\Psi_{\mathrm{\mathsf{NN}}}^{\Theta(\mu)}(x, t) \coloneqq \sum_{i=1}^n c_i(\mu) \psi_{\mathrm{\mathsf{NN}}}^{\mu^i}(T_{\mu,\mu^i}(x,t)),
\label{eq:tgpt_ansatz_pinn}
\end{equation}
the loss function of the TGPT-PINN is set to be the same loss as the full PINN, consisting of three parts: the residual of the PDE and the losses corresponding to the initial value condition and the boundary condition: 
\begin{equation}\label{eq:loss}
\begin{aligned}
\mathcal{L}_{\mathrm{PINN}}^{\mathrm{TGPT}}&(\Theta(\mu))=\frac{1}{\left|\mathcal{C}_o^r\right|} \sum_{(\bm{x}, t) \in \mathcal{C}_o}\left\|\lambda \cdot\left(\frac{\partial}{\partial t}\left( \Psi_{\mathrm{\mathsf{NN}}}^{\Theta(\mu)}\right)(\bm{x}, t)+\mathcal{F}\left( \Psi_{\mathrm{\mathsf{NN}}}^{\Theta(\mu)}\right)(\bm{x}, t)\right)\right\|_2^2 \\
&+\varepsilon_b \cdot \frac{1}{\mid \mathcal{C}_{\partial}^r \mid} \sum_{(\bm{x}, t) \in \mathcal{C}_{\partial}}\left\|\mathcal{G}\left( \Psi_{\mathrm{\mathsf{NN}}}^{\Theta(\mu)}\right)(\bm{x},t)\right\|_2^2+\varepsilon_i \cdot \frac{1}{\left|\mathcal{C}_i^r\right|} \sum_{\bm{x} \in \mathcal{C}_i}\left\| \Psi_{\mathrm{\mathsf{NN}}}^{\Theta(\mu)}(\bm{x}, 0)-u_0(\bm{x})\right\|_2^2.
\end{aligned}
\end{equation}
The online collocation sets $\cC_o^r \subset \Omega \times [0,T]$, ${\cC_\partial^r \subset \partial \Omega \times [0, T]}$ and ${\cC_i^r \subset \Omega}$ need not be related to their full PINN counterparts $\cC_o, \cC_\partial$ and $\cC_i$, but in this paper we take them to be the same sets for simplicity. The training of $\Theta(\mu)$ is accomplished through standard automatic differentiation and back propagation. See Ref. \cite{chen2024gpt} for more details on this step, such as precomputations for fast training of the reduced network, and the non-intrusiveness of the (T)GPT-PINN.

\subsubsection{Offline training of the TGPT-PINN}

\begin{algorithm}[htbp]
    \caption{TGPT-PINN for parametric PDE: Offline stage}
    \label{alg:gptpinn_greedy}
    {\bf Input: }{A random (or given) $\mu^1$, training set $\Xi_{\rm train} \subset \calD$, full PINN.}
\begin{algorithmic}[1]
\State Train a full PINN at $\mu^1$ to obtain $\Psi_{\mathsf{NN}}^{\mu^1}$. Precompute quantities necessary for $\nabla_\Theta(\mu) \mathcal{L}_{\text{PINN}}^{\text{TGPT}}$ at collocation nodes $\cC_o^r$, $\cC_\partial^r$, and $\cC_i^r$. Set $n=2$.
\While{\textit{stopping criteria not met,}}
\State Train the ($n-1$)-neuron TGPT-PINN at $\mu$ for all $\mu \in \Xi_{\rm train}$ and record the indicator $\Delta_{\mathsf{NN}}^r(\Theta(\mu))$.
 \State Choose $\mu^n = \displaystyle
  \mbox{\rm arg}\hspace*{-1pt}\max_{\mu\in{\Xi_{\rm train}}}
 {\Delta_{\mathsf{NN}}^r(\mu)}$.
\State  Train a full PINN at $\mu^n$ to obtain $\Psi_{\mathsf{NN}}^{\mu^n}$. Precompute quantities necessary for $\nabla_\Theta \mathcal{L}_{\text{PINN}}^{\text{TGPT}}$ at collocation nodes $\cC_o^r$, $\cC_\partial^r$, and $\cC_i^r$.
\State Update the TGPT\_PINN by adding a neuron to the hidden TGPT-PINN layer to construct the $n$-neuron TGPT-PINN.
\State Set $n \leftarrow n+1$.
\EndWhile
\end{algorithmic} 
 {\bf{Output:}~}  $N$-neuron TGPT-PINN, with {$N$ being the terminal index.}
\end{algorithm}
With the online solver described above, the offline training amounts to the application of the greedy algorithm outlined in Algorithm \ref{alg:gptpinn_greedy}.  
The meta-network adaptively ``learns'' the parametric dependence of the system and ``grows'' the TGPT-PINN hidden layer and enriches the transform layer one neuron and transformation at a time. At every step, we select the parameter value that is worst approximated by the current meta-network. Specifically, we first randomly select, in a discretized parameter domain $\Xi_{\rm train} \subset \mathcal{D}$, one parameter value $\mu^1$ and train the associated (highly accurate) PINN $\Psi_{\mathsf{NN}}^{\mu^1}$. The algorithm then decides how to ``grow'' its meta-network by scanning the entire discrete parameter space $\Xi_{\rm train}$ and, for each parameter value, training this reduced network (of 1 neuron $\Psi_{\mathsf{NN}}^{\mu^1}$). As it scans, it records an error indicator $\Delta_{\mathsf{NN}}^r(\Theta(\mu))$ at every location. The next parameter value $\mu^2$ is the one generating the largest error indicator. The algorithm then proceeds by training a full PINN at $\mu^2$ and therefore grows its hidden PINN layer into two neurons with customized (but pre-trained) activation functions $\Psi_{\mathsf{NN}}^{\mu^1}$ and $\Psi_{\mathsf{NN}}^{\mu^2}$. This process is repeated until the stopping criteria is met which can be either that the error indicator is sufficiently small or a pre-selected size of the reduced network is met.

\section{Numerical results}

\label{sec:numerics}

In this section, we test the proposed TGPT-PINN and report  numerical results. This is done with two types of experiments. First, to separate the influence of the PINN solver which is analogous to the ``truth'' discretization in traditional RBM methods, we test the case when $u(\bm{x},\mu)$ is a known function with different types of regularity. Next, we treat $u(\bm{x},\mu)$ as a solution to a parameterized PDE, numerically computed using a PINN solver. The code for all these examples are published on GitHub at \href{https://github.com/DuktigYajie/TGPT-PINN}{https://github.com/DuktigYajie/TGPT-PINN}.

\subsection{TGPT-PINN vs EIM as function approximators}\label{ssec:example-function}
We first test the case when 
\[
u(\bm{x},\mu): \Omega \times \mathcal{D} \longrightarrow \mathbb{R}
\]
is analytically given, and compare the TGPT-PINN with the EIM \cite{Barrault_Nguyen_Maday_Patera} which is the standard approach in model order reduction for approximating nonlinear and nonaffine functions. Table \ref{tab:tgpt_functions} lists the results when both approaches are applied to $7$ functions of different regularities.

\begin{table}[!htb]
\centering
  \resizebox{\textwidth}{!}{
\begingroup
\renewcommand{\arraystretch}{1.5} 
\begin{tabular}{cccccccc}
  \rowcolor{white} & & & &\multicolumn{2}{c}{EIM} & \multicolumn{2}{c}{TGPT-PINN} \\
  \rowcolor{white} \multirow{-2}{*}{Function} & \multirow{-2}{*}{$\Omega$} & \multirow{-2}{*}{$\mathcal{D}$} & \multirow{-2}{*}{$\left|\Xi_{\mathrm{train}}\right|$} &$\#_{\text{basis}}$ &$L_2$ error&$\#_{\text{basis}}$&$L_2$ error\\\toprule
  $\sin(x+\mu)$&$[\pi,\pi]$&$[-5, 5]$& $201$ & 3&7.6e-15&1&1.2e-14\\
  \rowcolor{gray!30}
  $\sin(\mu x)$&$[-\pi,\pi]$&[1, 2]& $201$ & 17&1.8e-15&1&9.9e-13\\
  $\sin(\mu_1(x+\mu_2))$&$[-\pi,\pi]$&$[-5, 5]^2$& $20 \times 20 $ & 21&2.0e-14&1&9.9e-11\\
  \rowcolor{gray!30}
  & & & & 192 & 4.2e-04 & &\\
\rowcolor{gray!30}
  \multirow{-2}{*}{$\max(\sin(x+\mu),0)$} & \multirow{-2}{*}{$[-\pi,\pi]$} & \multirow{-2}{*}{$[-5, 5]$} & \multirow{-2}{*}{$201$} & 193 & 2.3e-15 & \multirow{-2}{*}{1} & \multirow{-2}{*}{9.7e-15} \\
  \multirow{2}{*}{$|x+\mu|$}&\multirow{2}{*}{$[-10,10]$}&\multirow{2}{*}{$[-5, 5]$}& \multirow{2}{*}{$201$} & 101&8.1e-02&\multirow{2}{*}{1}&\multirow{2}{*}{4.9e-15}\\
  &&&&102&5.6e-14&&\\
  \rowcolor{gray!30}
  &&&&100&1.9e-01&&\\
  \rowcolor{gray!30}
  \multirow{-2}{*}{$\psi\left(\frac{x}{0.4+\mu}-1\right)$}&\multirow{-2}{*}{[-1,1]}&\multirow{-2}{*}{$[-1,1]$}& \multirow{-2}{*}{$201$} & 101 &6.9e-16&\multirow{-2}{*}{10}&\multirow{-2}{*}{1.5e-05}\\
  $\frac{1}{\sqrt{(x-\mu_1)^2+(y-\mu_2)^2}}$&$[-1,1]^2$&$[-1,-0.01]^2$& $21 \times 21$ &150&3.8e-15&1&9.9e-13\\
\end{tabular}
\endgroup
  }
  \caption{Results for \Cref{ssec:example-function,sssec:function-example-1,sssec:function-example-2,sssec:function-example-3}: Comparison of the EIM and the proposed TGPT-PINN as function approximators. The parameter space $\mathcal{D}$ is discretized to $\Xi_{\mathrm{train}}$ using $|\Xi_{\mathrm{train}}|$ equispaced points.}
\label{tab:tgpt_functions}
\end{table}

The first three of $7$ functions in Table \ref{tab:tgpt_functions} are smooth functions. We observe that the TGPT-PINN is able to capture the parameter dependence exactly and, as a result, approximates each of these functions by one neuron to machine precision. Because the parametric dependence for these functions is smooth, the EIM is also able to approximate the functions to machine accuracy. However, it requires many more basis functions since it is a linear reduction approach (while the TGPT-PINN is nonlinear). We discuss the remaining four functions, which have different types of regularity, in the sections below.

\subsubsection{Functions with moving kinks}\label{sssec:function-example-1}

\begin{figure}[htbp]
\centering
\includegraphics[width=0.48\textwidth]{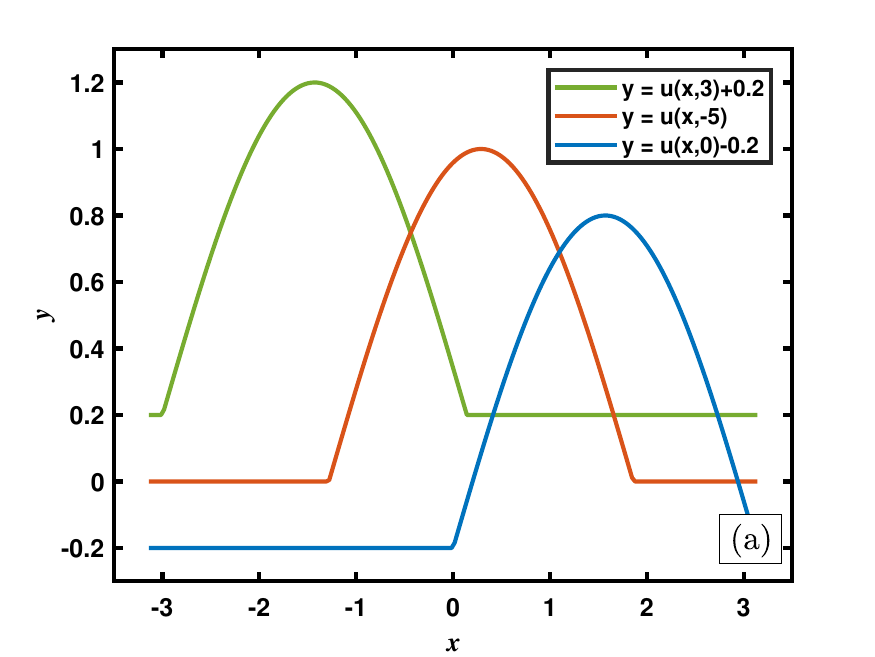}
\includegraphics[width=0.48\textwidth]{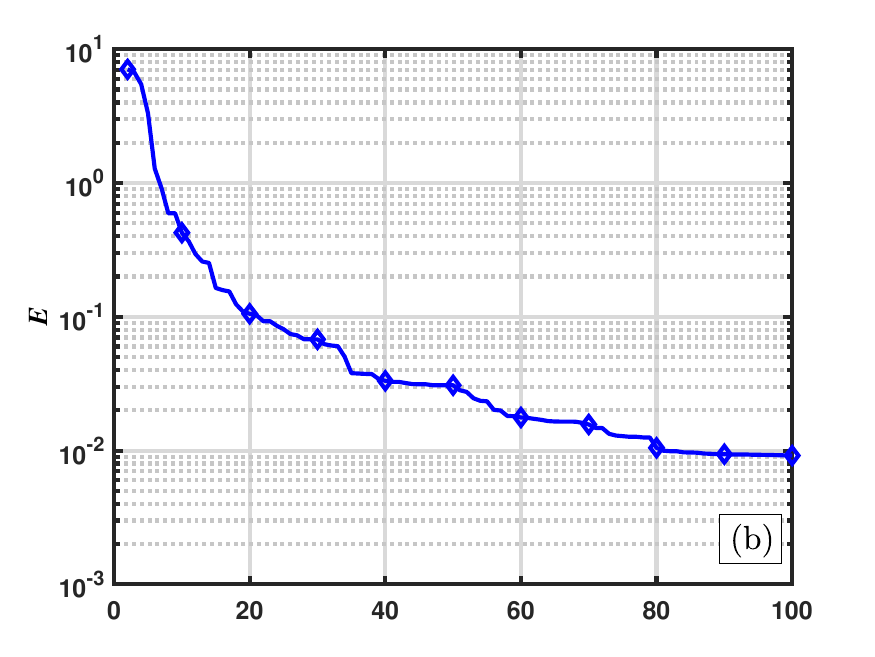}\\
\includegraphics[width=0.48\textwidth]{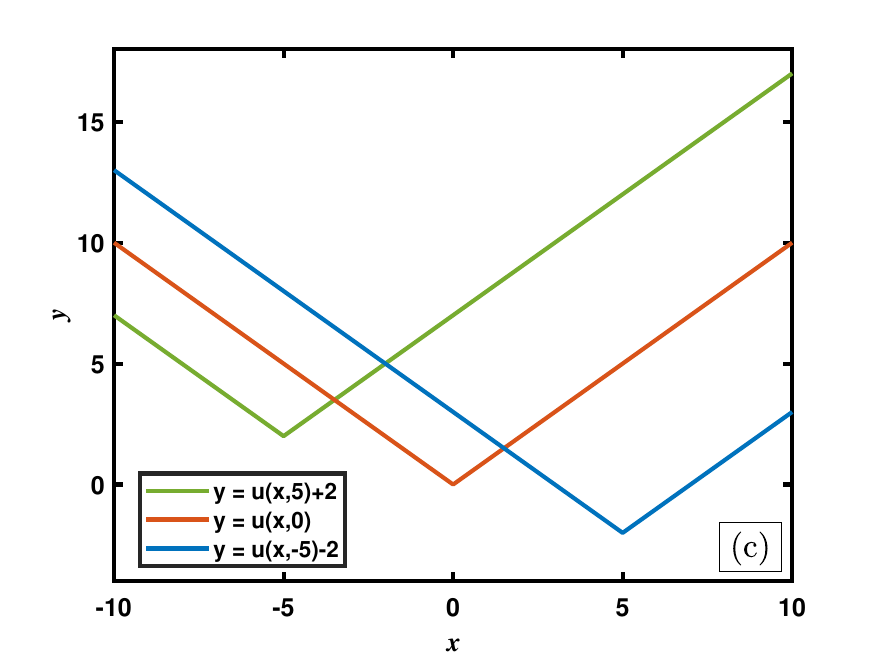}
\includegraphics[width=0.48\textwidth]{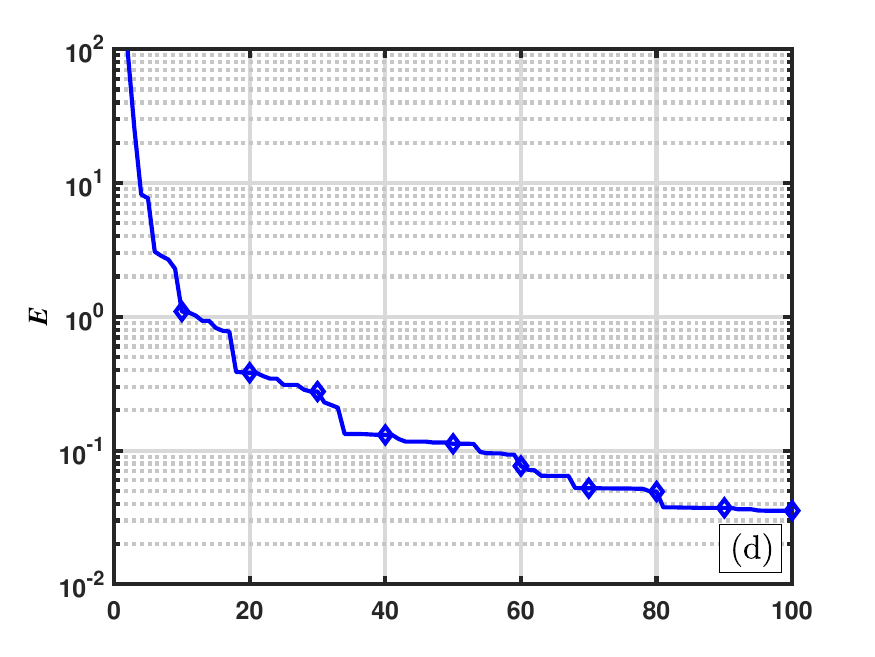}
  \caption{Results from Section \ref{sssec:function-example-1} paired with rows 4 and 5 of Table \ref{tab:tgpt_functions}: Snapshots and EIM histories of convergence for functions with moving kinks. (a, b) $u(x,\mu)=\max(\sin(x+\mu),0)$; and (c, d) $u(x,\mu)=|x+\mu|$ .}
\label{fig:kinky_functions}
\end{figure}

When $u(x,\mu)=\max(\sin(x+\mu),0)$ or $|x+\mu|$ (rows 4 and 5 of Table \ref{tab:tgpt_functions}), we have at least one kink (point of discontinuity of the derivative) and the location of the kink depends on $\mu$. See Figure \ref{fig:kinky_functions} for the convergence history for the EIM along with plots of the function snapshots for different $\mu$ values. We observe relatively slow convergence for the EIM, reaching two digits with $99$ snapshots. On the other hand, the TGPT-PINN only relies on one snapshot to achieve machine accuracy (see rows 4 and 5 of Table \ref{tab:tgpt_functions}) for each of these two parametric functions.

\subsubsection{Functions with moving discontinuities}
\label{sssec:function-example-2}
We next test the following example from \cite{Welper2017} which features a moving discontinuity as explained in Section \ref{sec:example},
$$
u(x, \mu):=\psi\left(\frac{x}{0.4+\mu}-1\right), \quad \psi(x):= \begin{cases}\exp \left(-\frac{1}{1-x^2}\right) & -1 \leq x<-\frac{1}{2}, \\ 0, & \text { else}.\end{cases}
$$
The summary of results is given by row 6 of Table \ref{tab:tgpt_functions}.  The function $\psi(x)$ is discontinuous at $x=-1/2$ which means that $u(x,\mu)$ is discontinuous at a parameter-dependent location $x=(\mu+0.4)/2$. Furthermore, we note that when $\mu > -0.4$, $x_\ast=(\mu+0.4)/2>0$ is a discontinuity point of $u(x,\mu)$ which is clamped to zero for $x > x_\ast$. On the other hand, when $\mu < -0.4$, then $x_\ast<0$ is the discontinuity point of $u(x,\mu)$, and $u$ is clamped to zero for $x < x_\ast$. 

\begin{figure}[htbp]
\centering
\includegraphics[width = 0.49\textwidth]{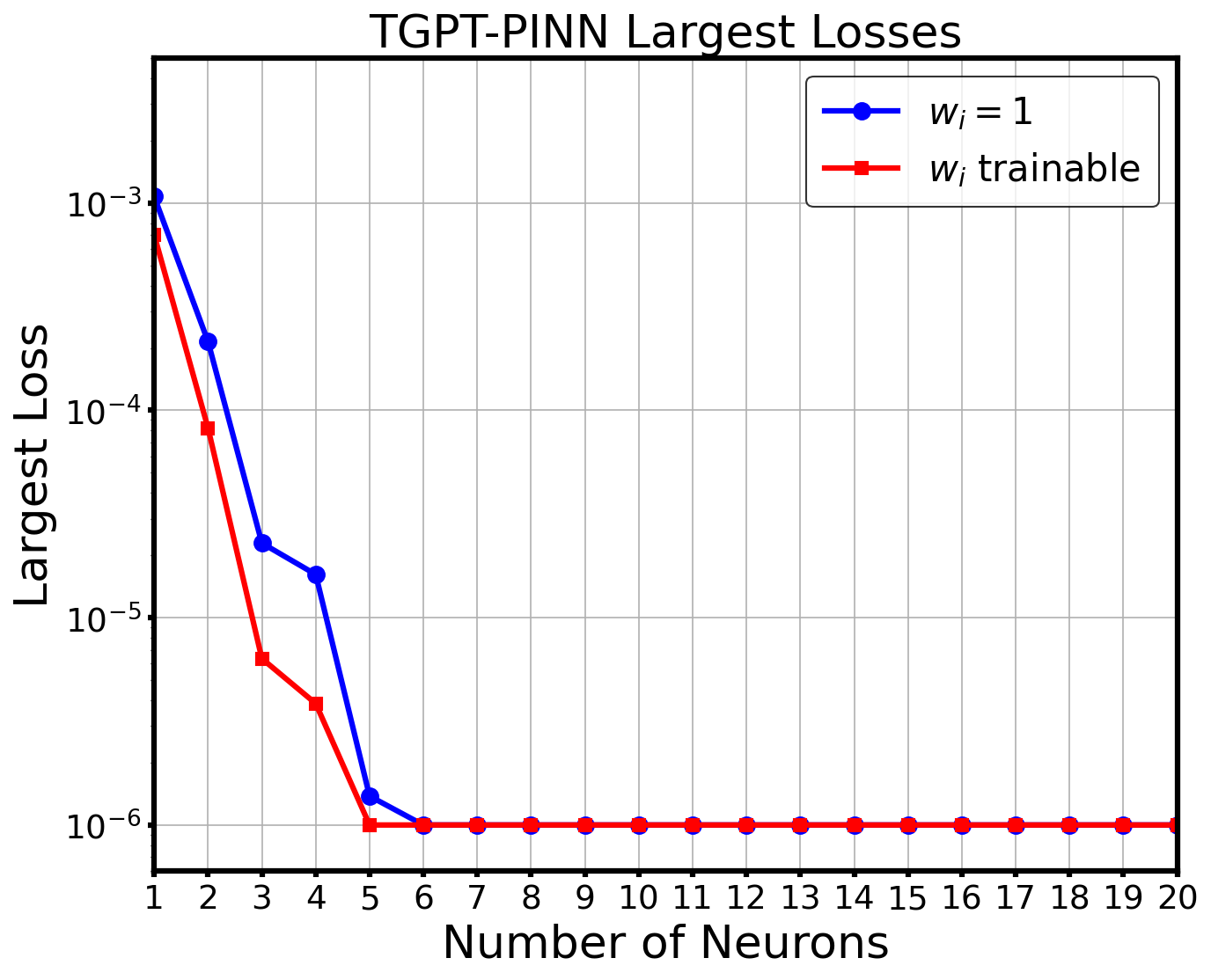}
\includegraphics[width=0.49\textwidth]{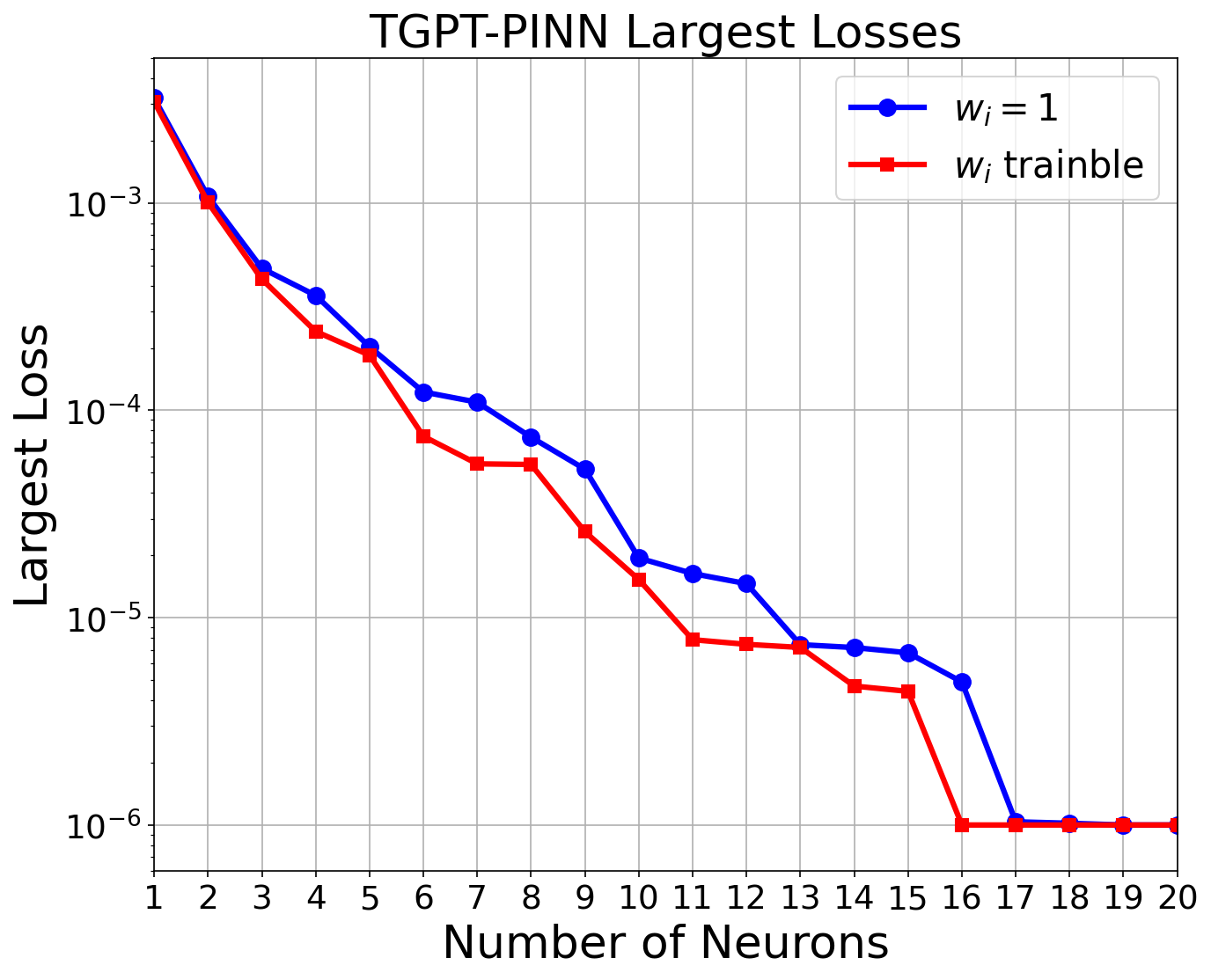}
\put(-360,25){\makebox(0,0){{(a)}}}
\put(-160,25){\makebox(0,0){{(b)}}}
\caption{Results from Section \ref{sssec:function-example-2} paired with row 6 of Table \ref{tab:tgpt_functions}: TGPT-PINN histories of convergence when the number of neurons increases for functions with a moving discontinuity, $u(x, \mu)=\psi\left(\frac{x}{0.4+\mu}-1\right)$. (a) $\mu \in [0,1]$ with 101 equispaced points; (b) $\mu \in [-1, 1]$ with 201 equispaced points.}

\label{fig:movingjump_hist}
\end{figure}
\begin{figure}[htbp]
\centering
\includegraphics[width=0.49\textwidth]{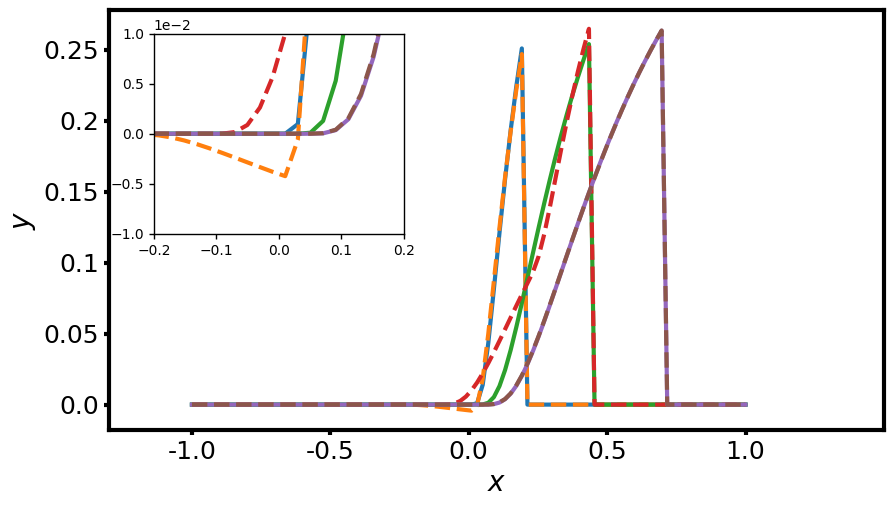}
\put(-15,25){\makebox(0,0){{(a)}}}
\includegraphics[width=0.49\textwidth]{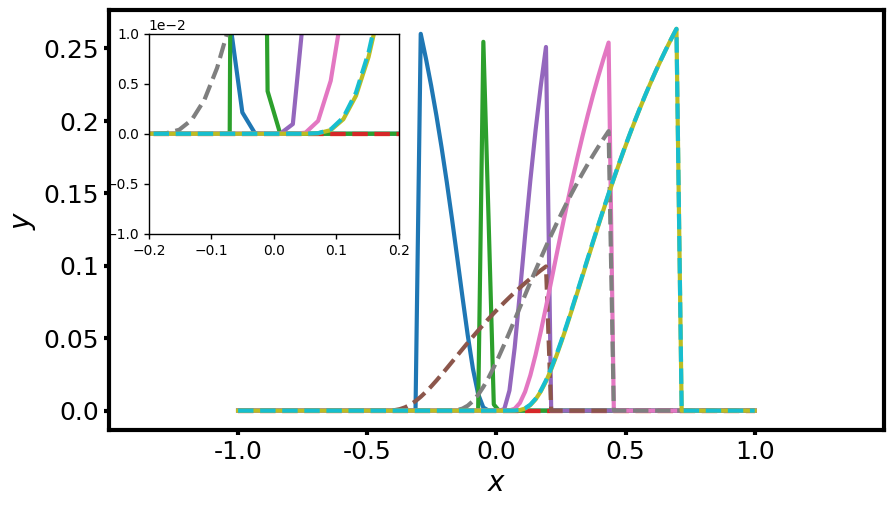}
\put(-15,25){\makebox(0,0){{(d)}}}\\
\includegraphics[width=0.49\textwidth]{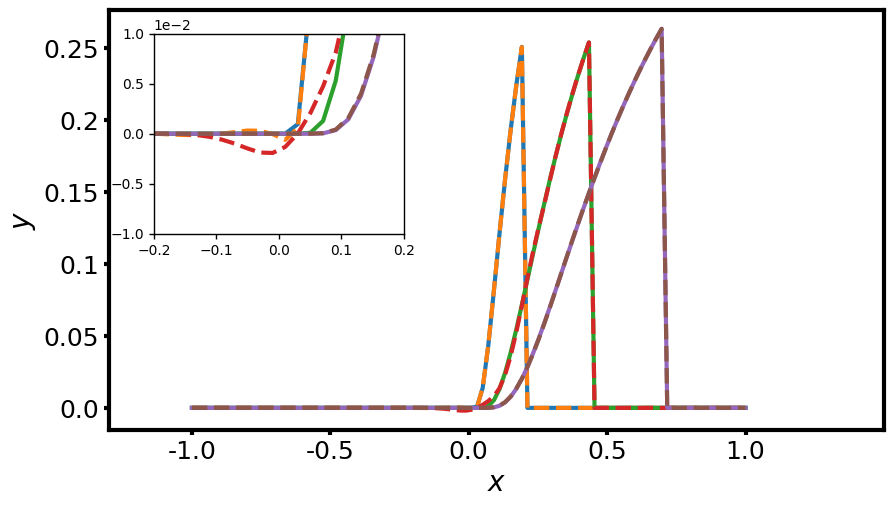}
\put(-15,25){\makebox(0,0){{(b)}}}
\includegraphics[width=0.49\textwidth]{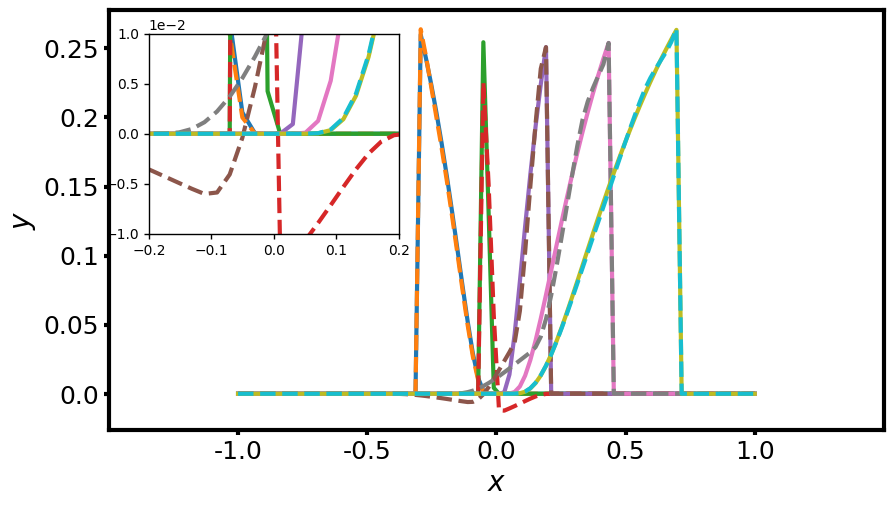}
\put(-15,25){\makebox(0,0){{(e)}}}\\
\includegraphics[width=0.49\textwidth]{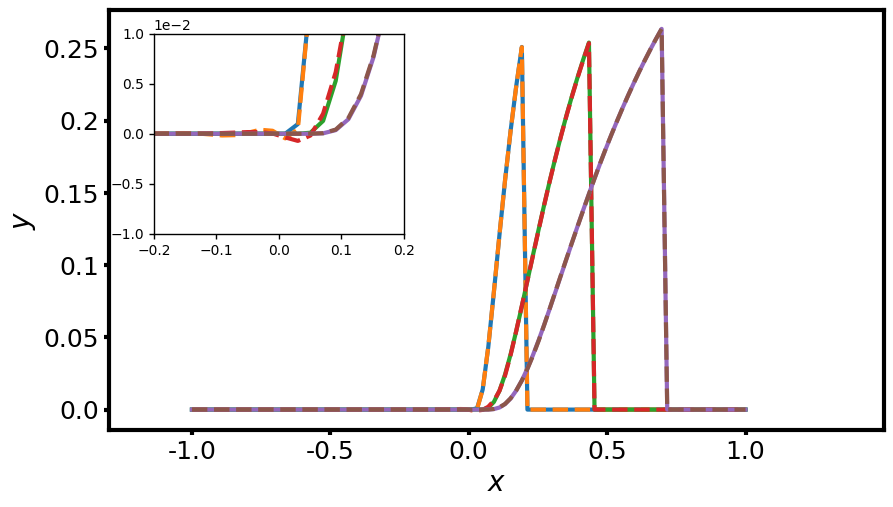}
\put(-15,25){\makebox(0,0){{(c)}}}
\includegraphics[width=0.49\textwidth]{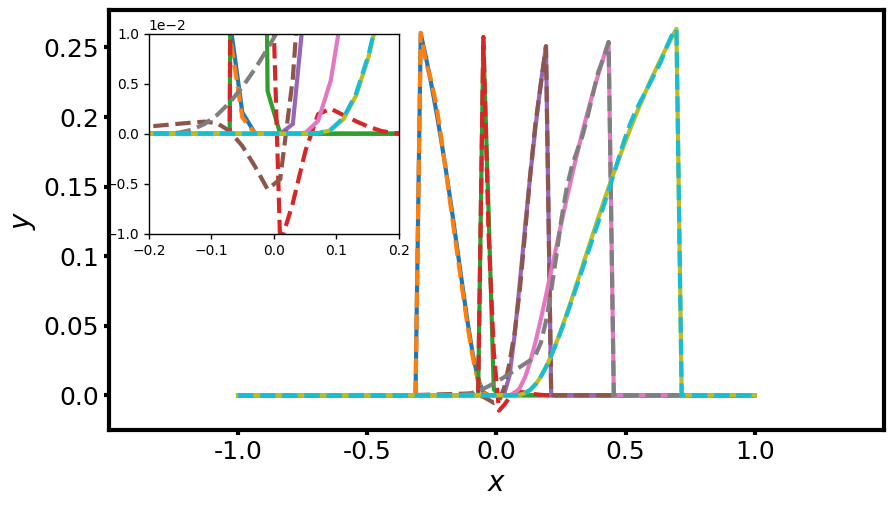}
\put(-15,25){\makebox(0,0){{(f)}}}
\caption{Results from Section \ref{sssec:function-example-2} paired with row 6 of Table \ref{tab:tgpt_functions}: Function with a moving discontinuity $u(x,\mu)=\psi\left(\frac{x}{0.4+\mu}-1\right)$ and its approximation with 1, 6, and 10 neurons (top to bottom). The solid lines describe the exact solutions and the dash lines are the TGPT-PINN solutions for different $\mu$. (a-c) training set $\mu \in [0,1]$ with 101 equispaced points and solutions of $\mu=0.005, 0.505$ and $0.995$ being displayed;  (d-f) training set $\mu \in [-1,1]$  with 201 equispaced points and solutions of $\mu=-0.995, -0.505, 0.005, 0.505$ and $0.995$ being displayed.}
\label{fig:movingjump_snap} 
\end{figure}

 We observe in row 6 of Table \ref{tab:tgpt_functions} that the EIM fails to capture this family of functions even with all but one snapshots from the entire discretized domain. On the other hand, the TGPT-PINN successfully reaches 6 digits of accuracy by 10 neurons. Since this example is more difficult than the previous ones, we conduct a further example to investigate the impact of the parameter space: We compare the TGPT-PINN performance using $\left(\mathcal{D}, \Xi_{\mathrm{train}}\right) = ([0,1], 101)$, and  $\left(\mathcal{D}, \Xi_{\mathrm{train}}\right) = \left([-1,1], 201\right)$. The latter is the setup in Table \ref{tab:tgpt_functions}. As we increase the number of neurons in the hidden layer, Figure \ref{fig:movingjump_hist} shows the change of the L2 norm of the error committed by the approximate TGPT-PINN solution. The exact function and its approximation with 1, 6, and 10 neurons are shown in Figure \ref{fig:movingjump_snap}. The visual agreement between the exact function and its TGPT-PINN approximation is confirmed by the computed error indicating an accuracy of $5$ to $6$ digits.

Next, we examine the convergent network parameters for the more challenging case with $\mu \in [-1,1]$ to confirm that it aligns with our theoretical understanding of this discontinuous function. We examine the function approximation with $10$ neurons 
\begin{align}\label{eq:temp}
u(x;\mu) \approx \sum_{i=1}^{10} c_i(\mu) u(w_i(\mu) x + b_i(\mu), \mu^i),
\end{align}
where in order to line up the discontinuity, we know the exact shift should be $b_i^{\rm exact} = \frac{\mu^i+0.4}{2} -w_i \frac{\mu+0.4}{2}$.
We list in Table \ref{tab:disc_paras} the results of $\{(\mu^i, w_i, b_i, c_i): i = 1, \cdots, 10\}$ for an unseen parameter value $\mu = 0.595$, i.e., this parameter was not in the training set for offline construction. The left half of the table shows the results when we fix $w_i = 1$ and the right part corresponds to when the $w_i$'s are trainable parameters. We see that the network training is highly effective as all the $b_i$'s coincide with their exact value in both cases. It is interesting to note that, for both cases, relatively larger discrepancies $|b^i-b_i^{\rm exact}|$ occur when $c_i$ is relatively small which is a testament to the robustness of the algorithm.

\begin{table}[!htb]
\centering
\begingroup
\aboverulesep = 0mm 
\belowrulesep = 0mm
\renewcommand{\arraystretch}{1.3} 
  \resizebox{\textwidth}{!}{
\begin{tabular}{rcccccaaaaa}
  & \multicolumn{5}{c}{$\mu = 0.595$, $w_i = 1$ fixed, final loss 1.23e-05} &
    \multicolumn{5}{a}{$\mu = 0.595$, $w_i$ trainable, final loss 8.12e-06} \\
  $i$ & $\mu^i$ & $w_i$ & $|b^i - b_i^{\mathrm{exact}}|$ & $b_i$ & $c_i$ & 
        $\mu^i$ & $w_i$ & $|b^i - b_i^{\mathrm{exact}}|$ & $b_i$ & $c_i$ \\\toprule
  $1$ & 1.0&1&1.5e-08&0.20&-2.34e-02&1.0&1.14& 2.98e-08 & 0.13&3.85e-01\\
  $2$ & -0.99&1&1.8e-07&-0.80&9.20e-42&-0.99&0.47& 1.79e-07 & -0.53&-5.17e-42\\
  $3$ & 0.03&1&0.0e-00&-0.28&-1.71e-1&-0.09&0.77& 5.96e-08 & -0.23&1.34e-01\\
  $4$ & -0.55&1&1.2e-07&-0.57&8.52e-41&-0.55&0.62& 1.19e-07 & -0.38&-3.45e-41\\
  $5$ & -0.29&1&3.0e-08&-0.44&-2.62e-02&-0.33&0.69& 0.00e-00 & -0.31&7.21e-03\\
  $6$ & 0.44&1&1.5e-08&-0.08&8.84e-01&0.31&0.90& 1.49e-08 & -0.10&9.47e-01\\
  $7$ & -0.47&1&1.2e-07&-0.53&-4.40e-41&-0.47&0.64& 8.94e-08 & -0.36&2.51e-41\\
  $8$ & -0.75&1&1.2e-07&-0.67&-4.97e-41&-0.72&0.56& 1.19e-07 & -0.44&2.15e-41\\
  $9$ & -0.17&1&3.0e-08&-0.38&-7.01e-02&-0.25&0.72& 8.94e-08 & -0.28&-2.79e-02\\
  $10$ & -0.33&1&3.0e-08&-0.46&8.50e-03&0.08&0.83& 7.45e-08 & -0.17&-4.44e-01
\end{tabular}
  }
\endgroup
  \caption{Results for Section \ref{sssec:function-example-2} paired with row 6 of Table \ref{tab:tgpt_functions}: Trained TGPT-PINN network parameters at $\mu = 0.595$ for the function with moving discontinuities.}
\label{tab:disc_paras}
\end{table}

\subsubsection{2D functions close to being degenerate}
\label{sssec:function-example-3}

\begin{figure}[htbp]
	\centering
\includegraphics[scale=0.5]{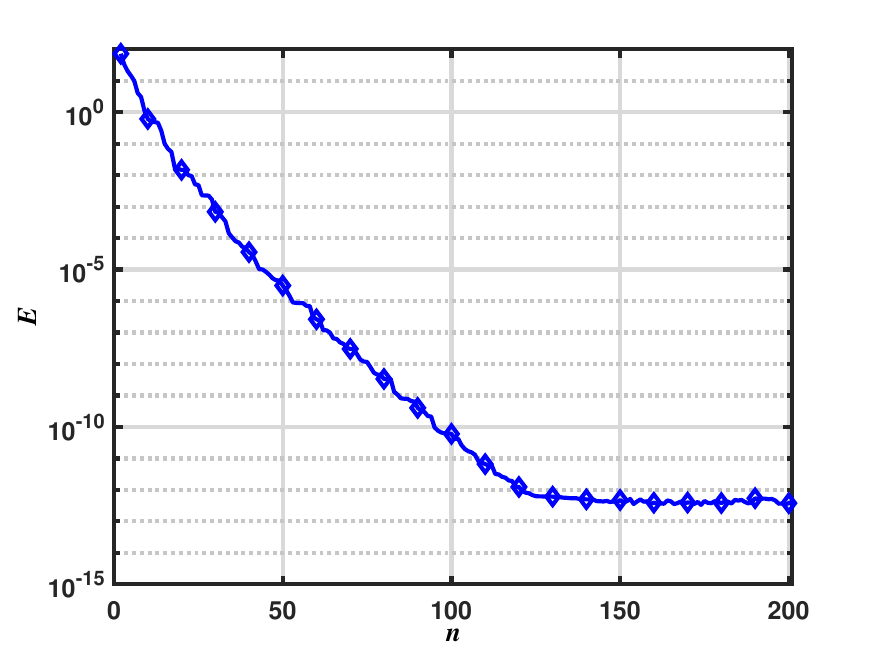}
\caption{Results for Section \ref{sssec:function-example-3} paired with row 7 of Table \ref{tab:tgpt_functions}: EIM history of convergence for the 2D functions close to being degenerate.}
\label{fig:EIM_2dDeg}
\end{figure}

\begin{figure}[htbp]
\subfigure[TGPT-PINN Loss]{\includegraphics[scale=0.31]{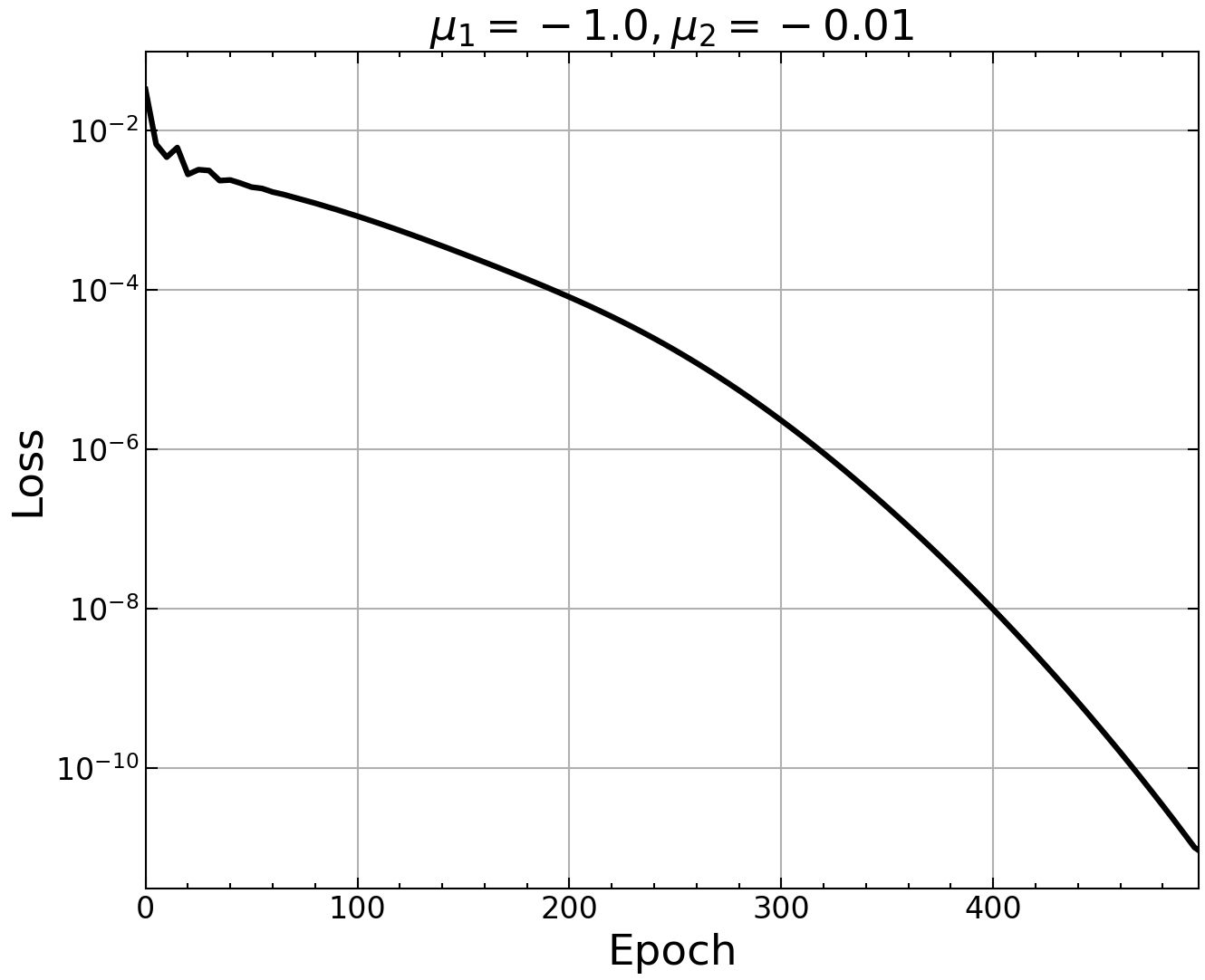}}
\subfigure[TGPT-PINN error]{\includegraphics[scale=0.3]{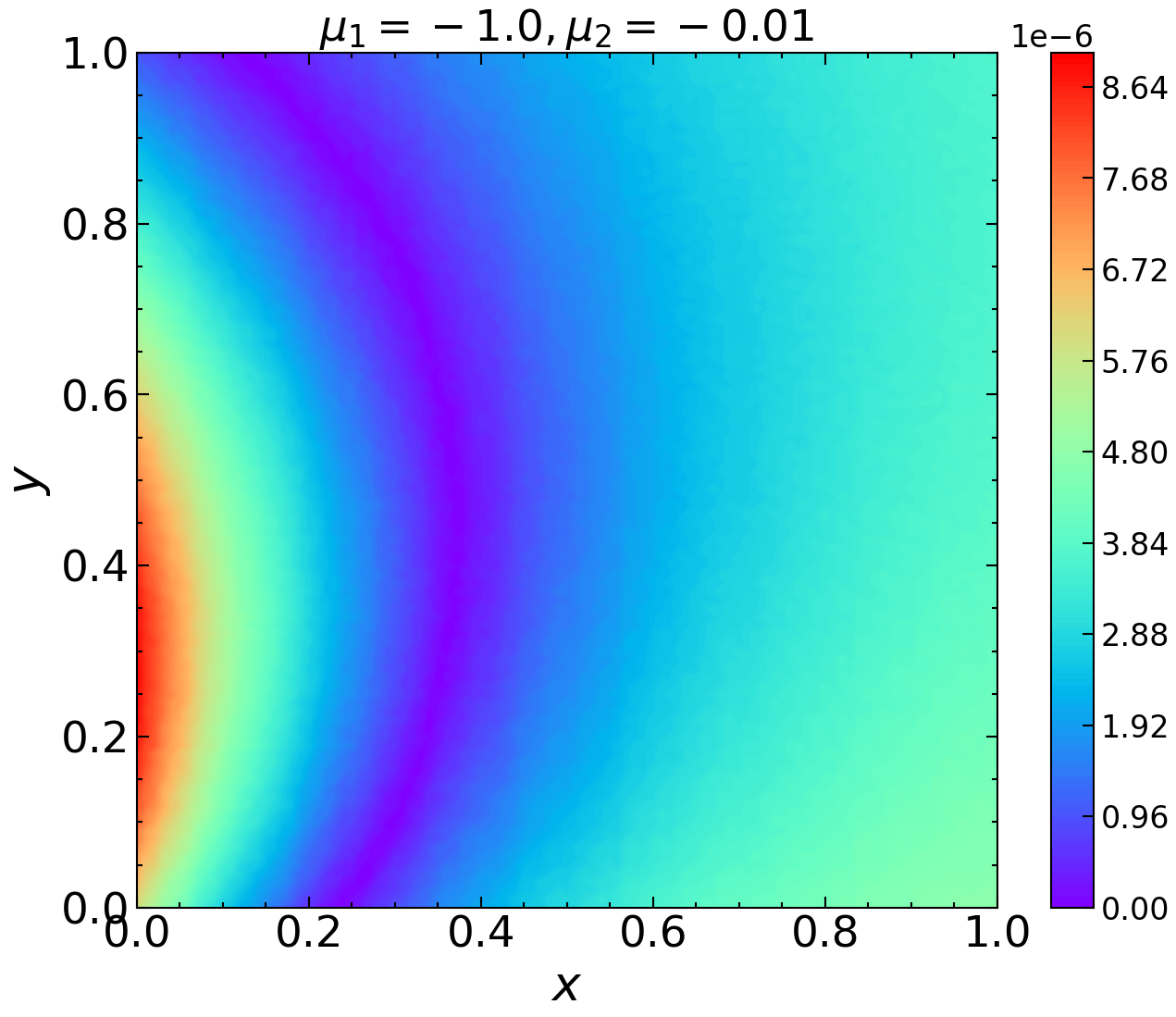}}\\
\subfigure[TGPT-PINN Loss]{\includegraphics[scale=0.31]{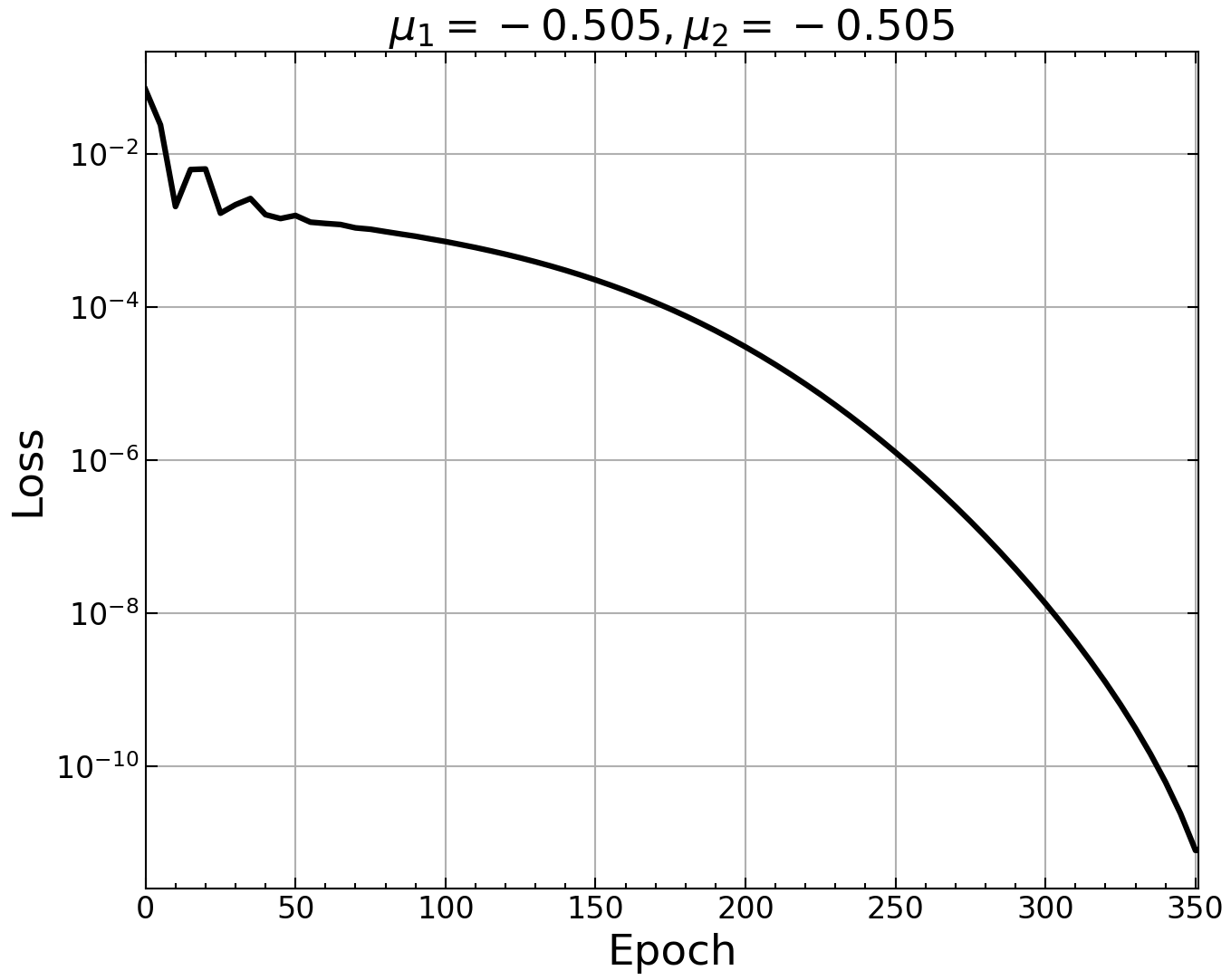}}
\subfigure[TGPT-PINN error]{\includegraphics[scale=0.3]{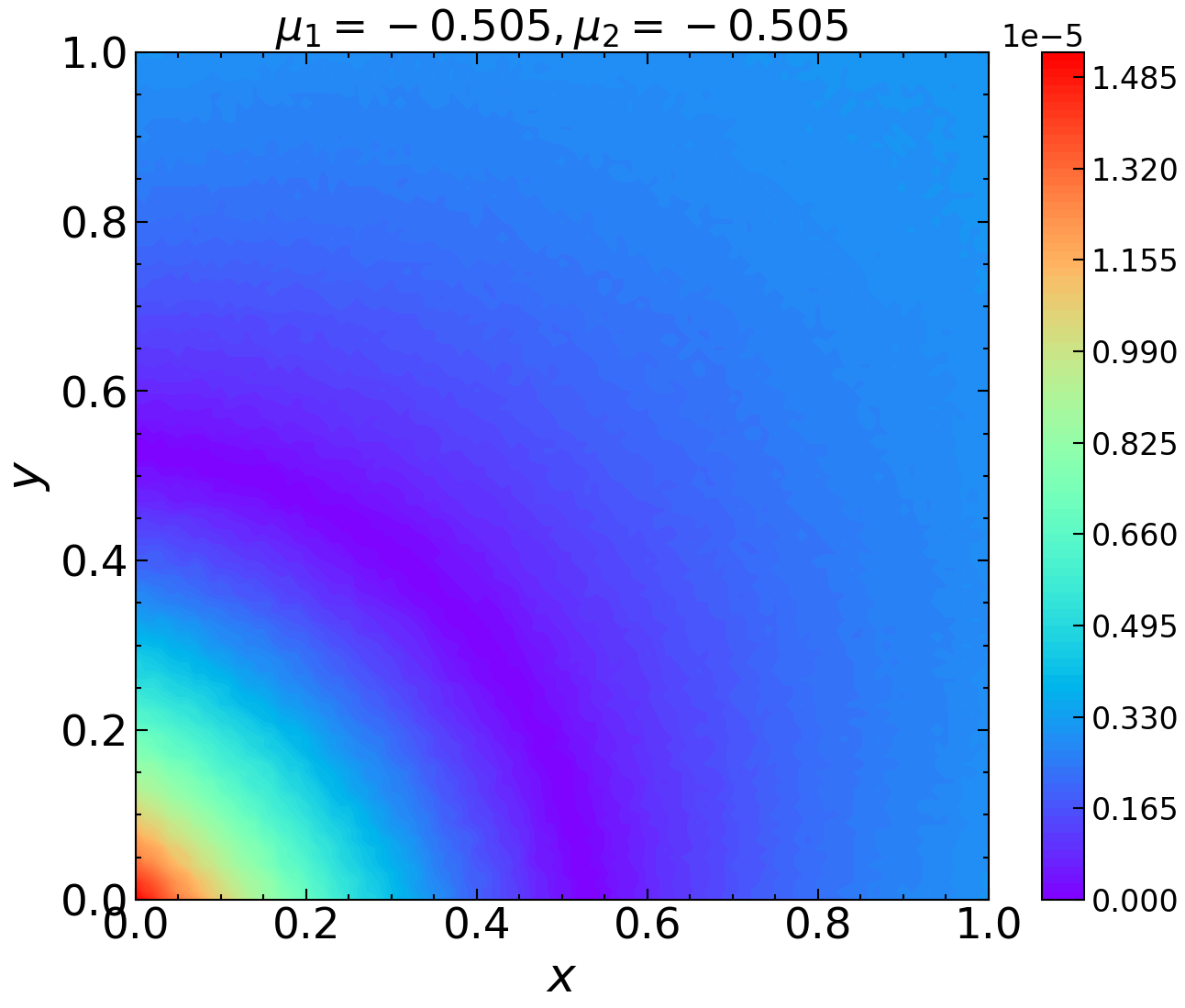}}
\caption{Results for Section \ref{sssec:function-example-3} paired with row 7 of Table \ref{tab:tgpt_functions} for the 2D functions close to being degenerate: TGPT-PINN error (right) committed by the TGPT-PINN solution with one neuron  at a corner and the center of the parameter domain. Shown on the left are the corresponding training histories.}
\label{fig:error_degenerate}
\end{figure}

We now consider the classical 2-dimensional nonlinear and nonaffine function from \cite{Grepl_Maday_Nguyen_Patera}, 
\[
u(x,y;\mu_1,\mu_2)=\frac{1}{\sqrt{(x-\mu_1)^2+(y-\mu_2)^2}} \mbox{ on } (x,y)\in [0,1]^2,
\]
parameterized by $(\mu_1,\mu_2)\in [-1,-0.01]^2$. Due to the setup of the physical and parameter domains, this family contains functions that are very close to being singular for certain parameter value. The EIM (see Figure \ref{fig:EIM_2dDeg}) needs about $120$ basis functions to reach an accuracy of $10^{-12}$ while the TGPT-PINN \textit{requires only one neuron,} $u(x,y;(-1,-1))$, to achieve the same accuracy for all parameters in the parameter space discretized by a $21 \times 21$ grid. The solution, the errors, and the training history corresponding to two locations of the parameter domain committed by the TGPT-PINN surrogate are plotted in Figure \ref{fig:error_degenerate}. 

\subsection{TGPT-PINN for parametric PDEs}\label{ssec:results-ppde}

In this section, we test the TGPT-PINN algorithm on three pPDEs, the linear transport equation parameterized by the wave speed (first example, \cref{sssec:results-ppde-transport}), the nonlinear 1D-reaction equation parameterized by the reaction coefficiant (second example, \cref{sssec:results-ppde-reaction}), and the 1D nonlinear reaction-diffusion equation parameterized by the viscosity and reaction coefficient (third example, \cref{sssec:results-ppde-rd}). The overall conclusions are as follows: 
\begin{itemize}[itemsep=0pt,topsep=0pt,leftmargin=15pt]
  \item {\Cref{sssec:results-ppde-transport}}: For wave-like problems when linear reduction does not work, the TGPT-PINN can achieve significant accuracy because transport-based parameter dependence can be efficiently approximated by the novel transform layer. For our particular example, \textit{only one neuron is needed to achieve machine precision accuracy}.
  \item {\Cref{sssec:results-ppde-reaction,sssec:results-ppde-rd}}: On problems when linear model reduction can work effectively, and in particular when the GPT-PINN performs well, the TGPT-PINN works using a very small number of snapshots (1 and 3, respectively), and achieves better accuracy with this small number of snapshots than the GPT-PINN with a much larger number of snapshots.
\end{itemize}
This vast improvement is a manifestation of the novel design of the TGPT-PINN: the addition of a trainable transform layer to a network of pre-trained networks.

\subsubsection{Transport equation}
\label{sssec:results-ppde-transport}
Consider the following example\cite{Welper2017},
\begin{align*}
  &u_t+\nu u_x=0, \quad \text { for } 0<t<2, x \in [-1,1], \\
& u(x, 0)=g(x):=\left\{\begin{array}{ll}
0, & -0.5 < x < 0.5, \\
1, & \text{else}.
\end{array} \right.
&
\end{align*}
The exact solution to the problem is given by $u(x,t)=g(x-\nu t)$, so that the parameter $\nu$ determines the direction and speed of propagation for the solution of the equation.

\noindent {\bf Sampling for the PINN loss function.} For fixed $\nu$, this problem is challenging to train a PINN on because of the discontinuity in $(x,t)$. To ameliorate this difficulty, we choose sampling points for the loss function in a standard randomized fashion over the $(x,t)$ domain, but we additionally place extra training points near the discontinuity. For  details, for the first term of the loss function \eqref{eq:loss} corresponding to the PDE residual in the domain interior, a $400 \times 400$ equidistant grid is created over the domain, and 10,000 points are randomly subsampled as the training set, with a relatively concentrated distribution near the vicinity of the discontinuity bands. We choose 60\% of the points by randomly subsampling within an $(x,t)$ distance of 0.1 of the disconinuity. 20\% more points are subsampled from equidistant points whose $(x,t)$ distance from the discontinuity takes values in $[0.1,0.2]$. The remaining 20\% of the points are subsampled from the remaining points. A similar procedure is done for the initial condition training grid, where a total of 1000 points are randomly subsampled from $10^4$ equispaced points in $x$, with a higher concentration of points subsampled close to the discontinuity.

\begin{figure}[htbp]
  \centering
  \subfigure[PINN solution]{\includegraphics[width=0.32\linewidth]{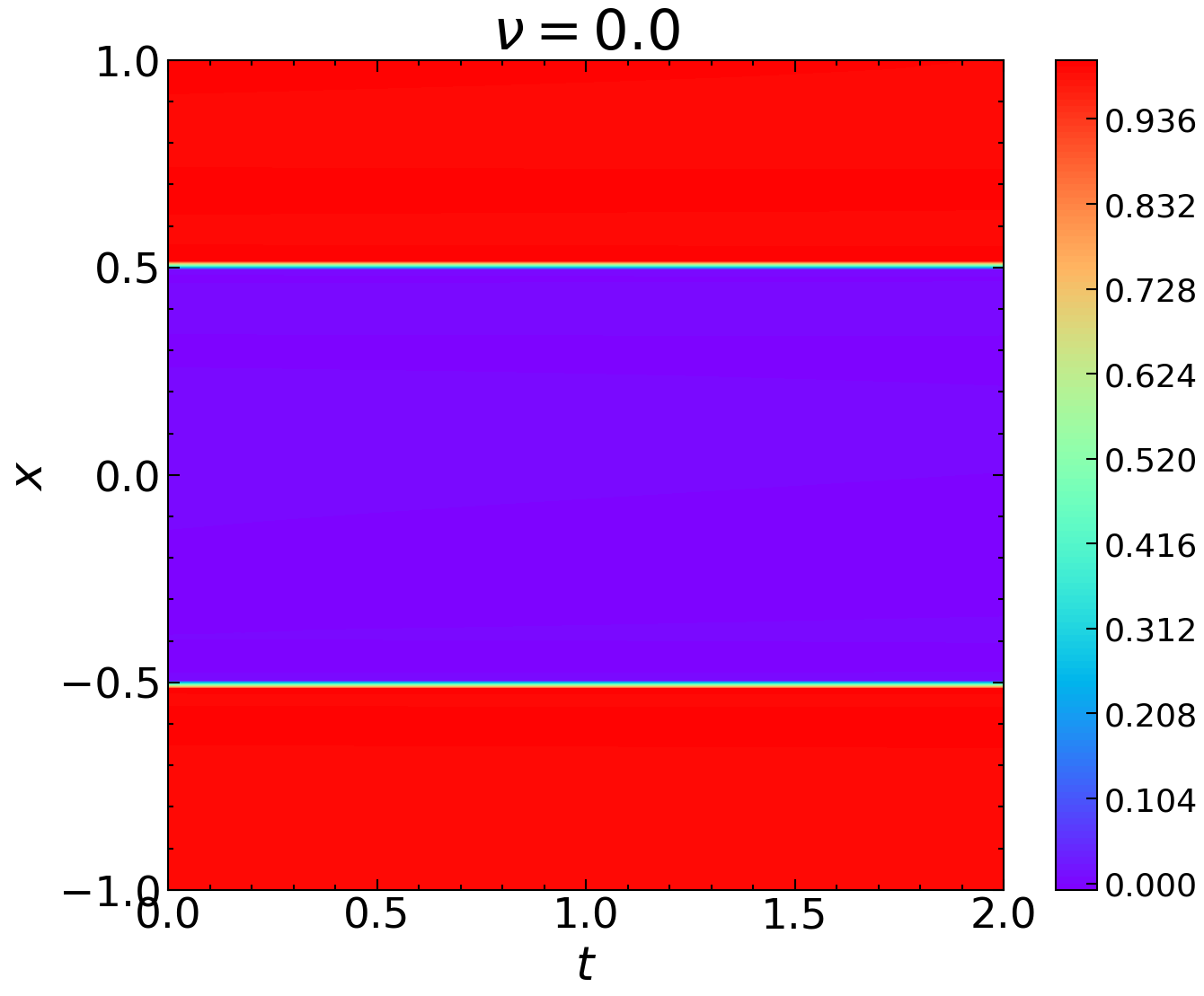}}
  \subfigure[PINN error]{\includegraphics[width=0.335\linewidth]{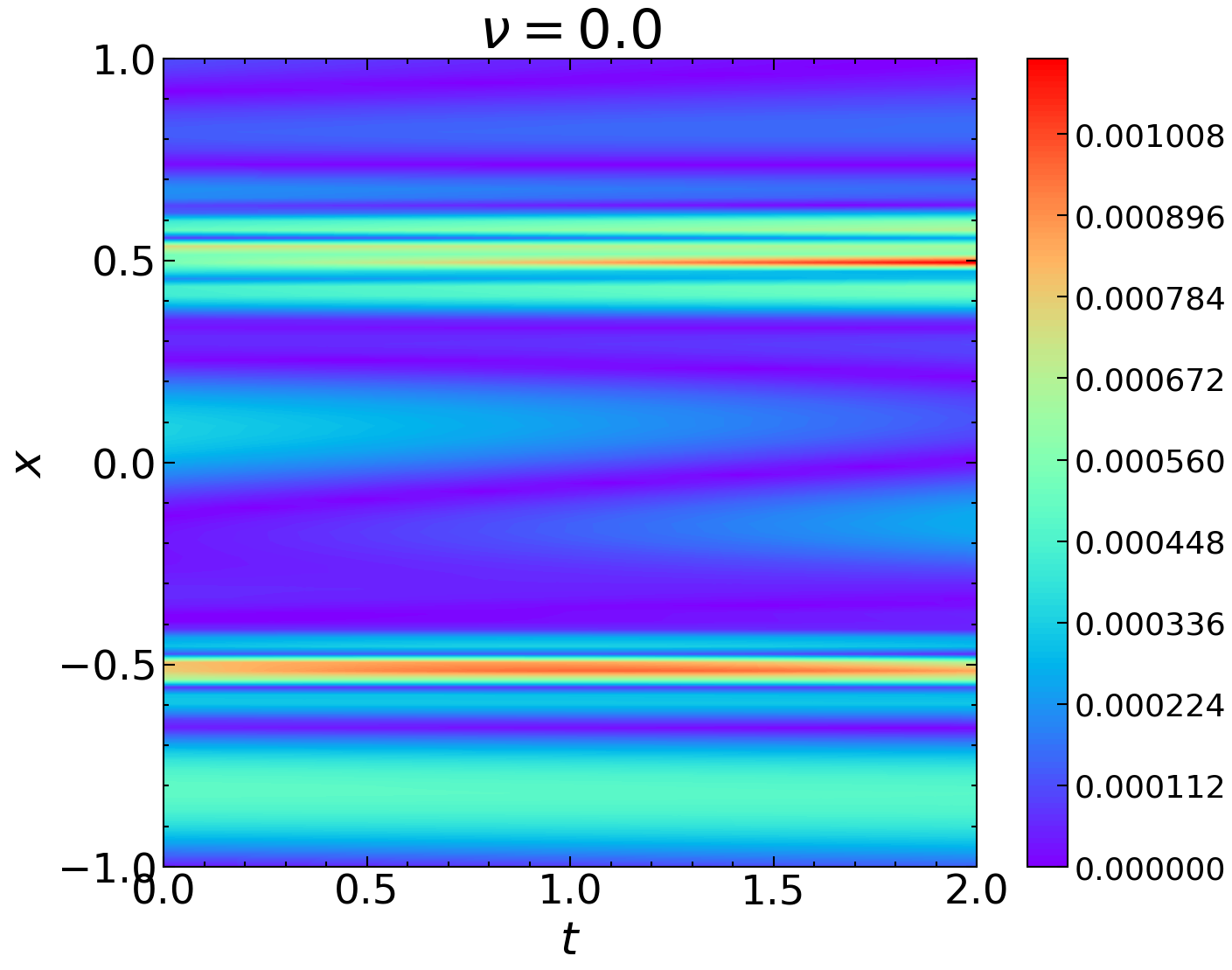}}
  \subfigure[PINN Loss]{\includegraphics[width=0.32\linewidth]{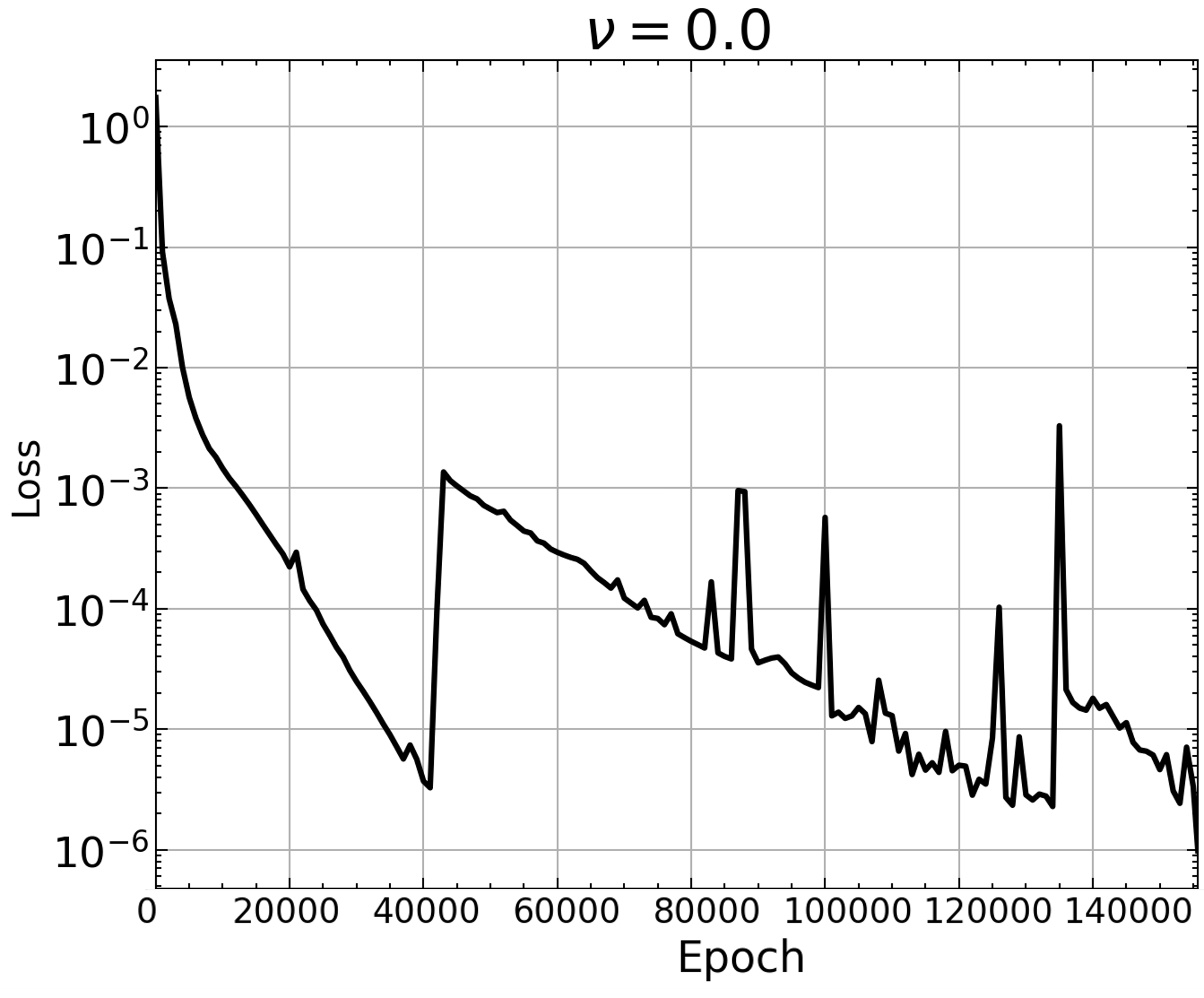}}
  \subfigure[TGPT-PINN solution]{\includegraphics[width=0.32\linewidth]{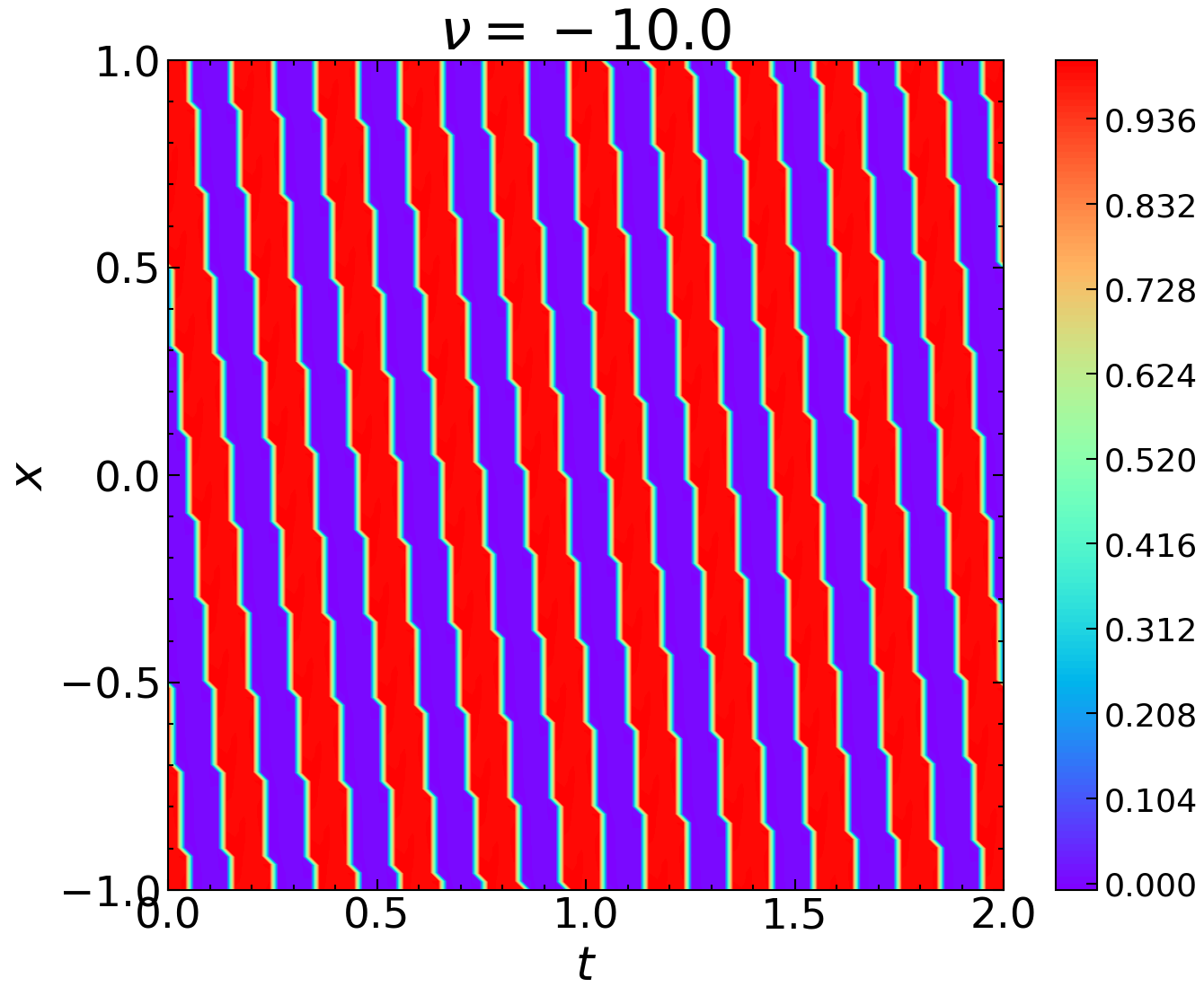}}
  \subfigure[TGPT-PINN error]{\includegraphics[width=0.33\linewidth]{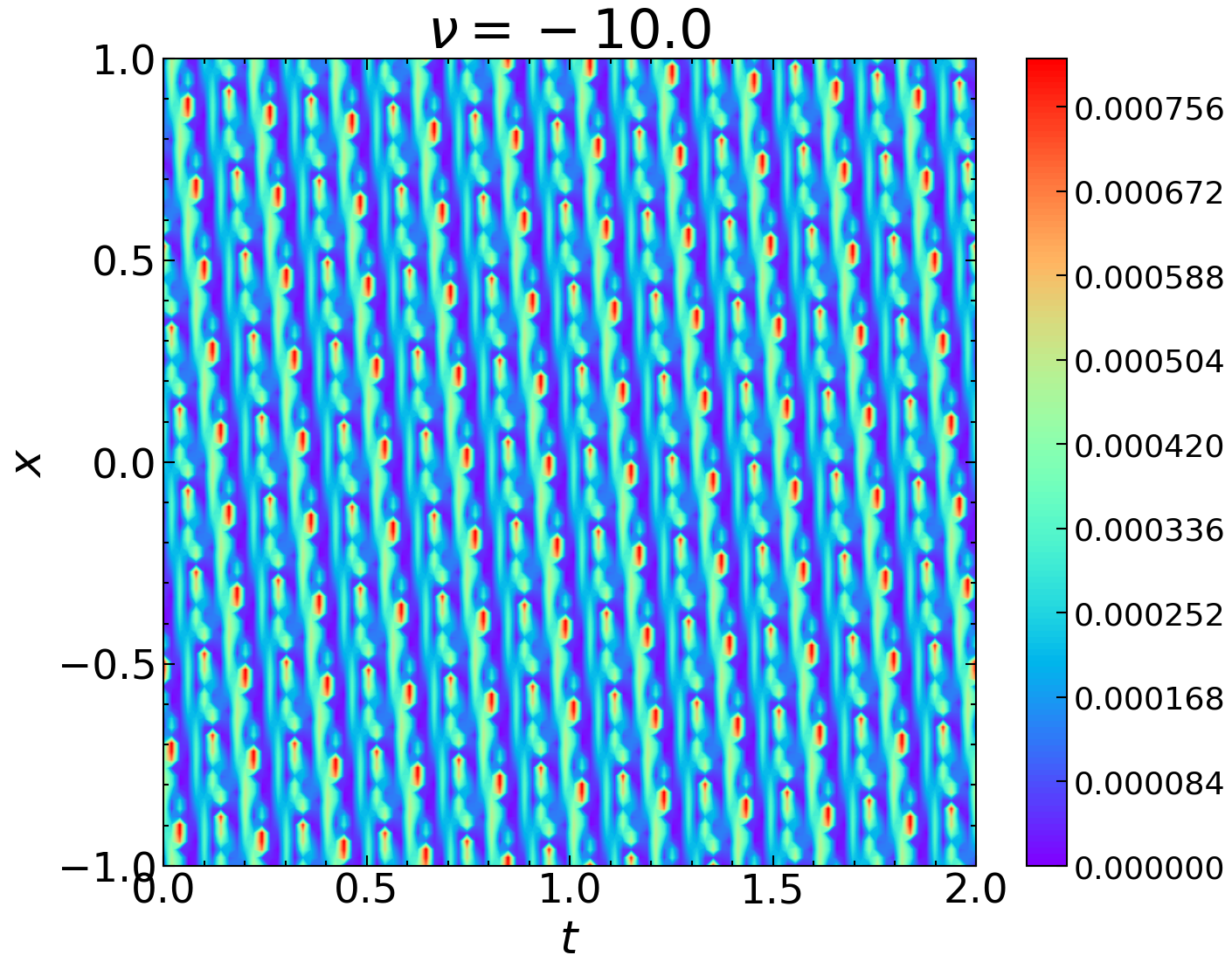}}
  \subfigure[TGPT-PINN Loss]{\includegraphics[width=0.33\linewidth]{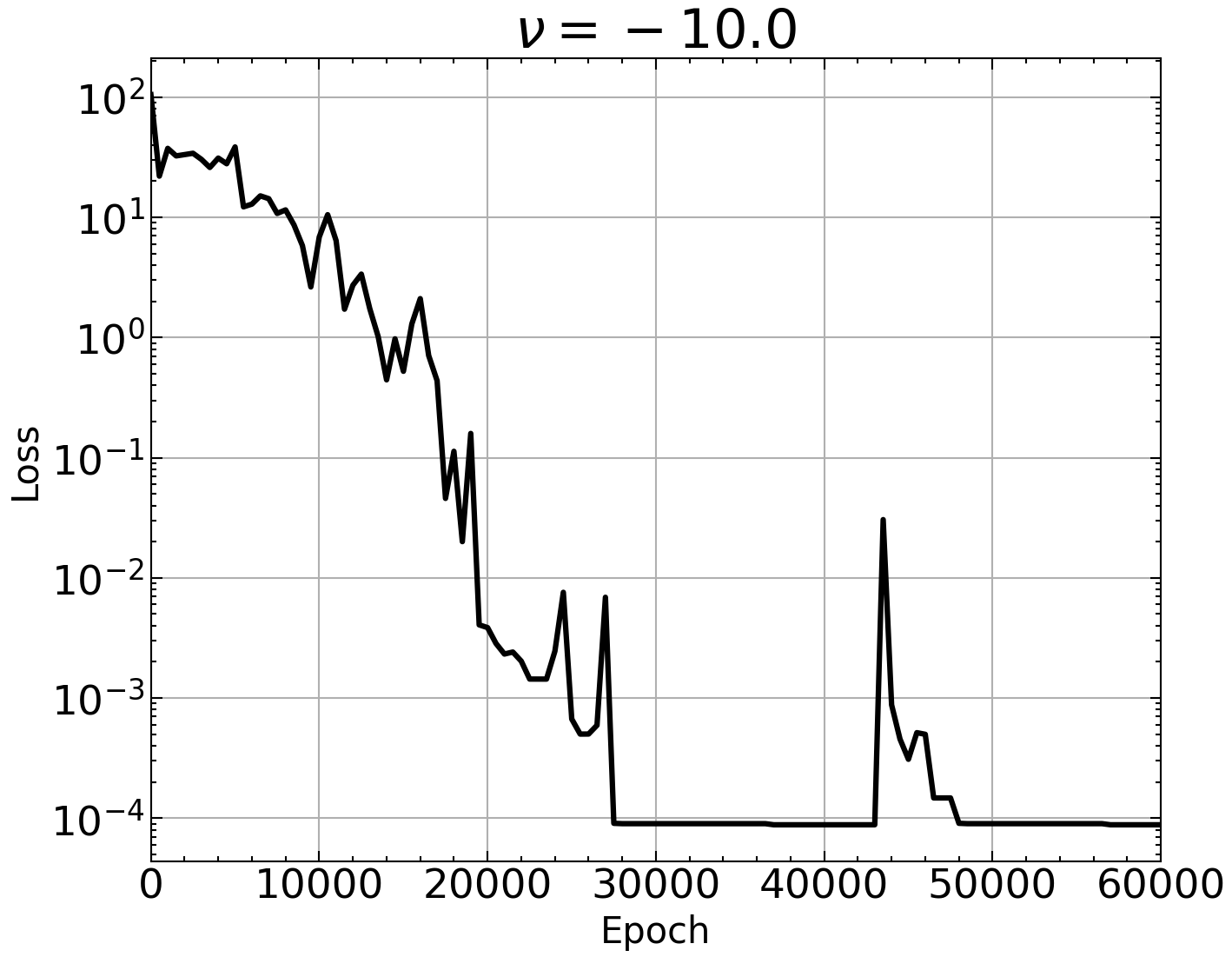}}
  \subfigure[TGPT-PINN solution]{\includegraphics[width=0.32\linewidth]{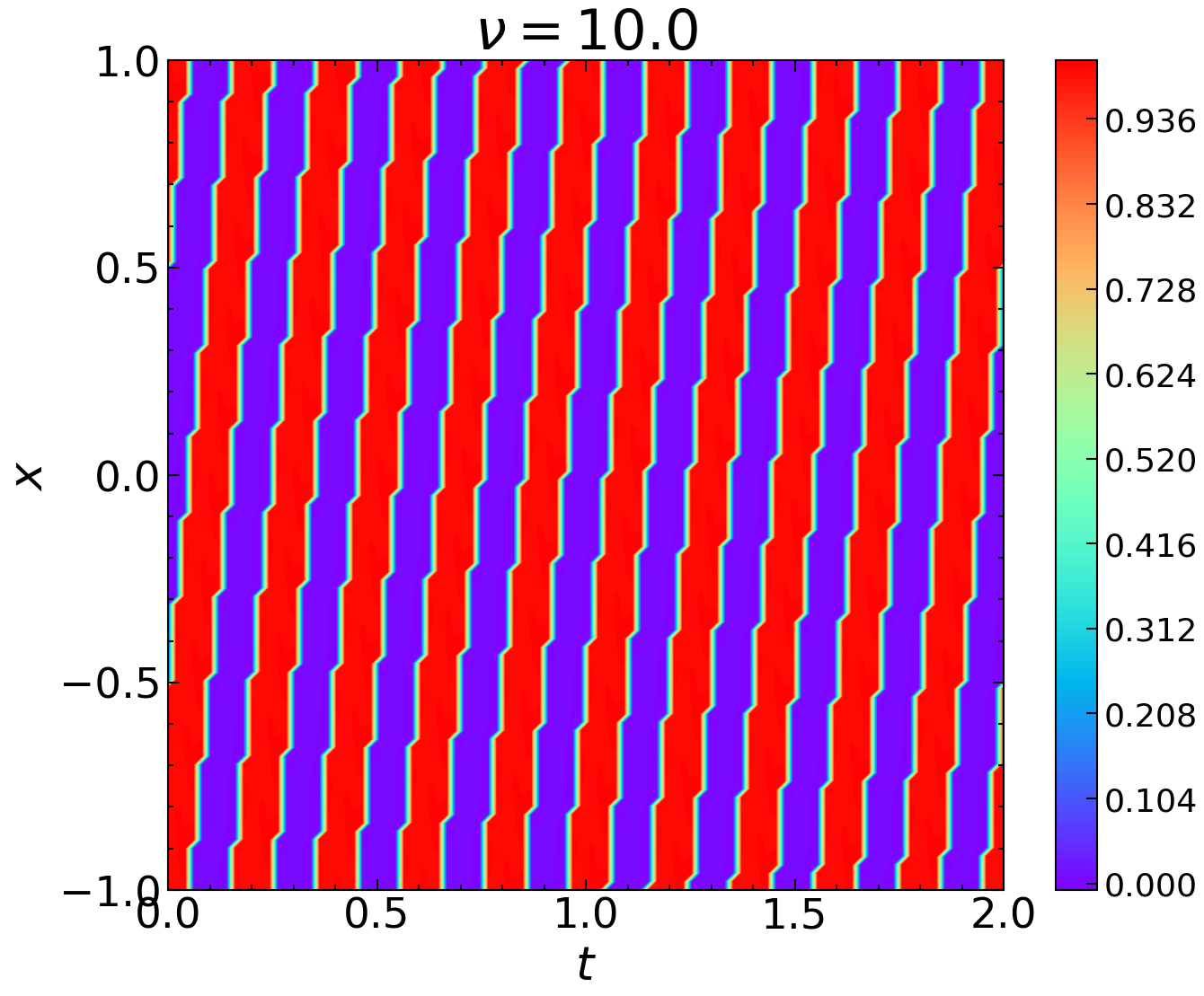}}
  \subfigure[TGPT-PINN error]{\includegraphics[width=0.33\linewidth]{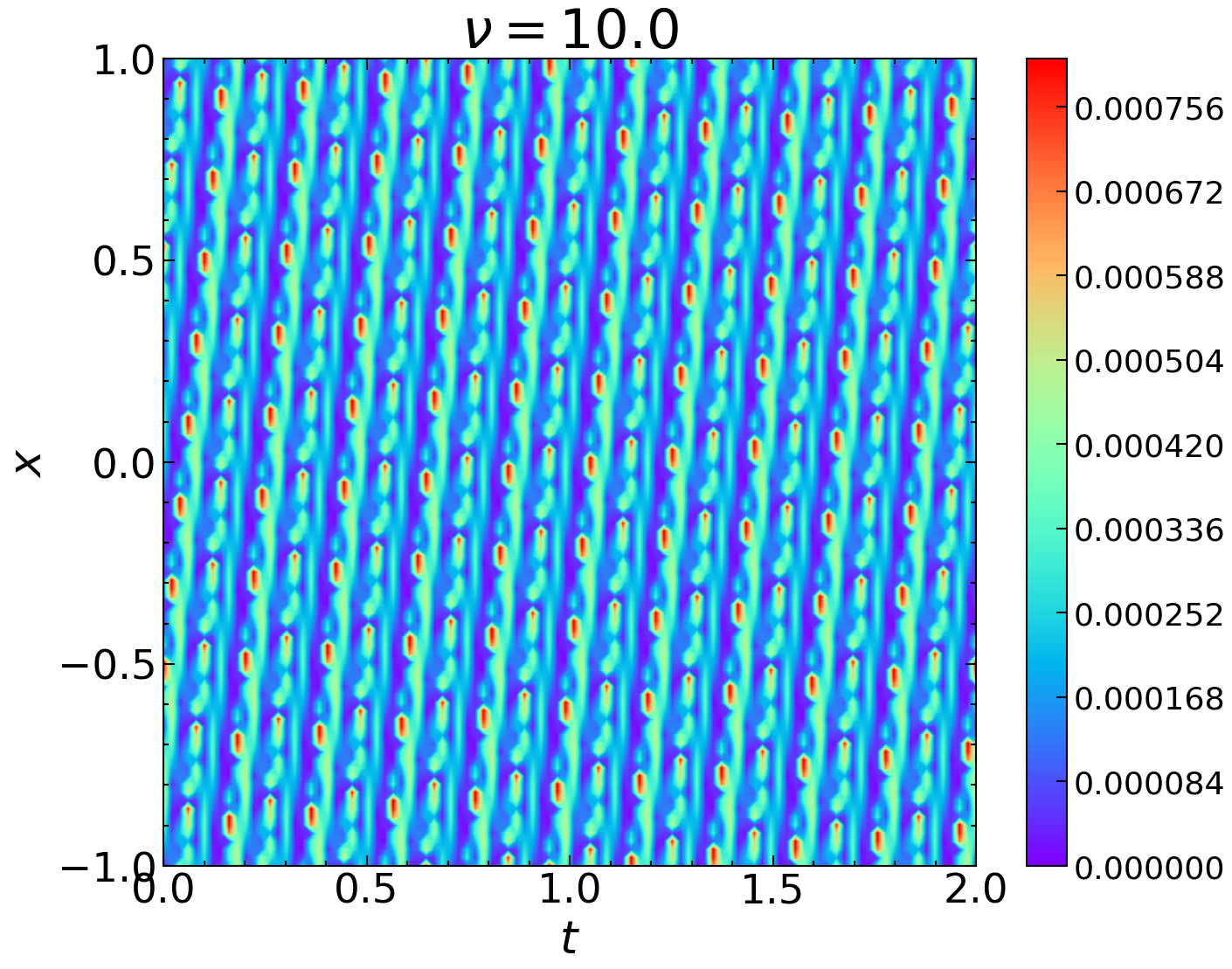}}
  \subfigure[TGPT-PINN Loss]{\includegraphics[width=0.322\linewidth]{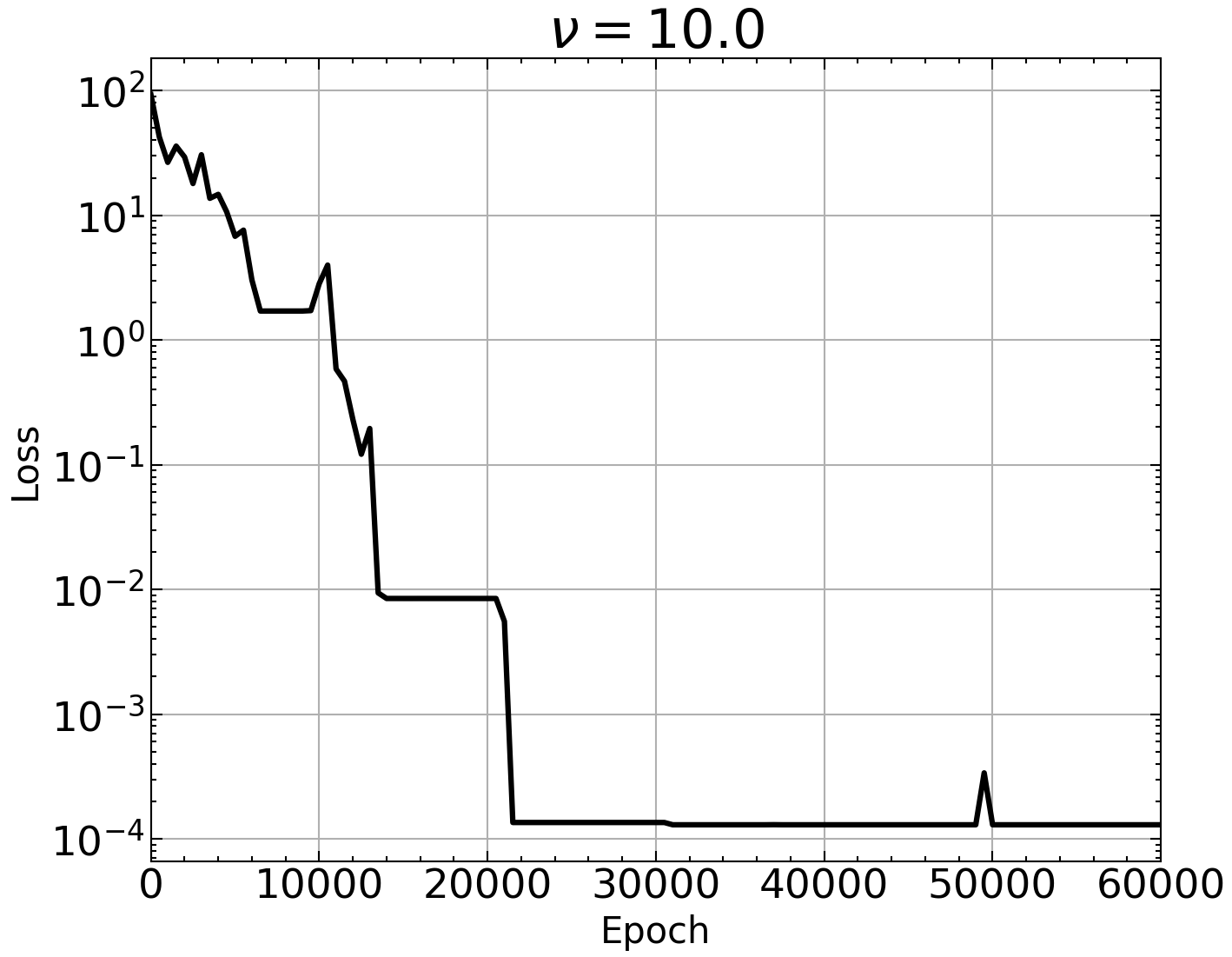}}
  \caption{Results for \cref{sssec:results-ppde-transport}. Top row: Full PINN for $\nu = 0.0$. Middle row: TGPT-PINN at $\nu=-10.0$. Bottom row: TGPT-PINN for $\nu = 10.0$. Both middle and bottom are approximated using the only PINN from the top row.}
  \label{fig:tgpt-pinn-transport}
\end{figure}

\noindent{\bf PINN results.} The full PINN for this transport equation is a fully connected neural network with [2, 20, 20, 20, 1] neurons across its layers. The activation function used is $tanh(\cdot)$. The PINN is trained with a learning rate of 0.001, a maximum of 4e5 iterations, a tolerance of 1e-6, and the Adam optimizer. We take $\nu=0.0$ and train a full PINN. The reduction in loss of the network in a single random experiment is shown in Figure\ref{fig:tgpt-pinn-transport} (a-c), showing that the PINN is relatively accurate.

\noindent {\bf TGPT-PINN results.} We conducted tests on a discretized parameter set of 41 equispaced samples over $\mu = \nu \in [-10, 10]$, with a stopping accuracy set to 1e-5, learning rate of 0.05, and a maximum iteration count of $10^5$. We use only $N=1$ snapshot, manually choosing $\mu^1 = 0$. The variation of loss, solutions, and errors under the two extreme parameter values, $\nu = \pm 10$ (the hardest parameter values to approximate) are shown in Figure \ref{fig:tgpt-pinn-transport} (d-i).

\subsubsection{1D-reaction PDE}\label{sssec:results-ppde-reaction}
The one-dimensional reaction problem is a hyperbolic PDE that is commonly used to model chemical reactions, 
\begin{equation}
\begin{gathered}
  \frac{\partial u}{\partial t}-\rho u(1-u)=0, \hskip 10pt (x,t) \in [0, 2\pi] \times [0, 1], \\
  u(x, 0)= h(x), \hskip 10pt  x \in [0, 2\pi], \\
  u(0, t)=u(2 \pi, t), \hskip 10pt t \geq 0,
\end{gathered}
\end{equation}
where $\rho > 0$ is the reaction coefficient. We take the initial condition as,
\begin{align}\label{eq:h-def}
  h(x) &= \exp \left(-\frac{(x-\pi)^2}{2(\pi / 4)^2}\right).
\end{align}
The equation has a simple analytical solution:
\begin{equation}
  u(x,t) =\frac{h(x) \exp (\rho t)}{h(x) \exp (\rho t)+1-h(x)}.
\label{eq:u_rea}
\end{equation}

\begin{figure}[htbp]
  \centering
  \subfigure[Exact solution]{\includegraphics[width=0.24\linewidth]{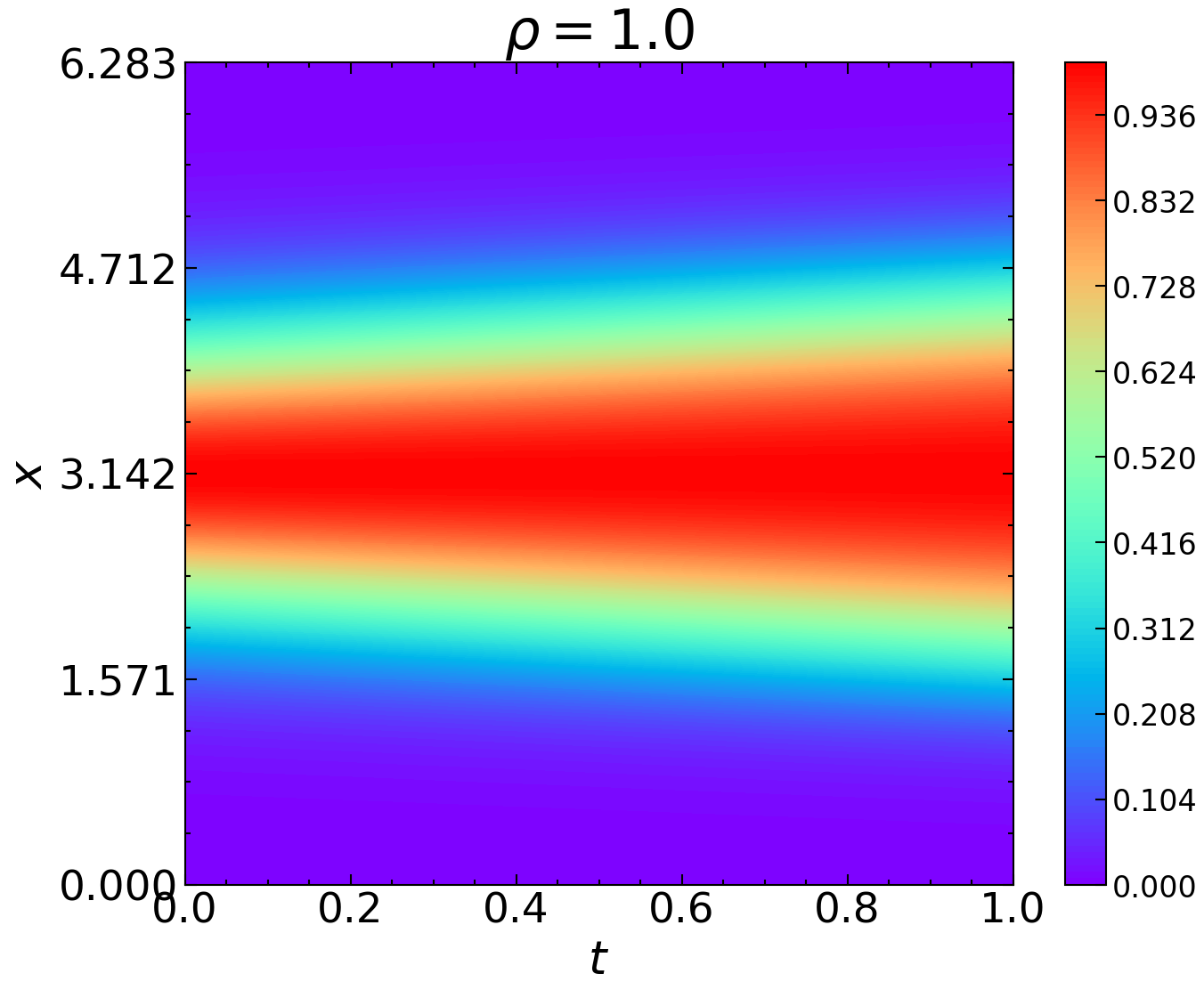}}
  \subfigure[PINN solution]{\includegraphics[width=0.24\linewidth]{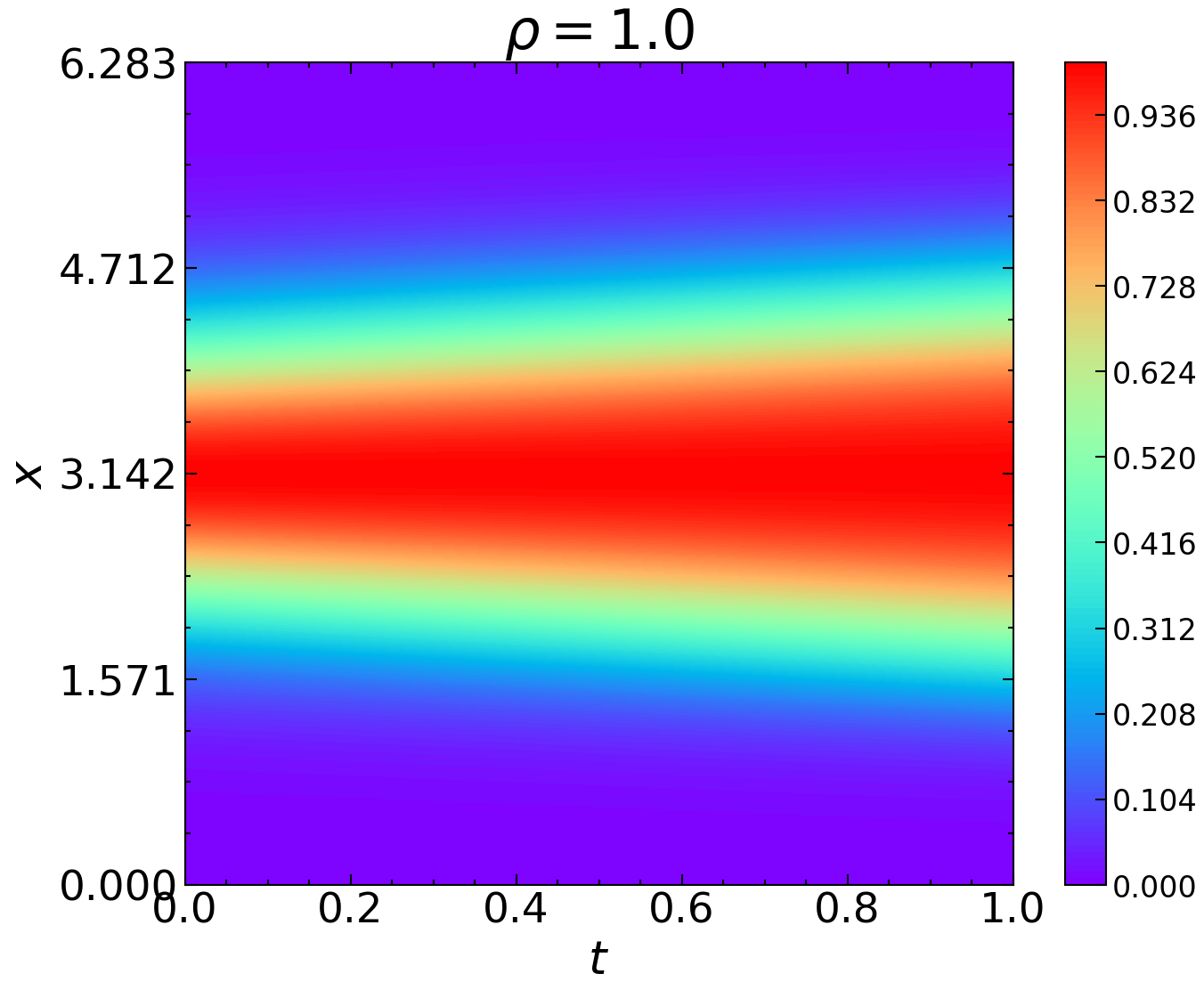}}
  \subfigure[PINN error]{\includegraphics[width=0.24\linewidth]{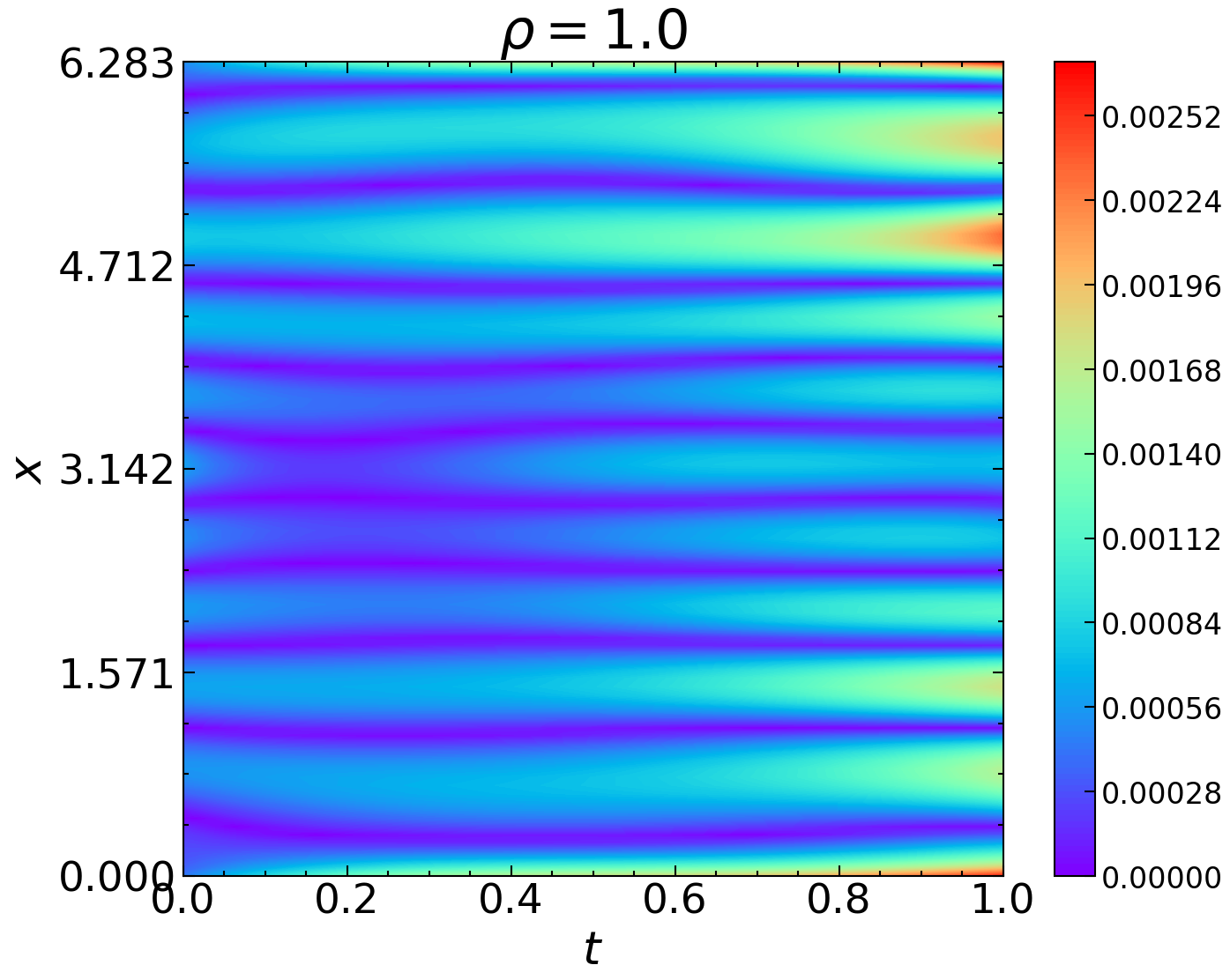}}
  \subfigure[PINN Loss]{\includegraphics[width=0.24\linewidth]{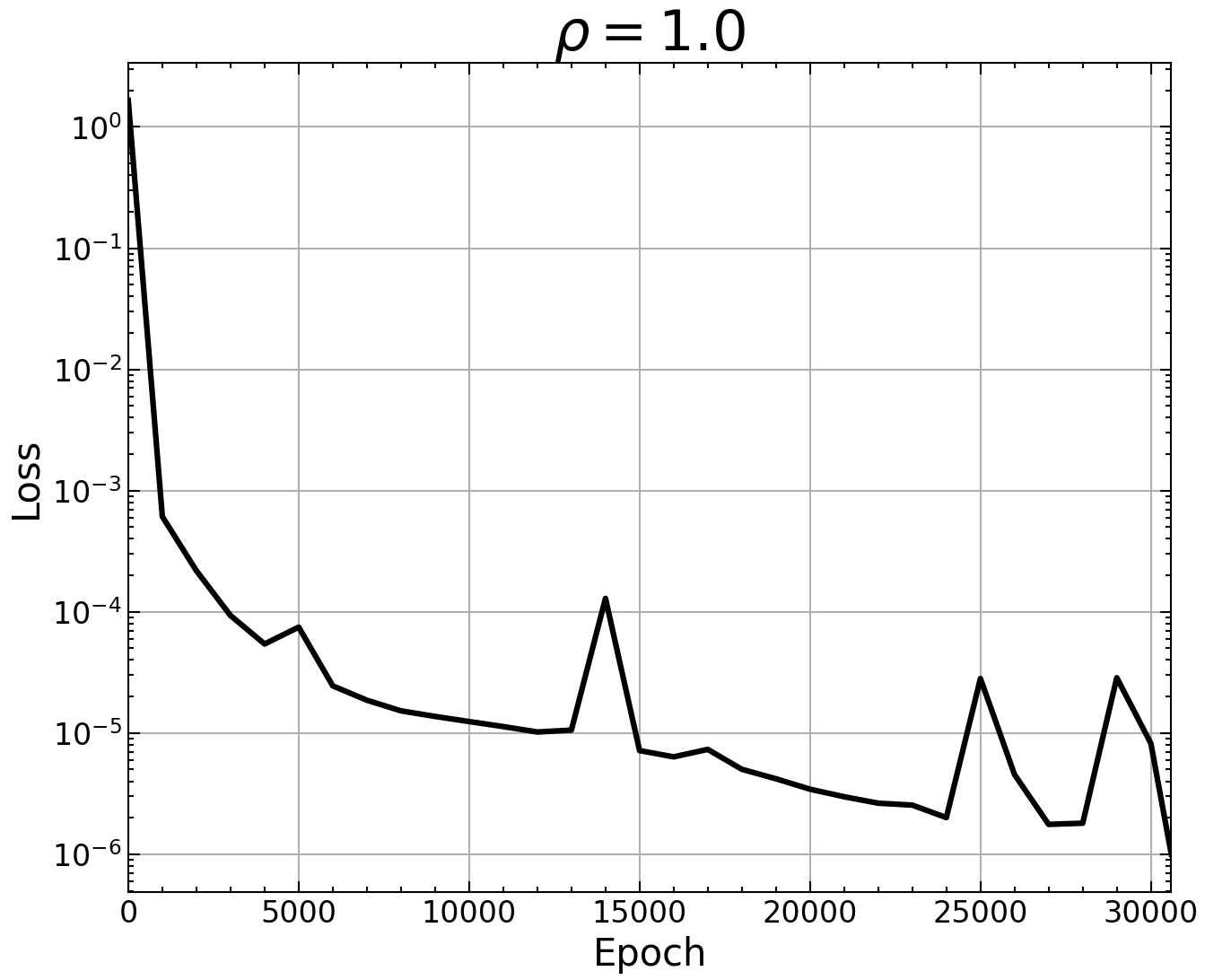}}
  \subfigure[Exact solution]{\includegraphics[width=0.24\linewidth]{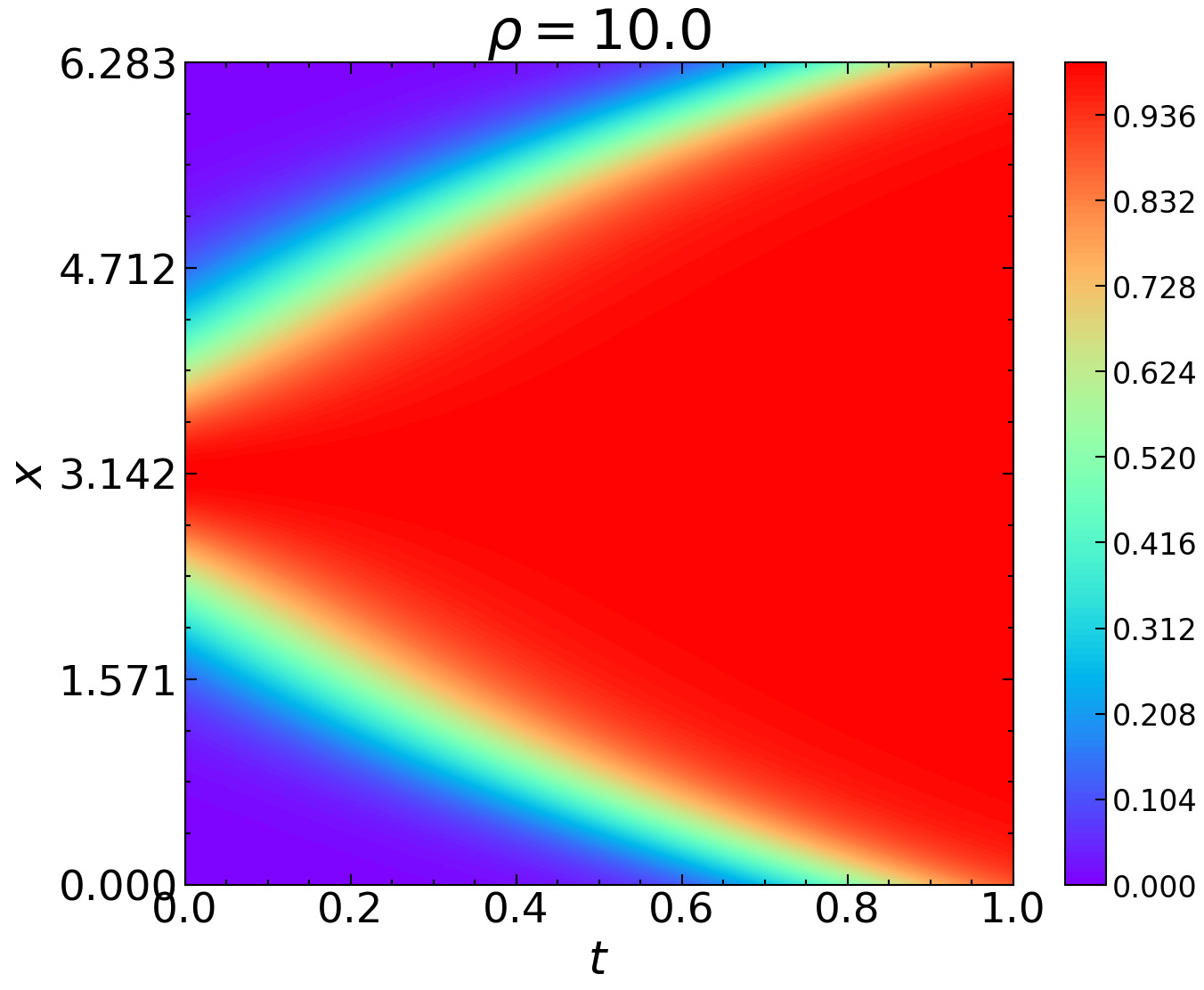}}
  \subfigure[PINN solution]{\includegraphics[width=0.24\linewidth]{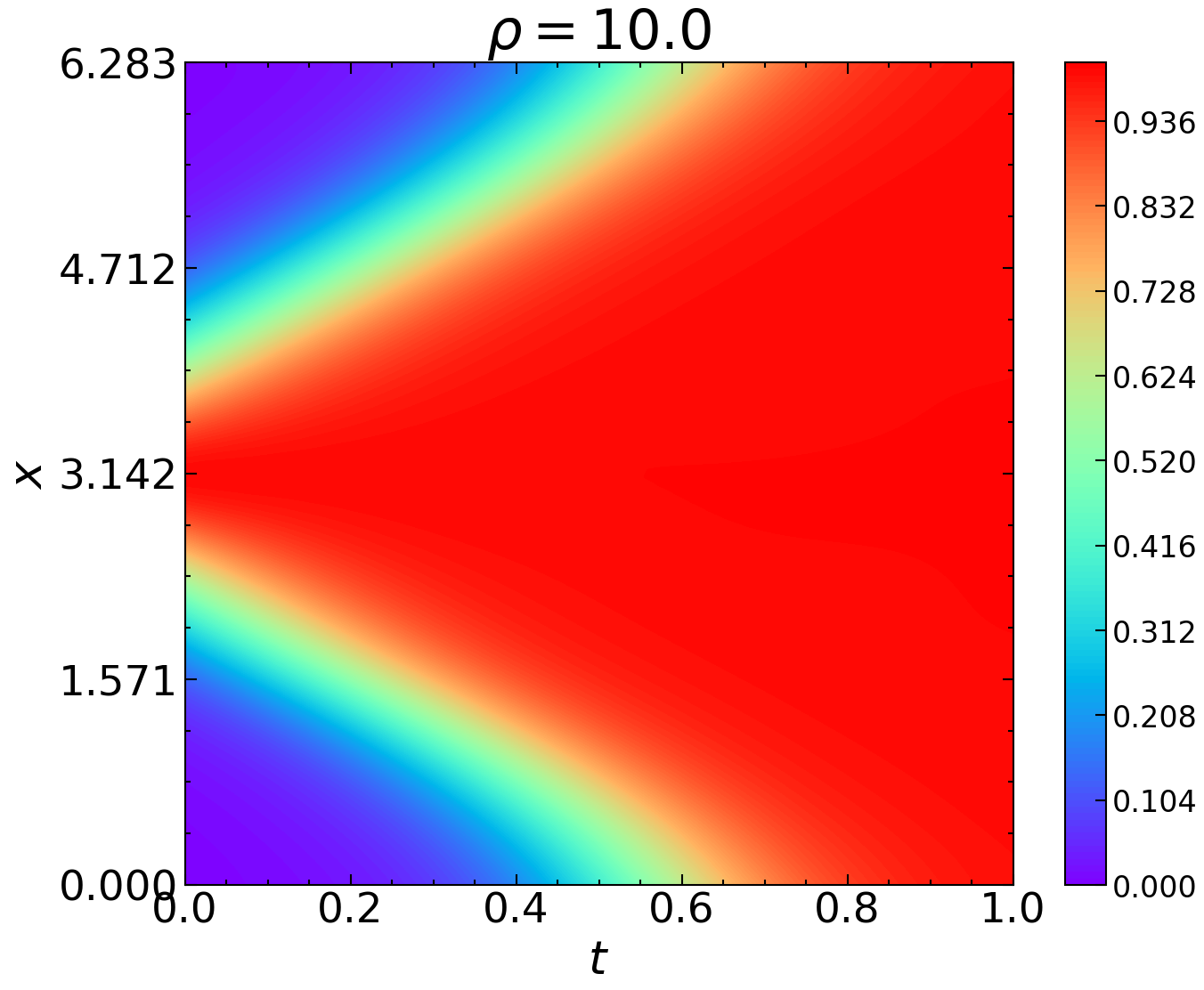}}
  \subfigure[PINN error]{\includegraphics[width=0.24\linewidth]{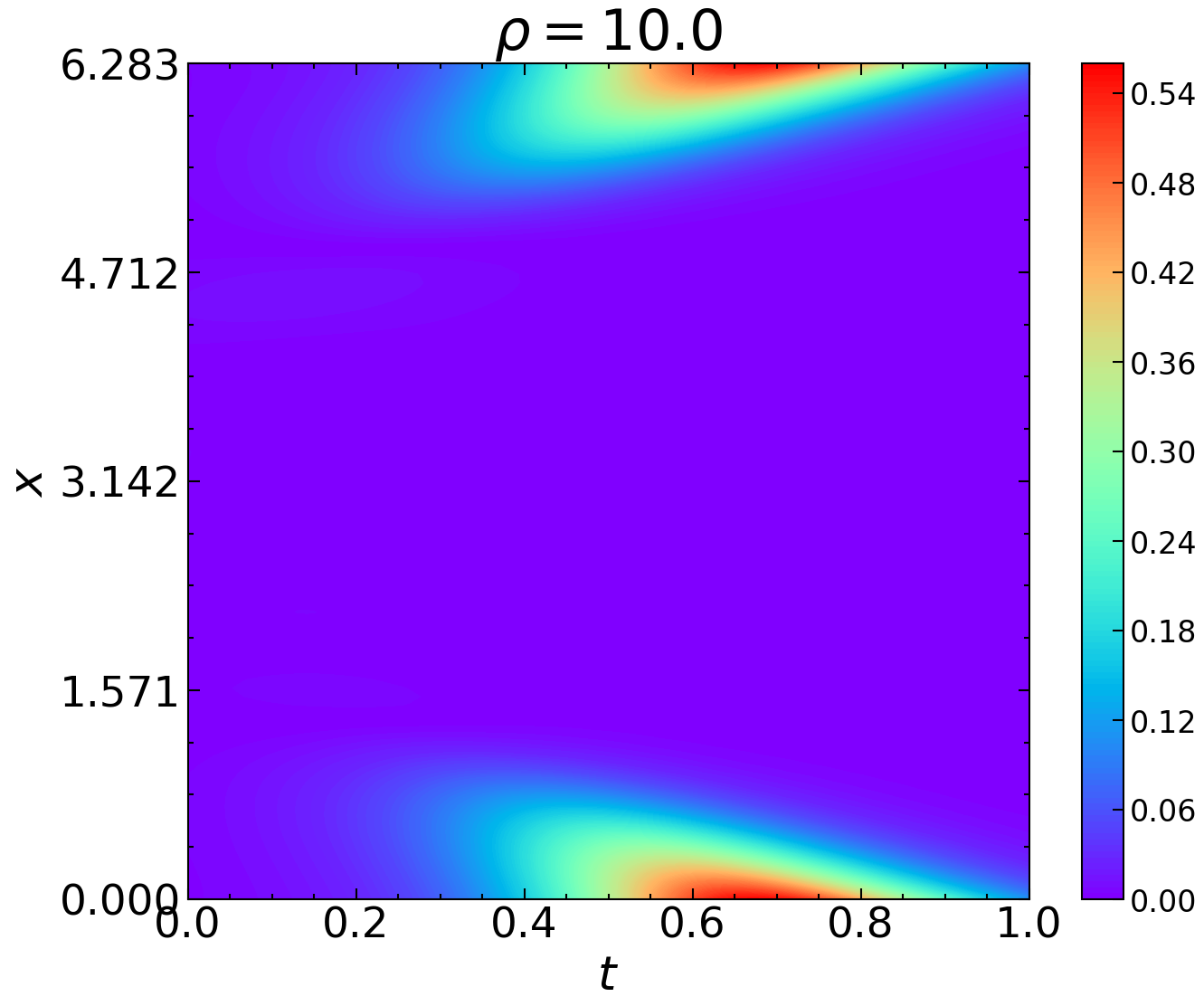}}
  \subfigure[PINN Loss]{\includegraphics[width=0.245\linewidth]{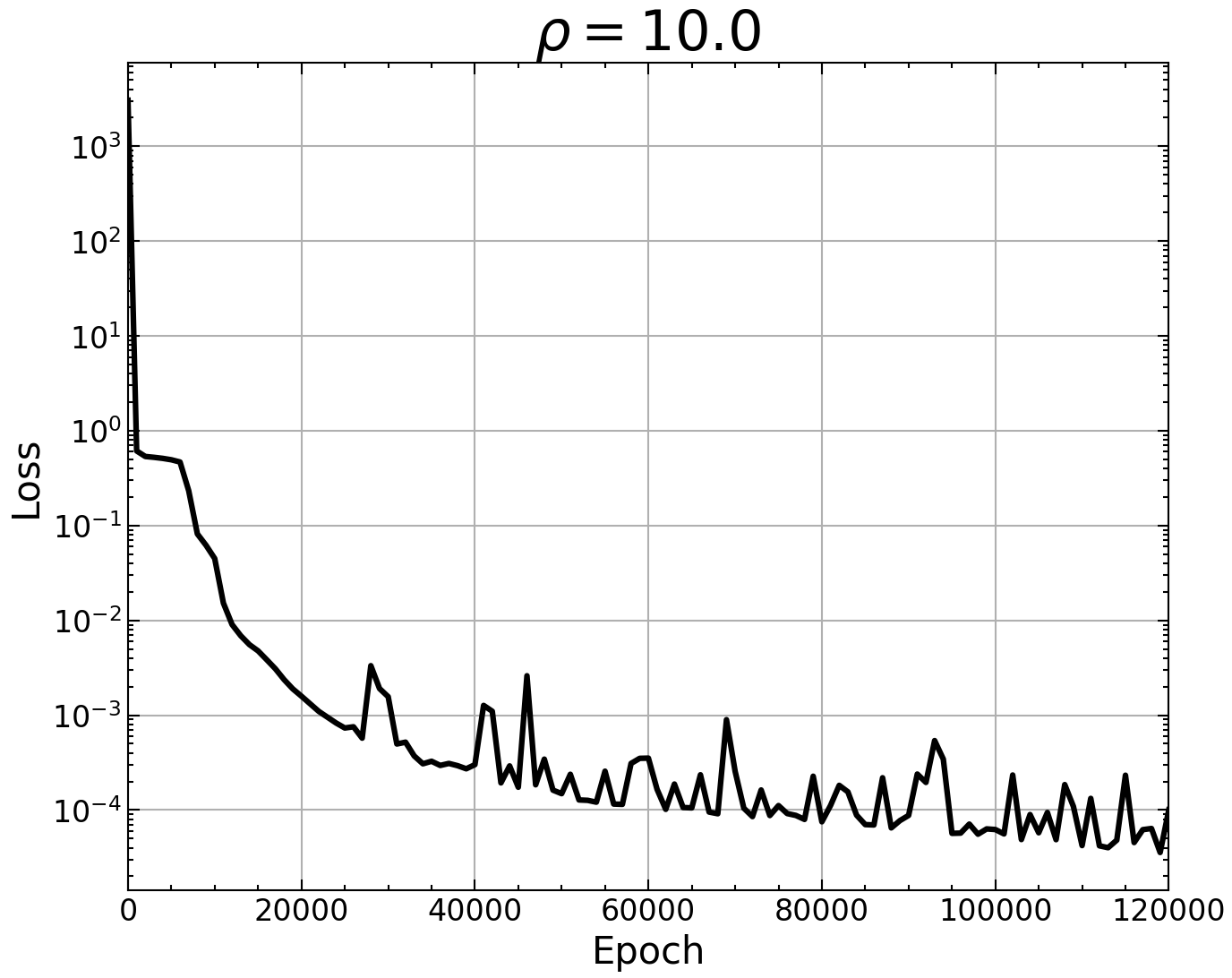}}
  \caption{Results for \cref{sssec:results-ppde-reaction}: PINN snapshots for $\rho = 1.0$ (top) and $10.0$ (bottom).}
  \label{fig:R-Full-PINN1}
\end{figure}

\noindent {\bf PINN results:}
The PINN we employ is a fully connected neural network with layer widths [2, 20, 20, 20, 1] with $\lambda \equiv 1$ in $\mathcal{L}_{\rm int}(u)$ of the PINN loss function \eqref{eq:continousloss}. The training samples are randomly selected, but we use an optimized activation function introduced in \cite{zhao2023pinnsformer}:
\begin{align}\label{eq:waveact}
{\rm WaveAct}(x)=w_1*\sin(x)+w_2*\cos(x),
\end{align}
where $w_1$ and $w_2$ are trainable hyperparameters.  Numerical experiments demonstrate that such an activation function is more suitable for solving reaction equations compared to traditional activation functions. Figure \ref{fig:R-Full-PINN1} presents the accuracy and loss reduction of the full PINN for parameters $\rho=1$ and $10$.
This problem is difficult to numerically solve for large values of the parameter $\rho$. One can observe this difficultly in Figure \ref{fig:R-Full-PINN1} where $\rho = 10$ requires substantially more epochs to successfully train; this is due to the more complex nature of the phase transition layer for larger $\rho$.

\begin{figure}[htbp]
  \centering
  \includegraphics[scale=0.2]{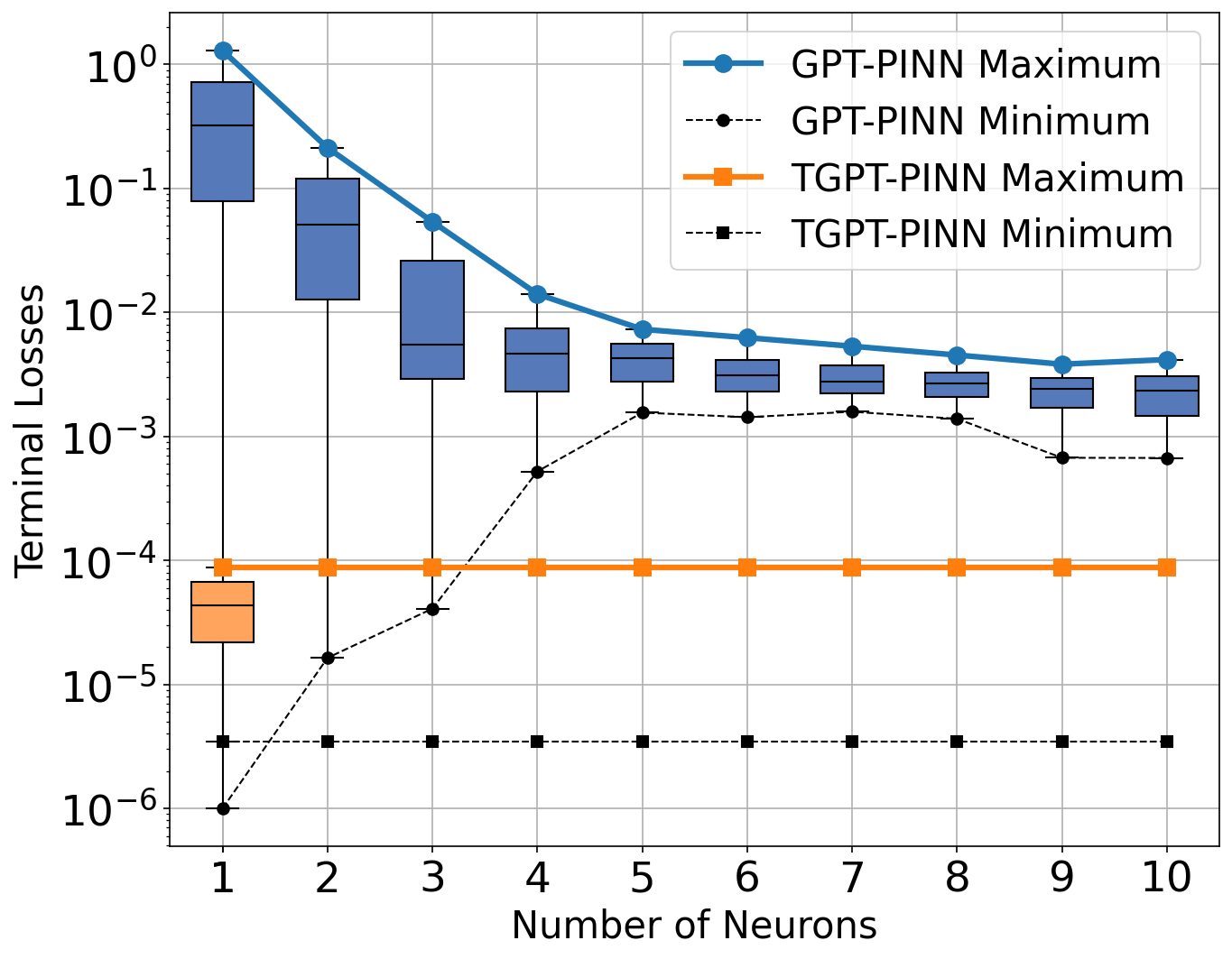}
  \includegraphics[scale=0.205]{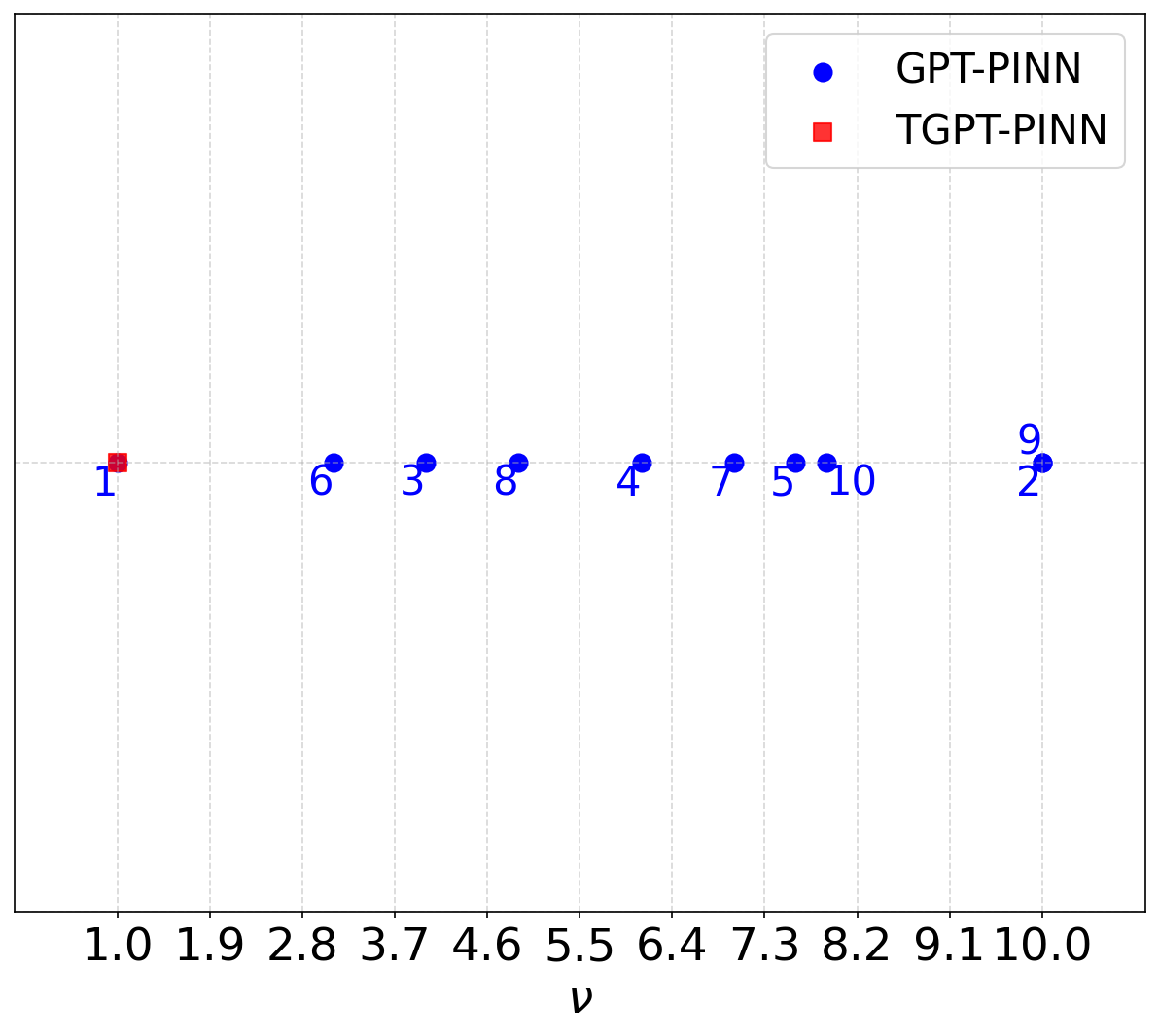}
  \includegraphics[scale=0.2]{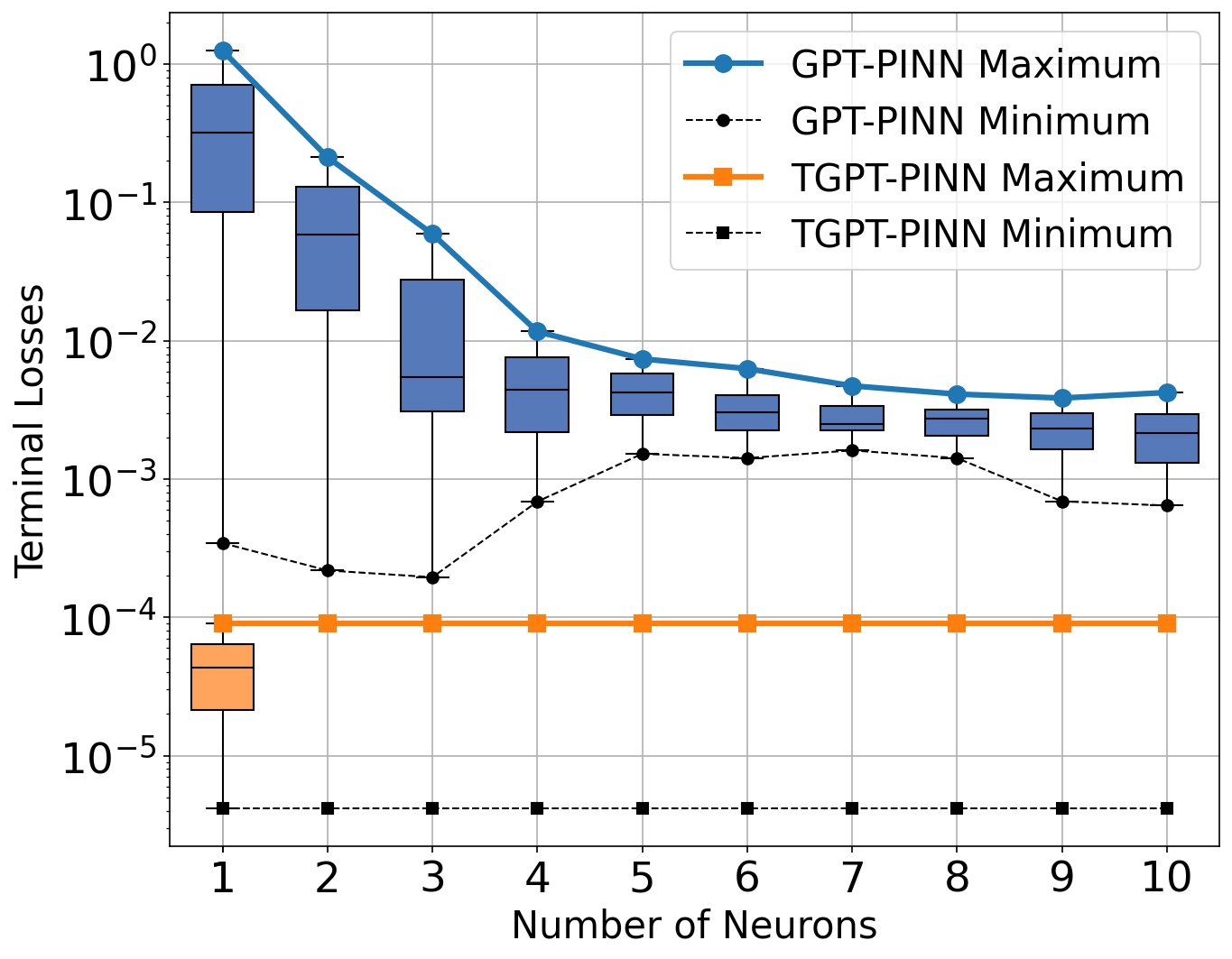}
\put(-270,18){\makebox(0,0){{(a)}}}
\put(-145,18){\makebox(0,0){{(b)}}}
\put(-10,18){\makebox(0,0){{(c)}}}
  \caption{Results for Section \ref{sssec:results-ppde-reaction}: Comparison ofthe  GPT-PINN and TGPT-PINN for the reaction equation in terms of histories of convergence (a, c) and parameter selection (b).}
  \label{fig:R-GPT-PINN}
\end{figure}

\begin{figure}[htbp]
  \centering
  \subfigure[Exact solution]{\includegraphics[width=0.24\linewidth]{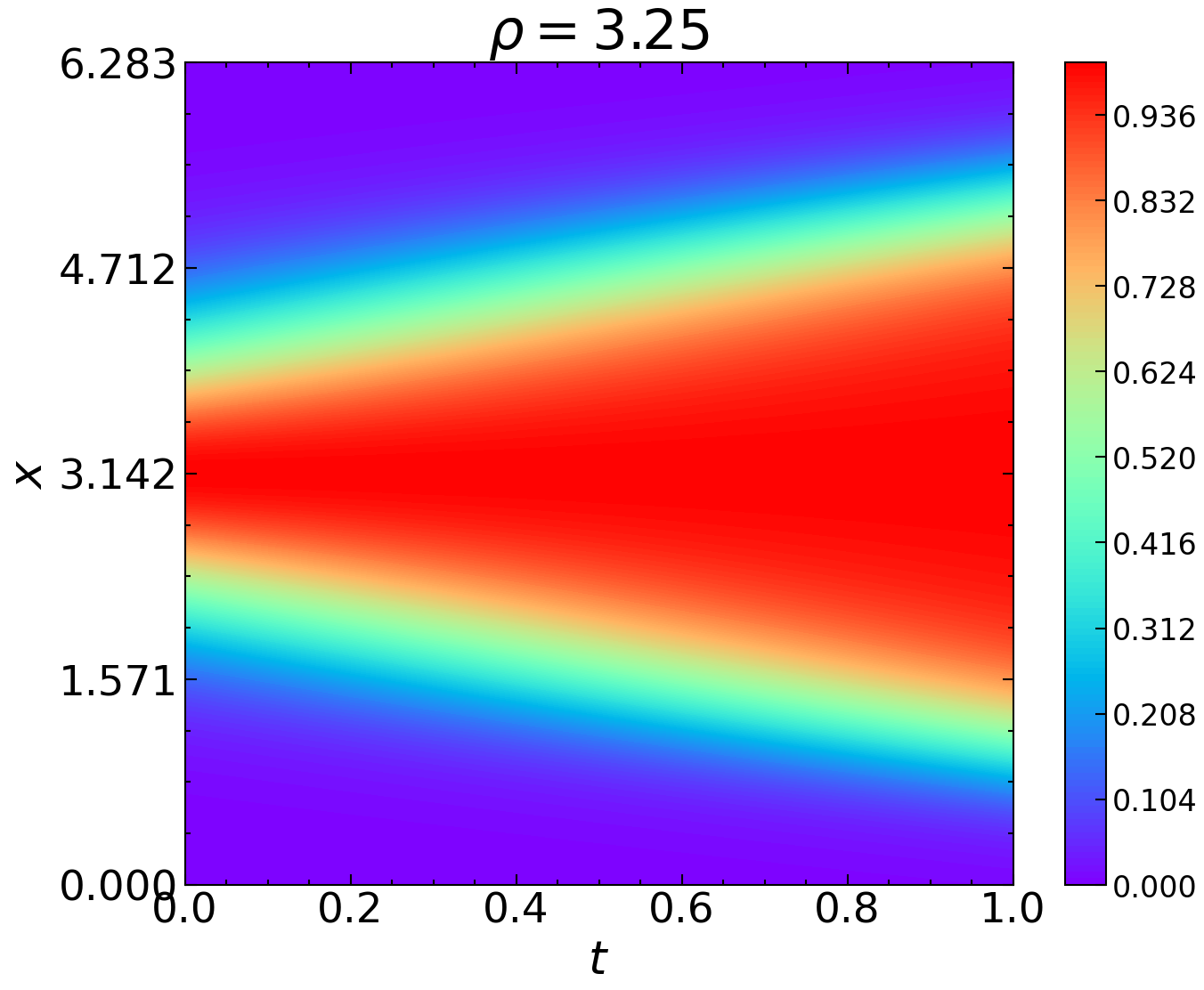}}
  \subfigure[GPT-PINN solution]{\includegraphics[width=0.24\linewidth]{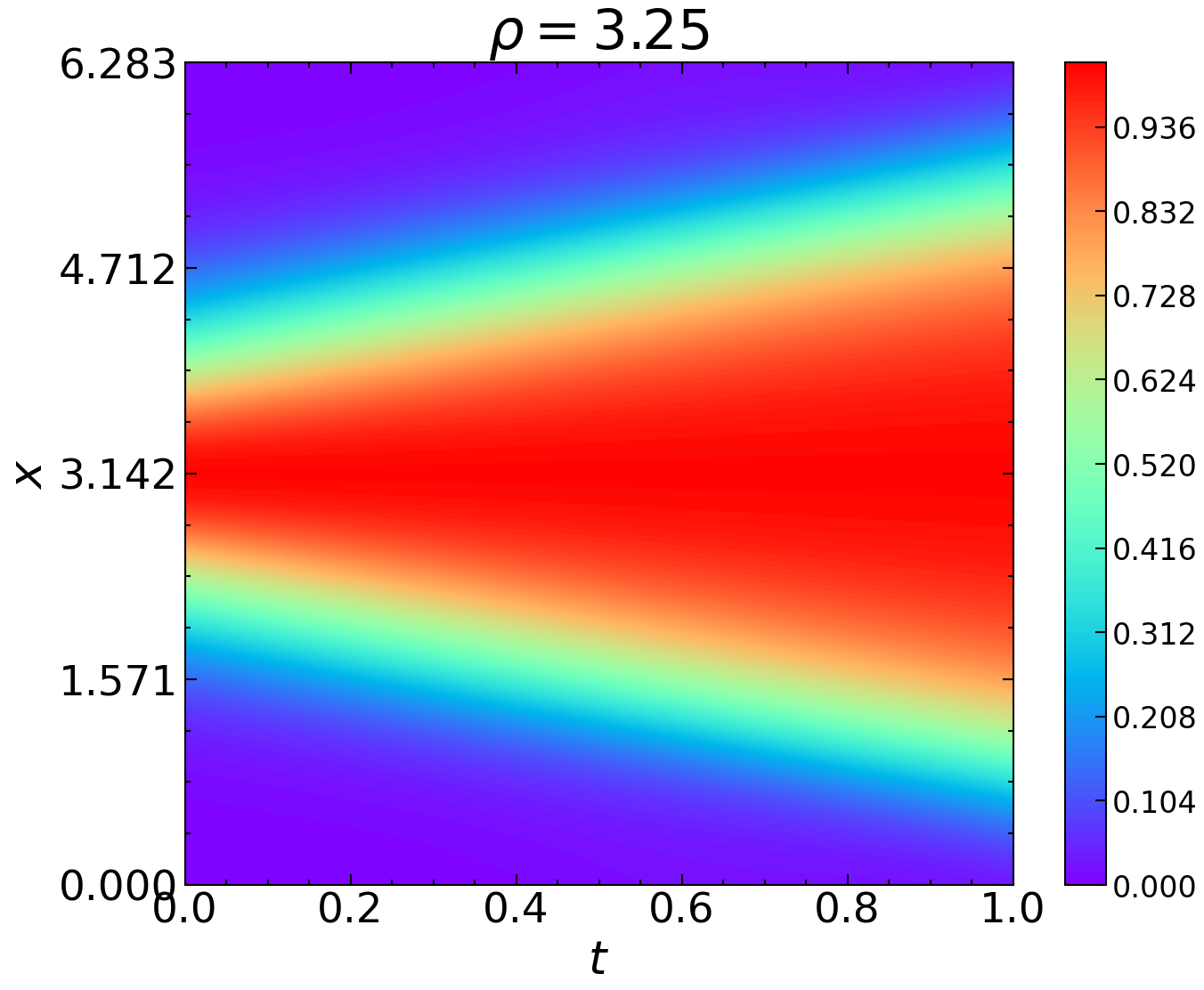}}
  \subfigure[GPT-PINN error]{\includegraphics[width=0.24\linewidth]{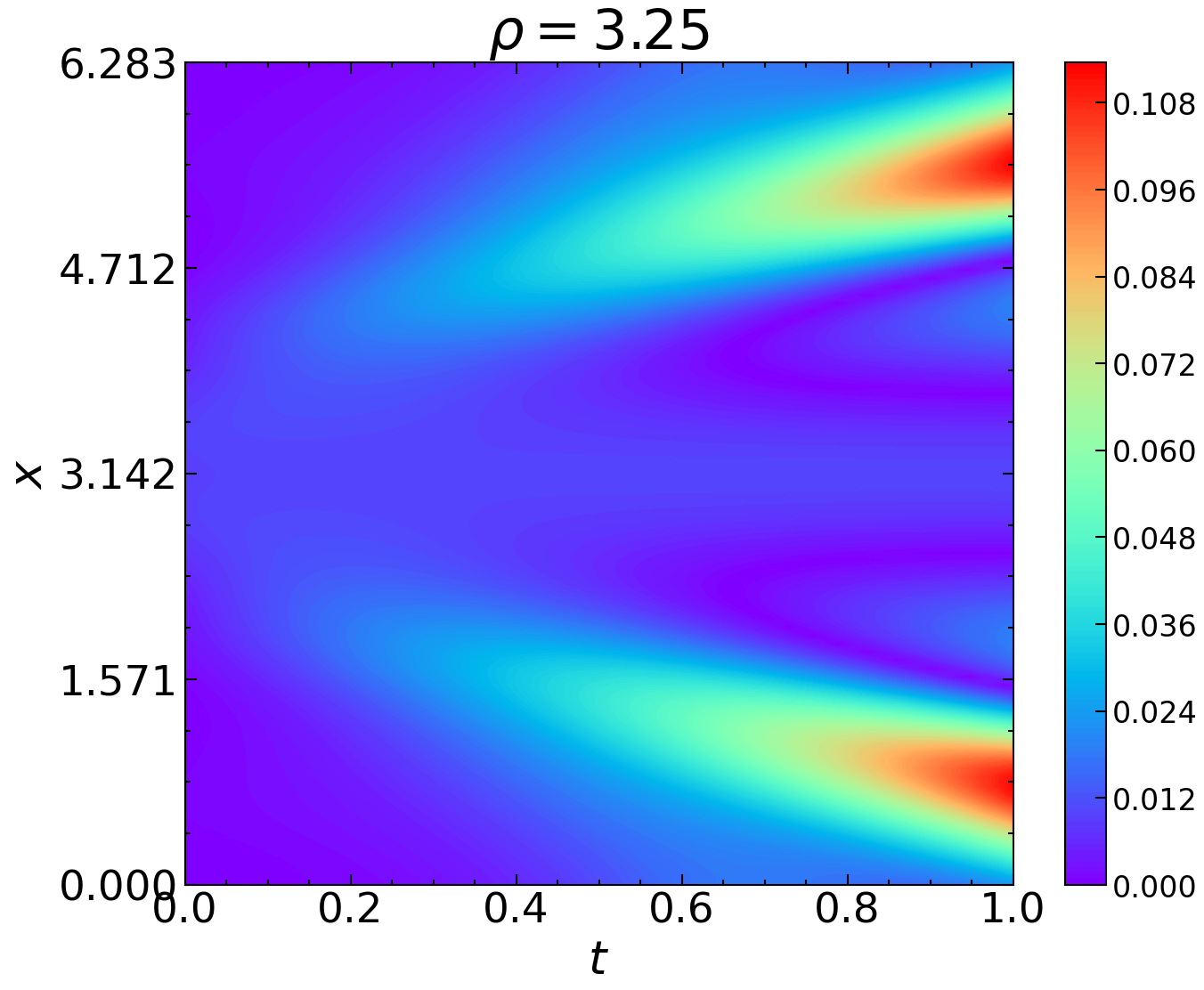}}
  \subfigure[GPT-PINN Loss]{\includegraphics[width=0.24\linewidth]{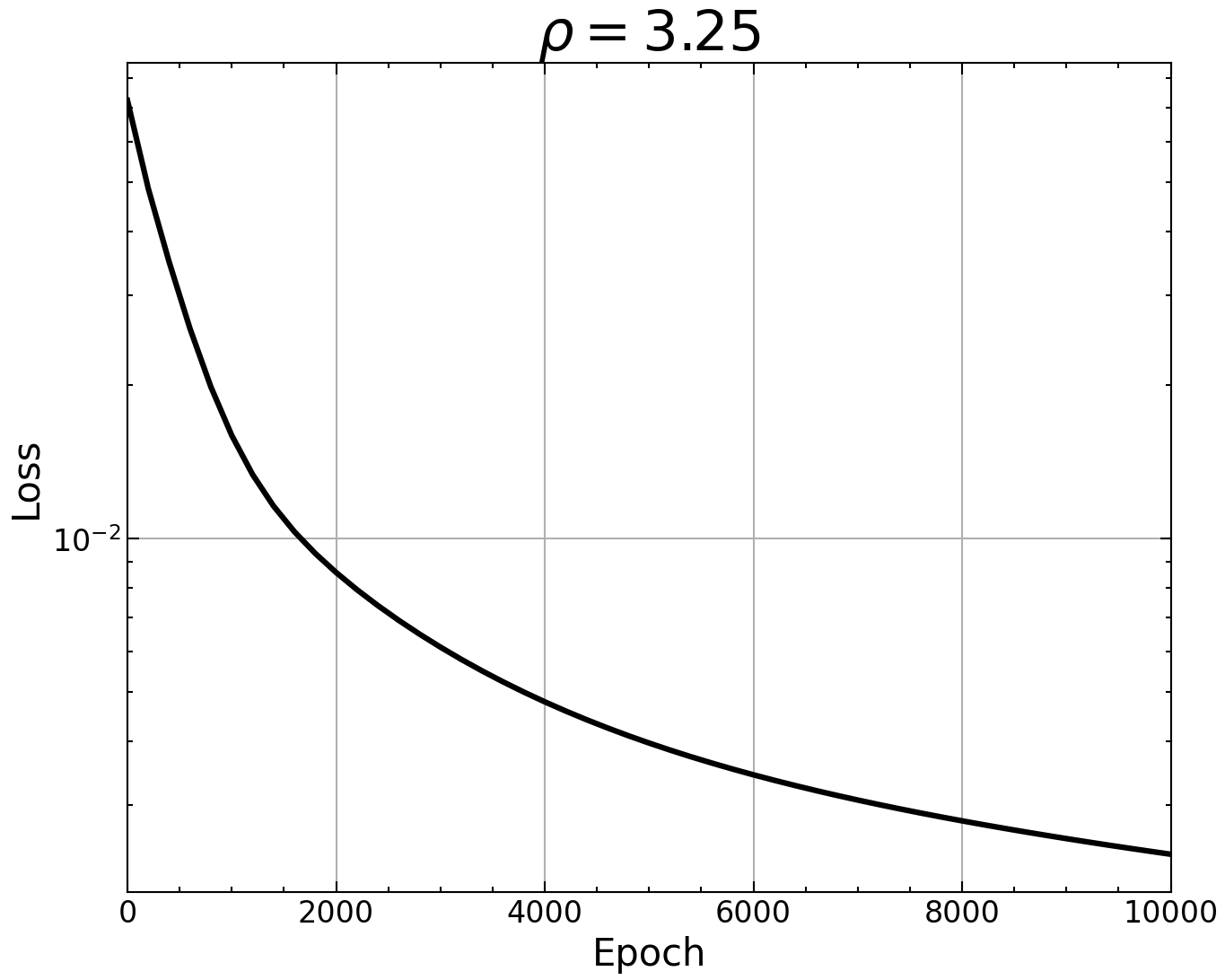}}

  \subfigure[Exact solution]{\includegraphics[width=0.24\linewidth]{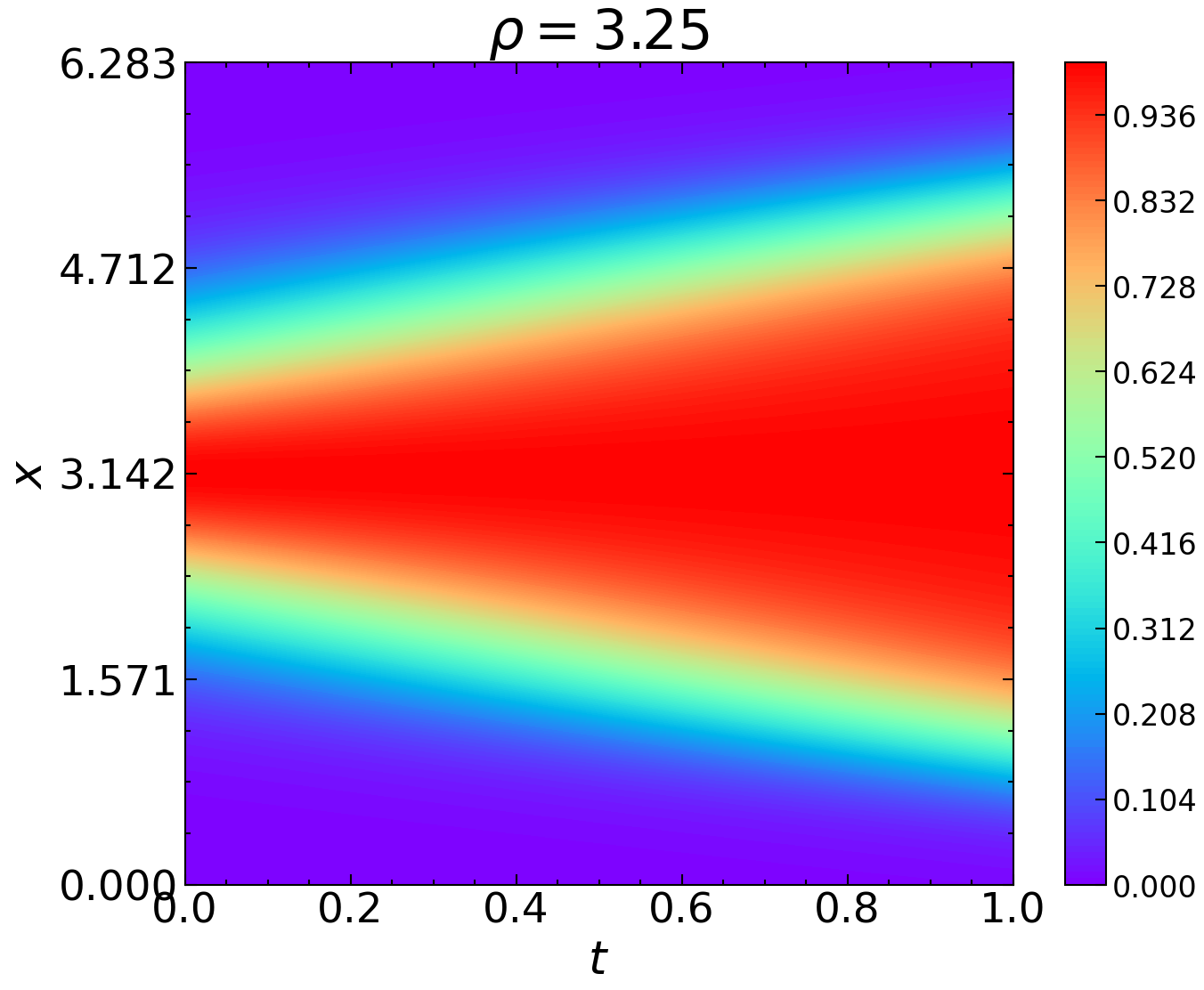}}
  \subfigure[TGPT-PINN solution]{\includegraphics[width=0.24\linewidth]{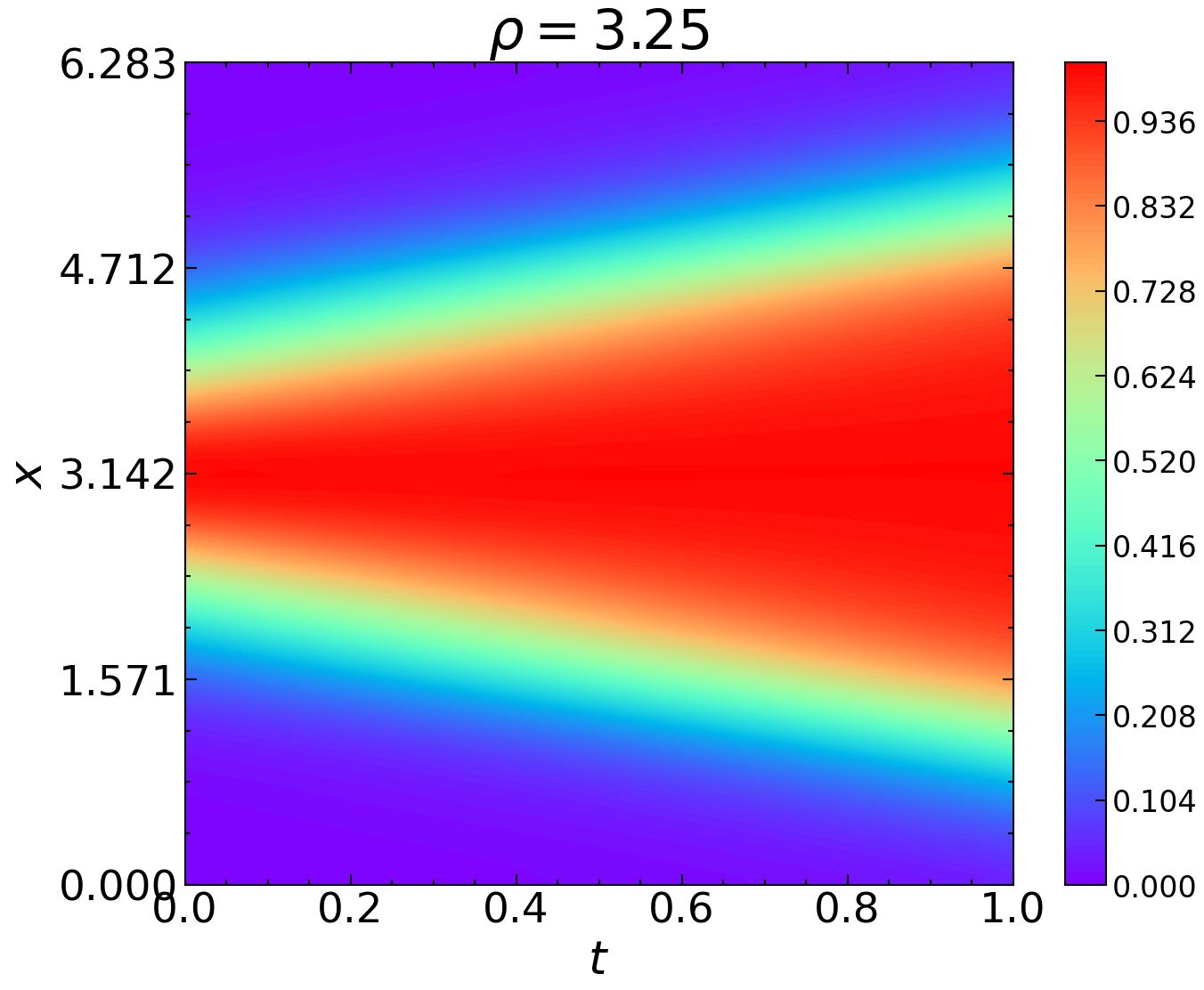}}
  \subfigure[TGPT-PINN error]{\includegraphics[width=0.24\linewidth]{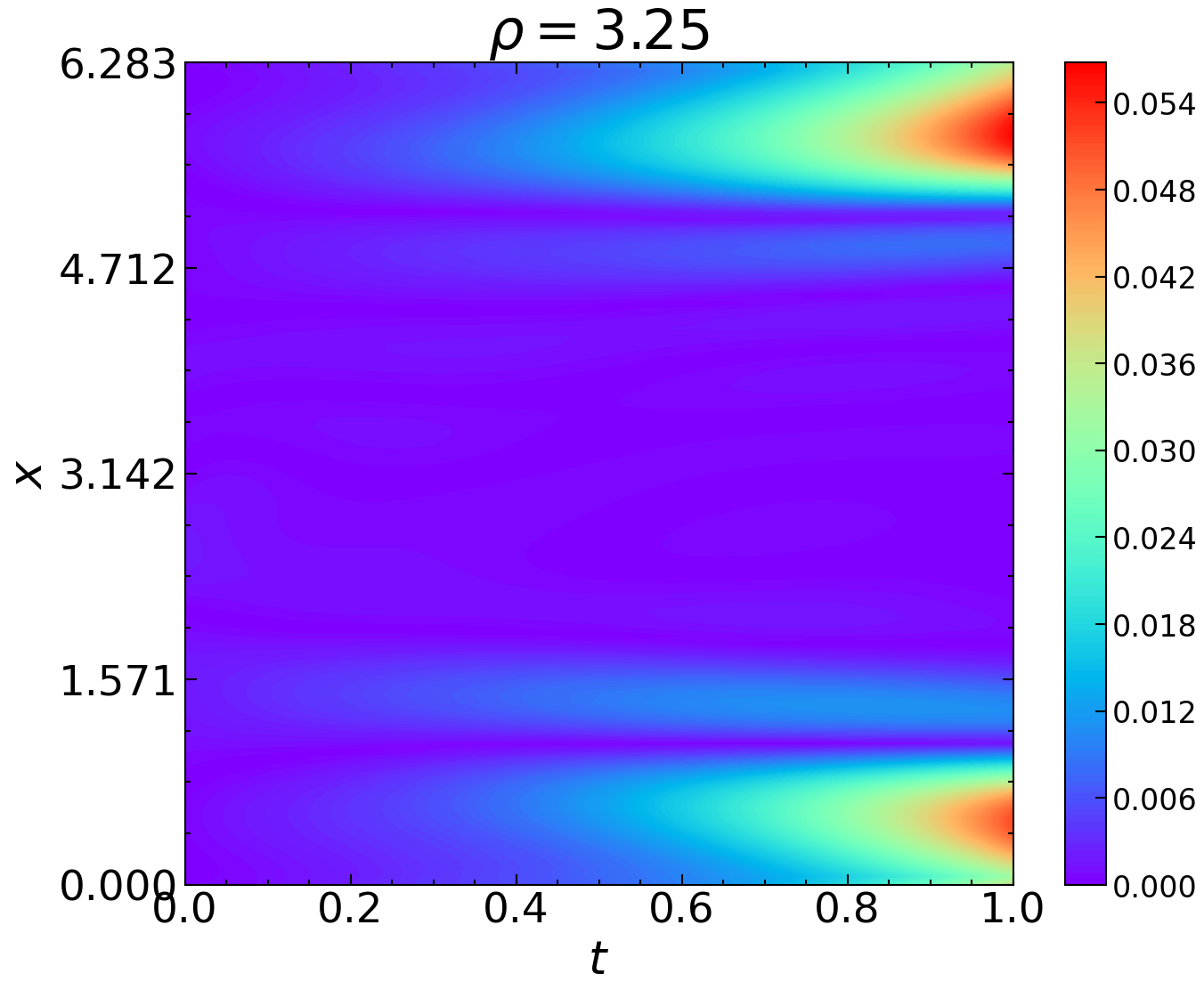}}
  \subfigure[TGPT-PINN Loss]{\includegraphics[width=0.24\linewidth]{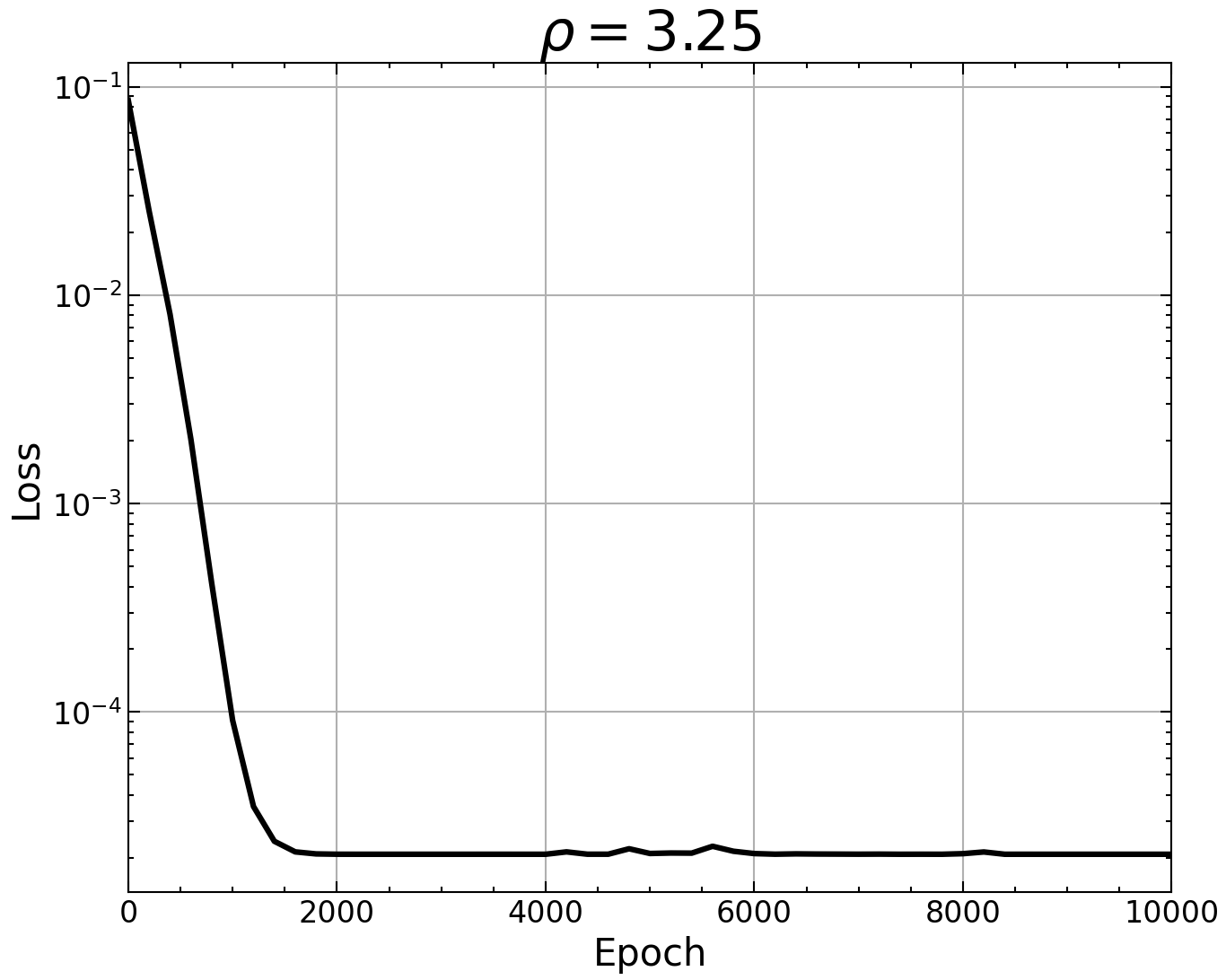}}
  \caption{Results for Section \ref{sssec:results-ppde-reaction}: GPT-PINN (top, with 10 neurons) and TGPT-PINN (bottom, with 1 neuron) results for $\rho = 3.25$.}
  \label{fig:R-gpt-pinn1}
\end{figure}

\begin{figure}[htbp]
  \centering
  \subfigure[Exact solution]{\includegraphics[width=0.24\linewidth]{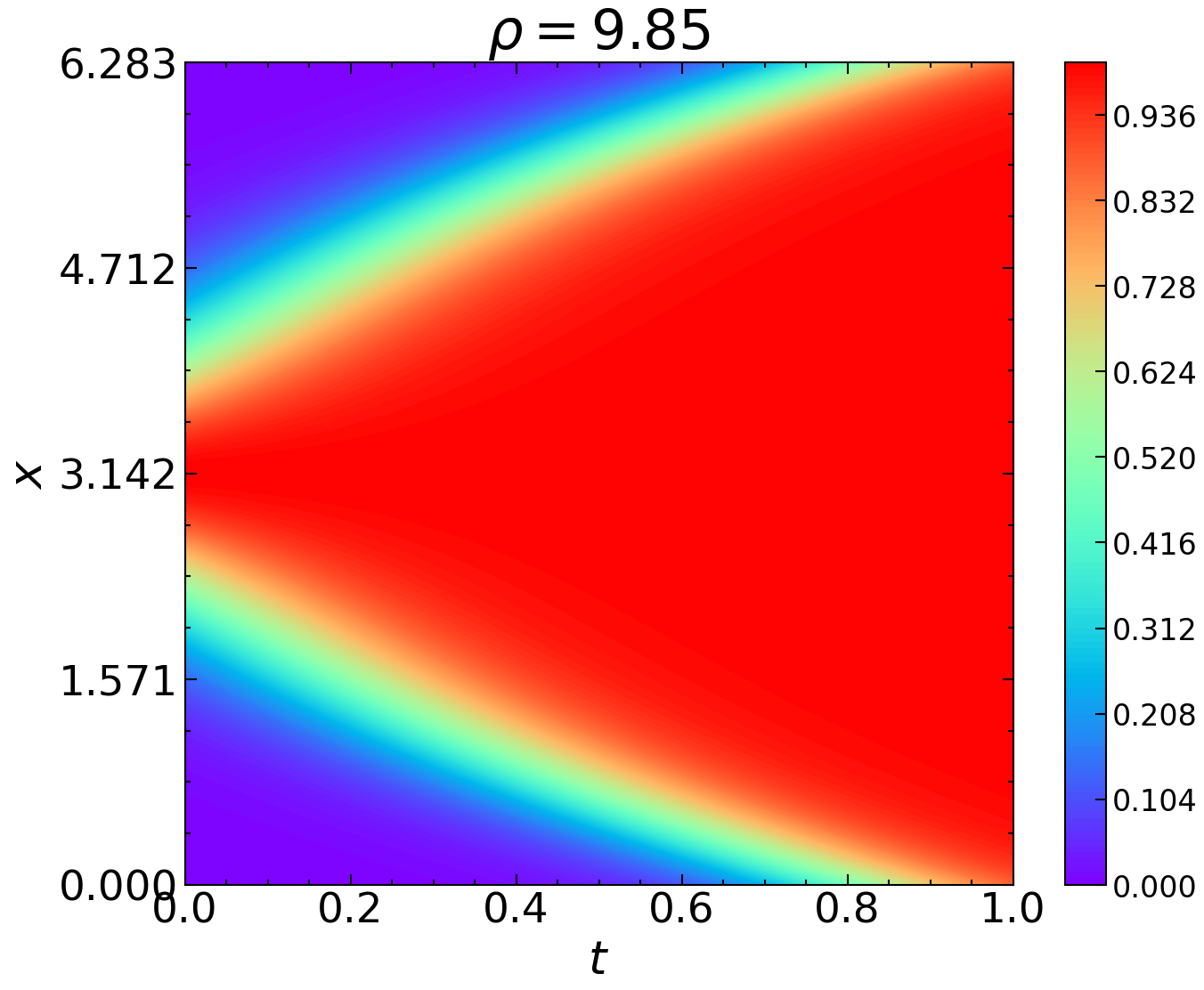}}
  \subfigure[GPT-PINN solution]{\includegraphics[width=0.24\linewidth]{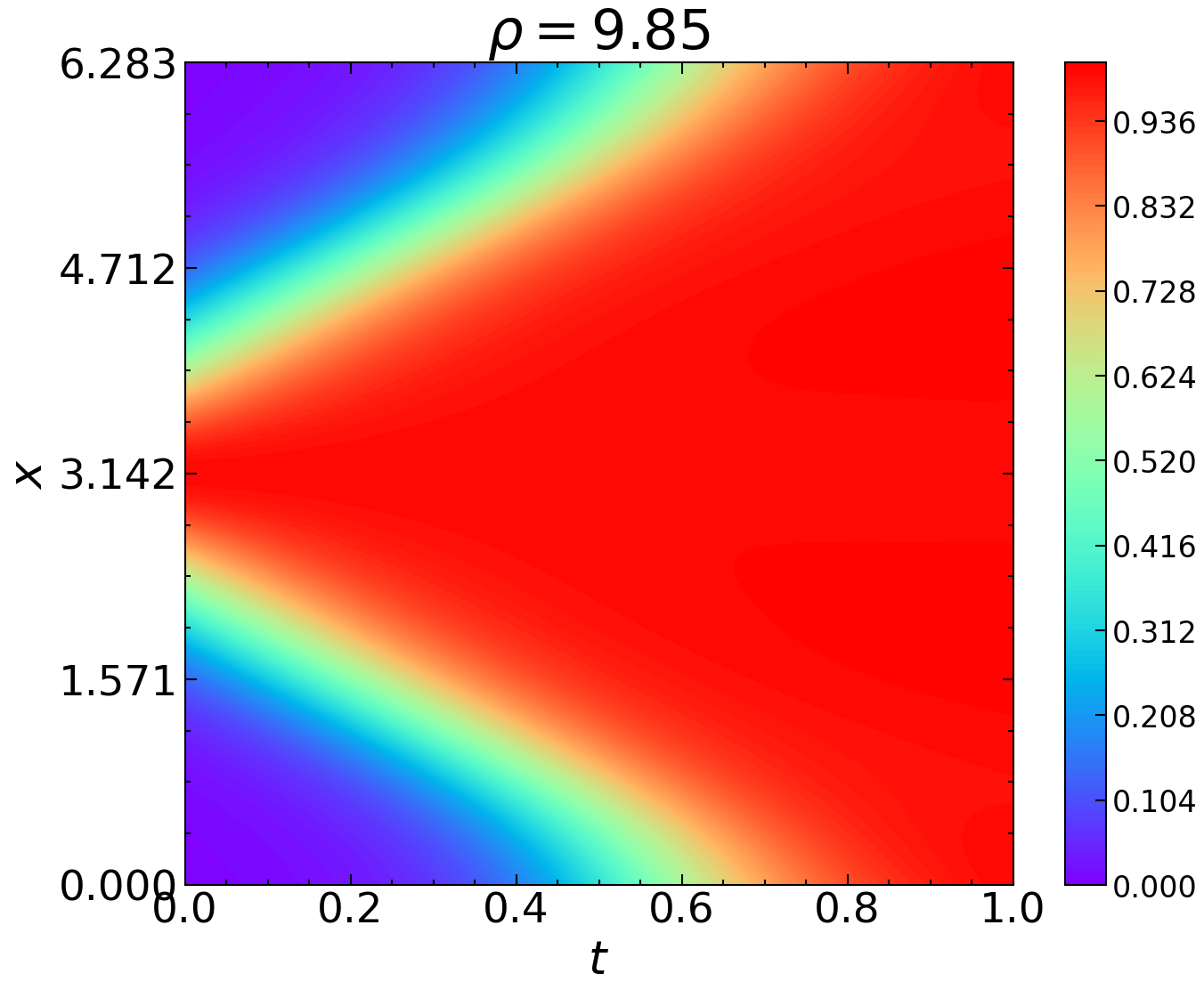}}
  \subfigure[GPT-PINN error]{\includegraphics[width=0.24\linewidth]{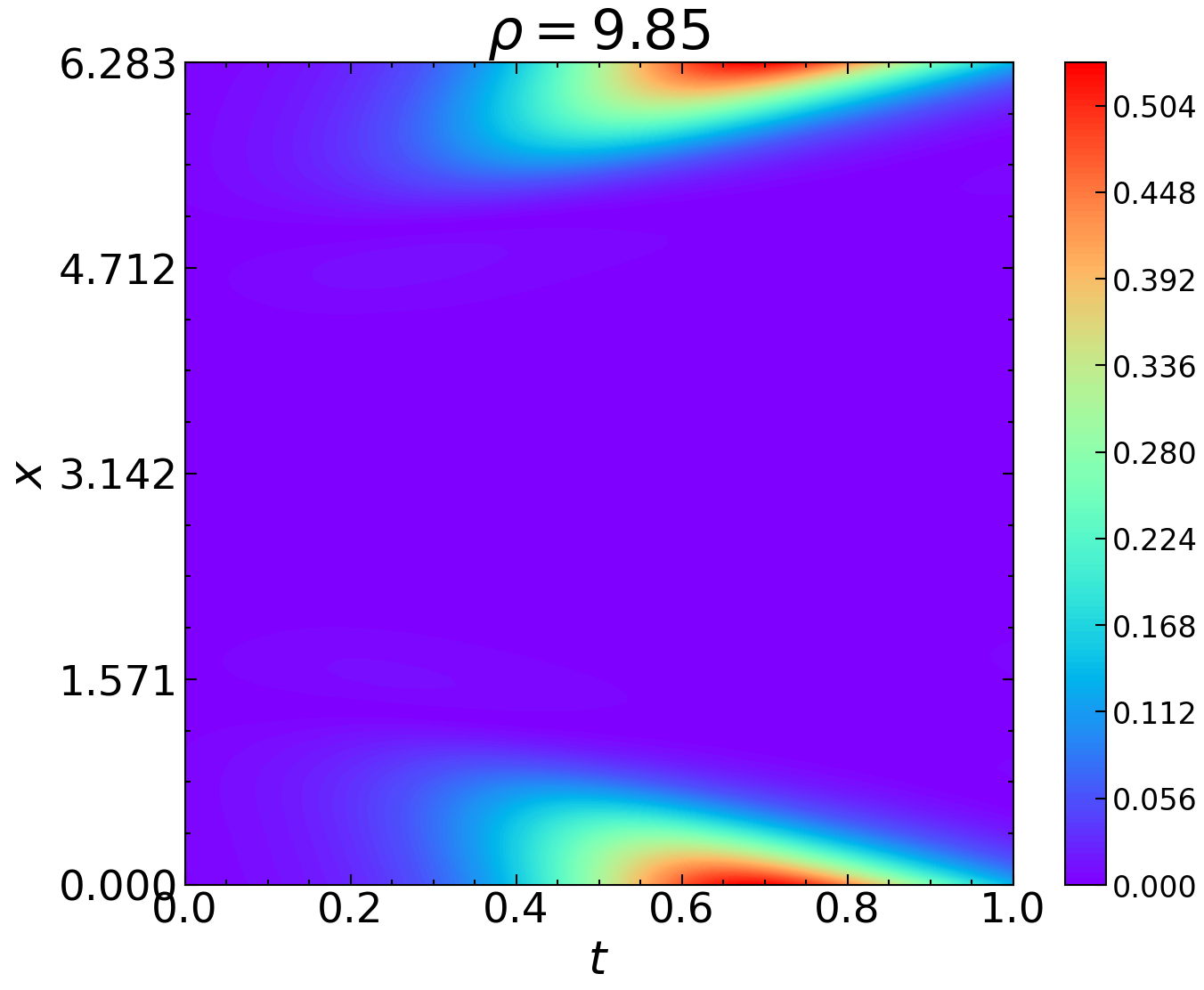}}
  \subfigure[GPT-PINN Loss]{\includegraphics[width=0.24\linewidth]{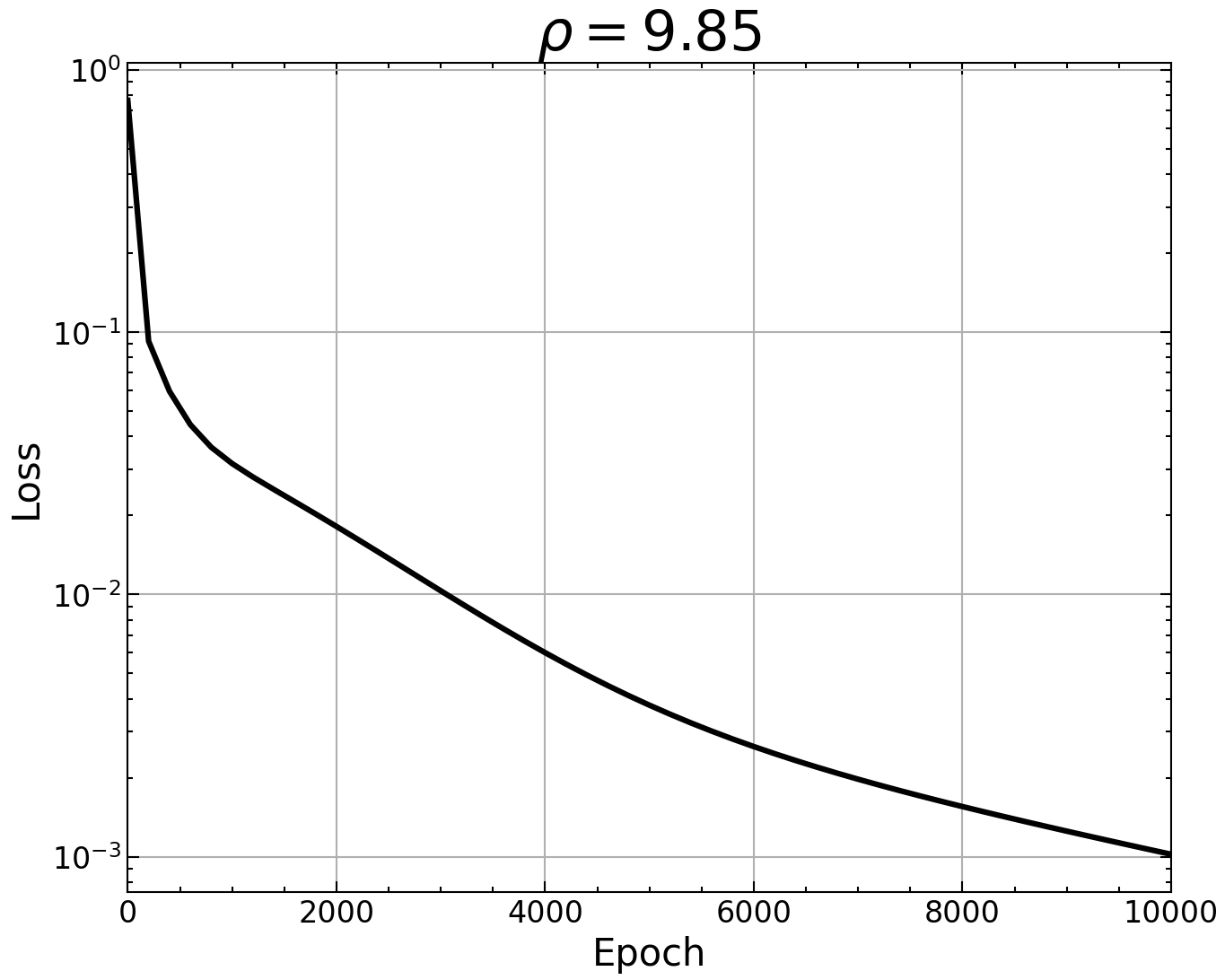}}

  \subfigure[Exact solution]{\includegraphics[width=0.24\linewidth]{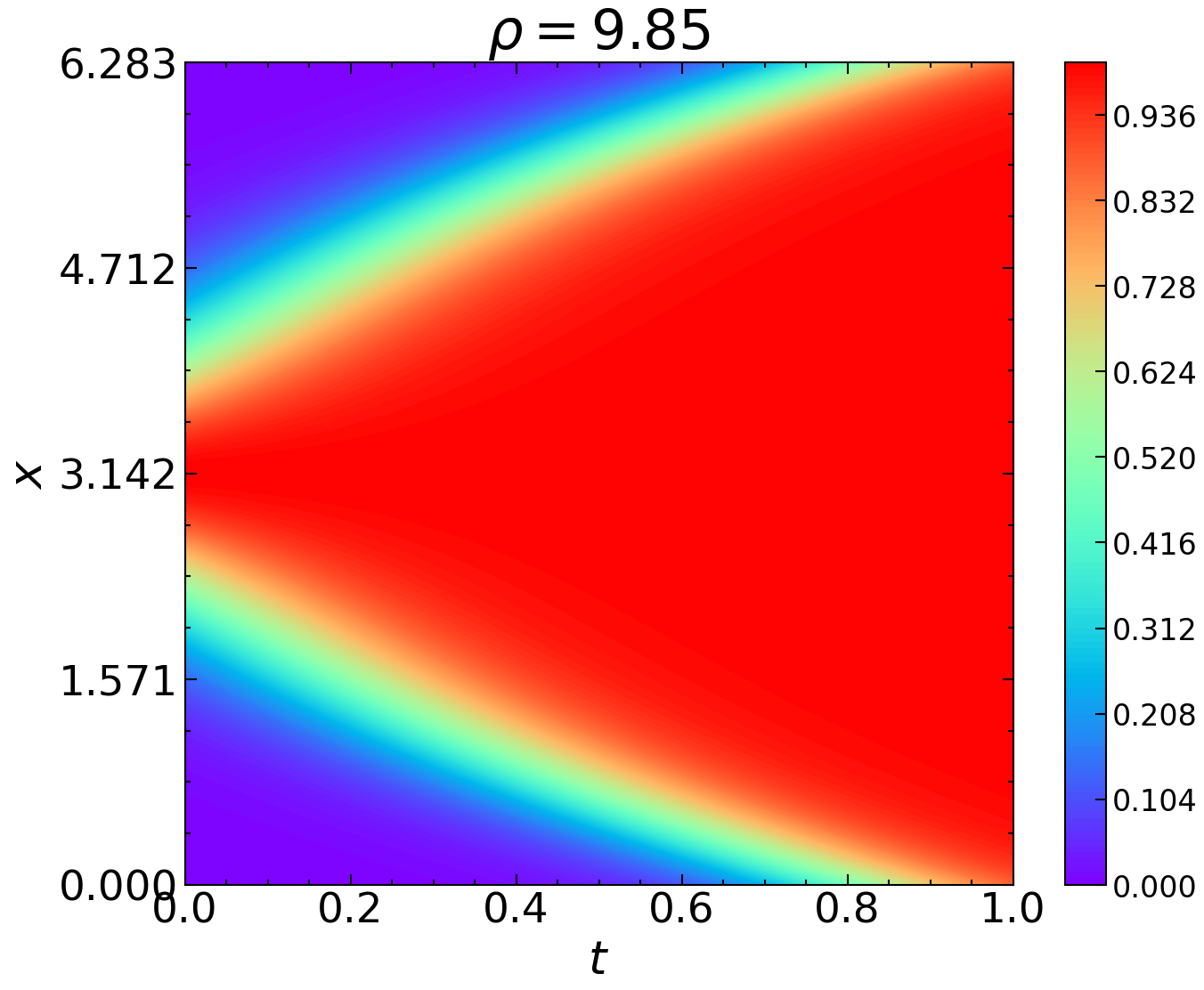}}
  \subfigure[TGPT-PINN solution]{\includegraphics[width=0.24\linewidth]{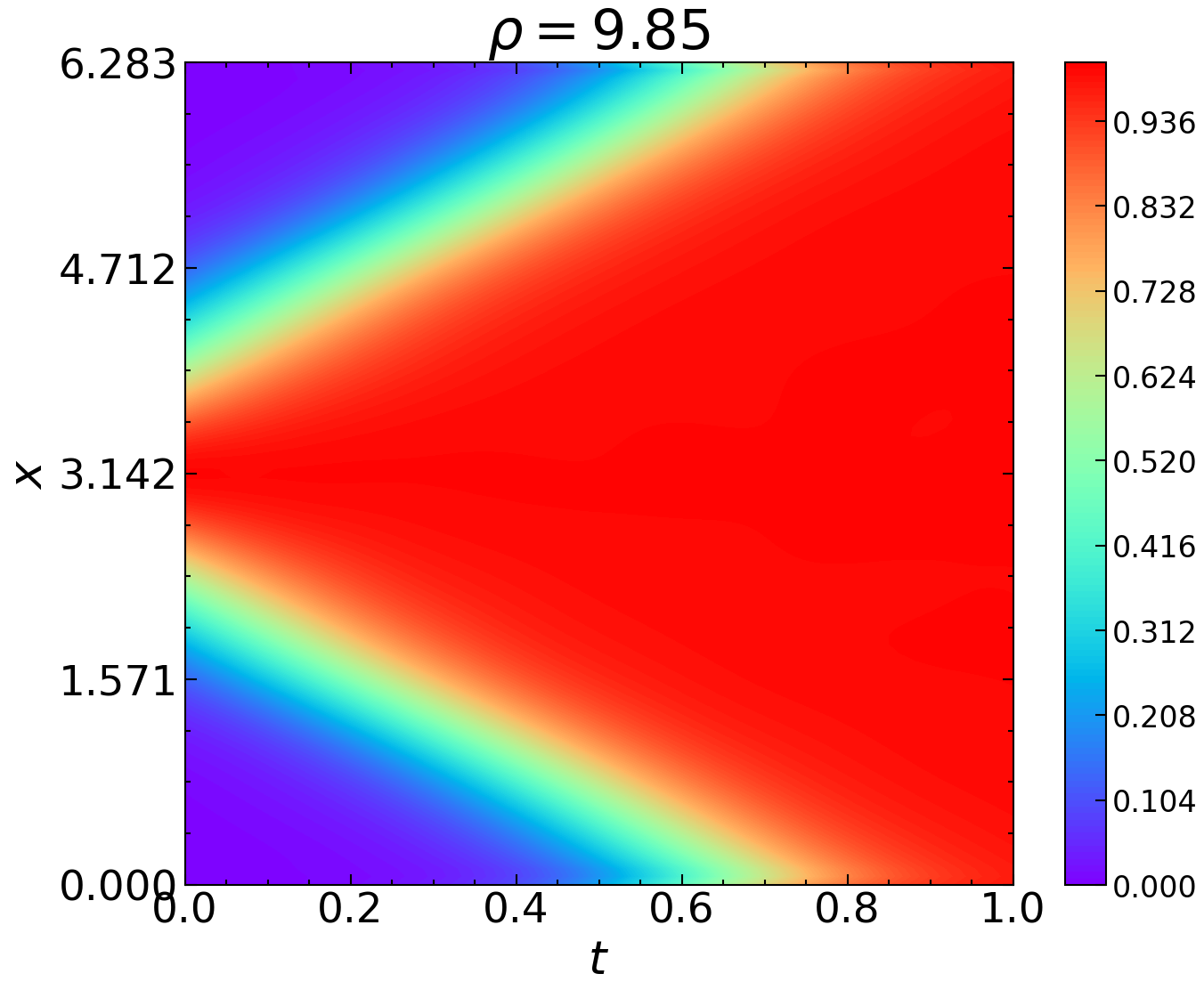}}
  \subfigure[TGPT-PINN error]{\includegraphics[width=0.24\linewidth]{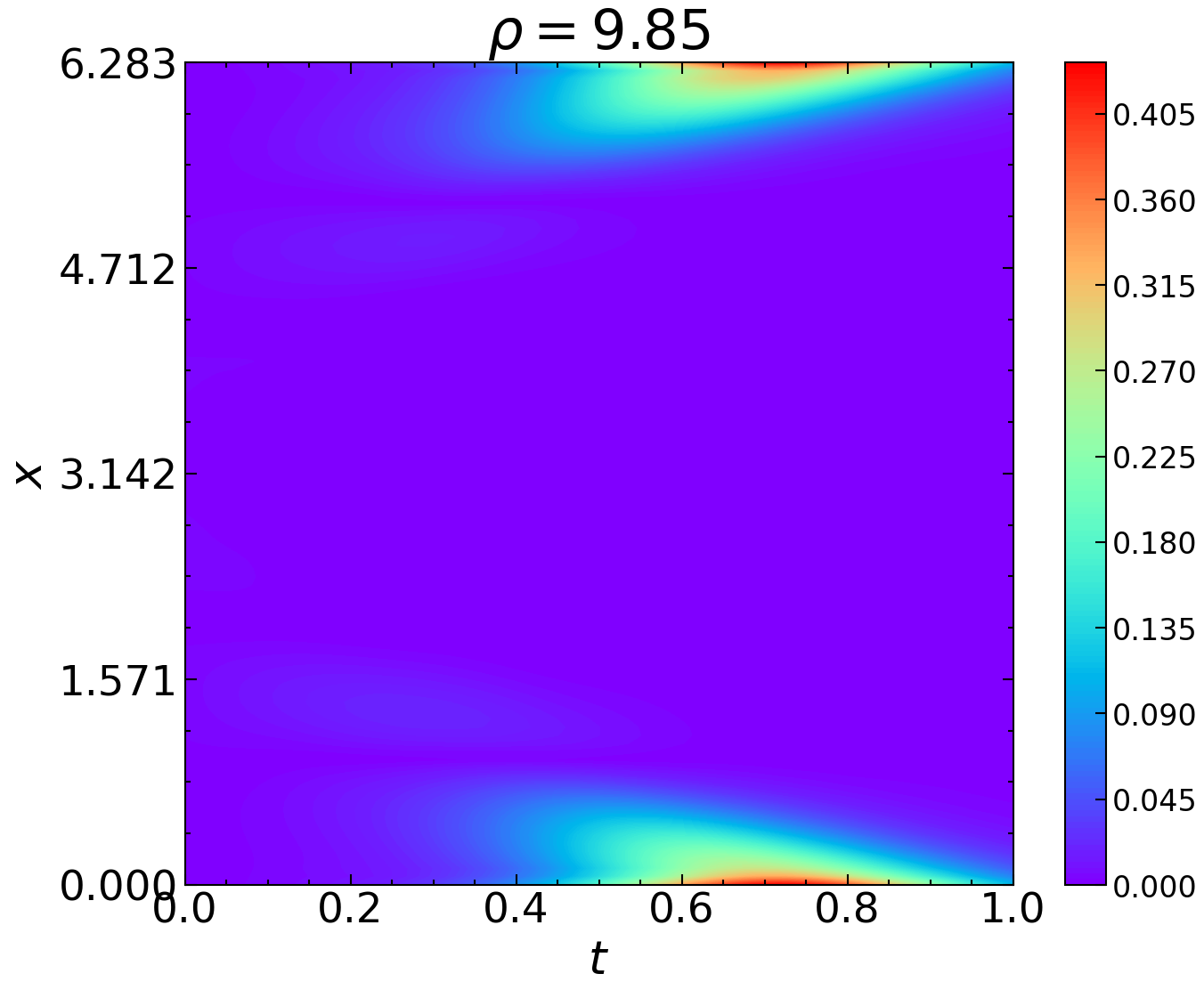}}
  \subfigure[TGPT-PINN Loss]{\includegraphics[width=0.24\linewidth]{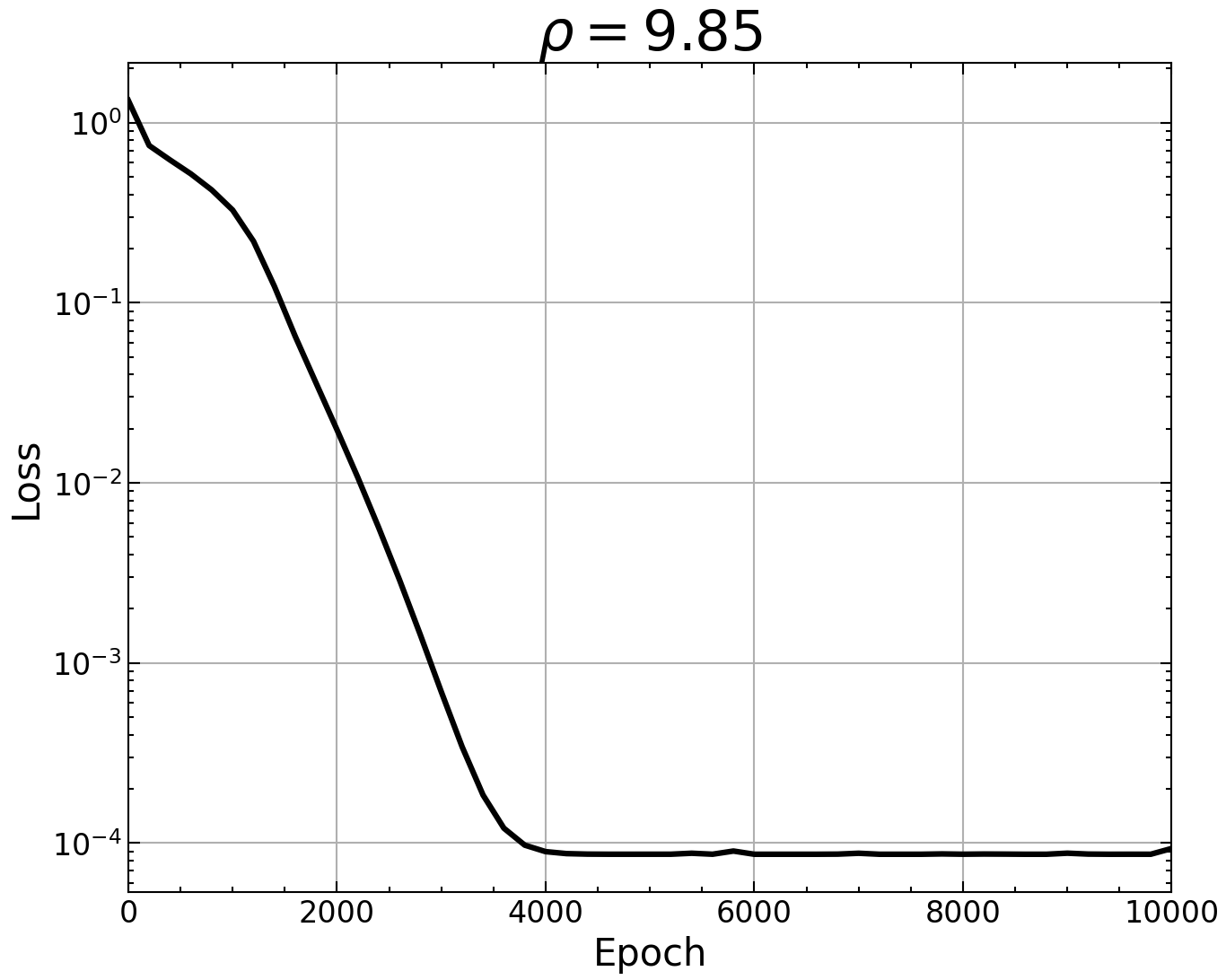}}
  \caption{Results for Section \ref{sssec:results-ppde-reaction}: GPT-PINN (top, with 10 neurons) and TGPT-PINN (bottom, with 1 neuron) results for $\rho = 9.85$.}
  \label{fig:R-gpt-pinn2}
\end{figure}

\noindent{\bf GPT/TGPT-PINN results:} We solve the problem using both the GPT-PINN and TGPT-PINN for $\rho \in [1,10]$, and compare the histories of convergence for the loss as the hyper-reduced neural networks gradually grow. The results of the loss and parameter selection are shown in \Cref{fig:R-GPT-PINN}. From the analytical solution of the problem, it can be observed that all solutions can be obtained through $\rho$-dependent $(x,t)$ alignment and shift operations from one reference solution. \Cref{fig:R-GPT-PINN} demonstrates this in practice: the TGPT-PINN with just one neuron can capture the parametric dependence and achieve better approximation than the (linear) GPT-PINN with $10$ neurons. This is true even though we have manually selected the single snapshot for the TGPT-PINN as the ``easiest'' (i.e., smallest) value of $\rho$ in the parameter domain. In \Cref{fig:R-gpt-pinn1} and \Cref{fig:R-gpt-pinn2} we show a comparison of solutions obtained by the GPT-PINN and TGPT-PINN for the relatively large values of $\rho = 3.25$ and $\rho = 9.85$.

\subsubsection{Reaction-diffusion PDE}\label{sssec:results-ppde-rd}
The reaction-diffusion system is where a diffusion operator is added to the reaction equation above. The system has the formulation with periodic boundary conditions as follows:
\begin{equation}
\begin{aligned}
  & \frac{\partial u}{\partial t}-\nu \frac{\partial^2 u}{\partial x^2}-\rho u(1-u)=0, \hskip 10pt (x,t) \in [0, 2\pi] \times [0, 1], \\
  & u(x,0) = h(x), \hskip 10pt x \in [0, 2\pi], \\
  & u(0, t) = u(2\pi t), \hskip 10pt t \geq 0,
\end{aligned}
\end{equation}
where $h(x)$ is as in \eqref{eq:h-def}, with $\nu > 0$ the diffusion coefficient, and $F$ is the Fourier transform operator.

This equation has the following analytical solution:
\begin{equation}\label{eq:rd-exact}
u(x, t)=F^{-1}\left(F\left(\frac{h(x) \exp (\rho t)}{h(x) \exp (\rho t)+1-h(x)}\right) e^{-\nu k^2 t}\right).
\end{equation}

\begin{figure}[htbp]
  \centering
  \subfigure[Exact solution]
  {\includegraphics[width=0.24\linewidth]{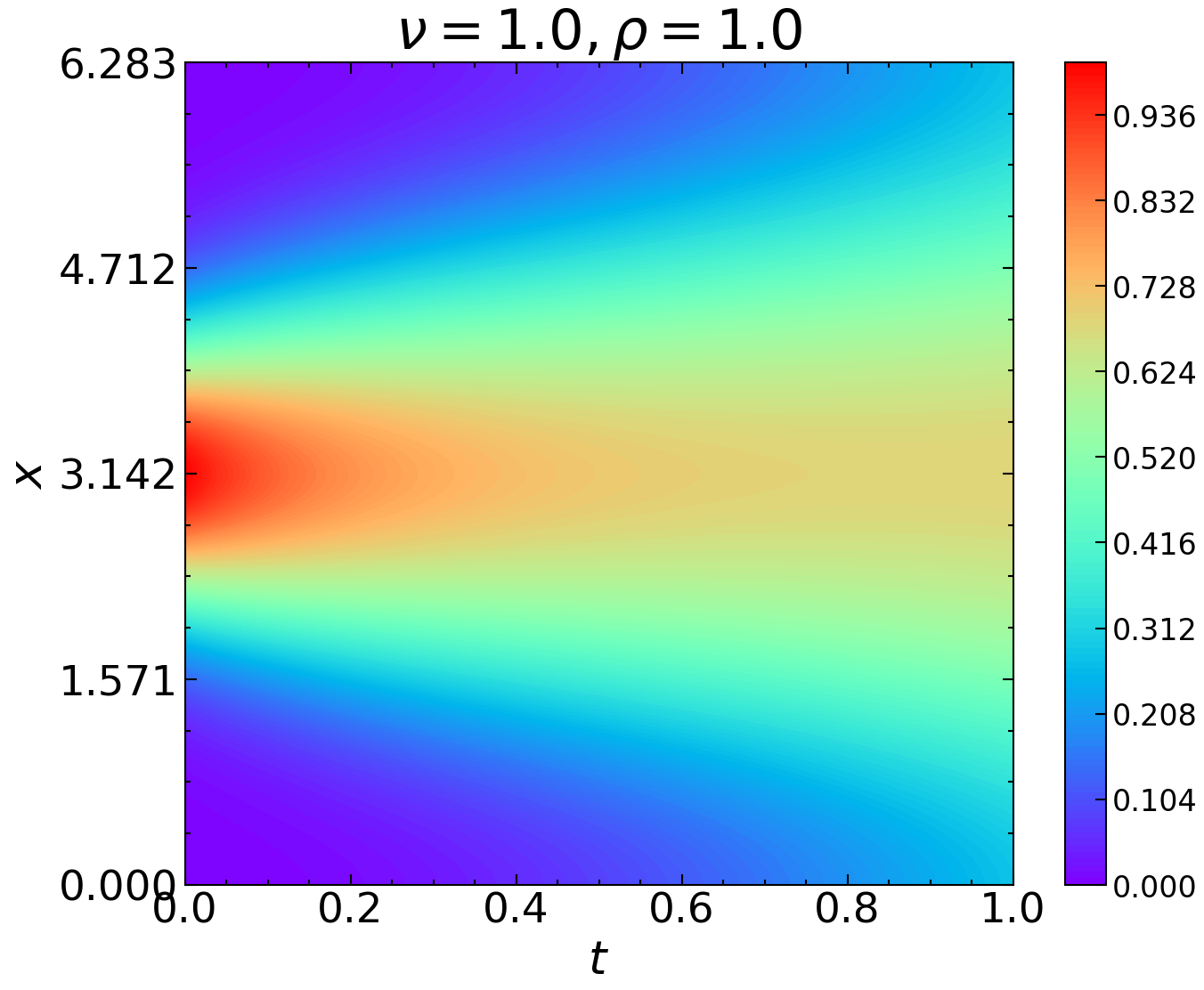}}
  \subfigure[PINN solution]{\includegraphics[width=0.24\linewidth]{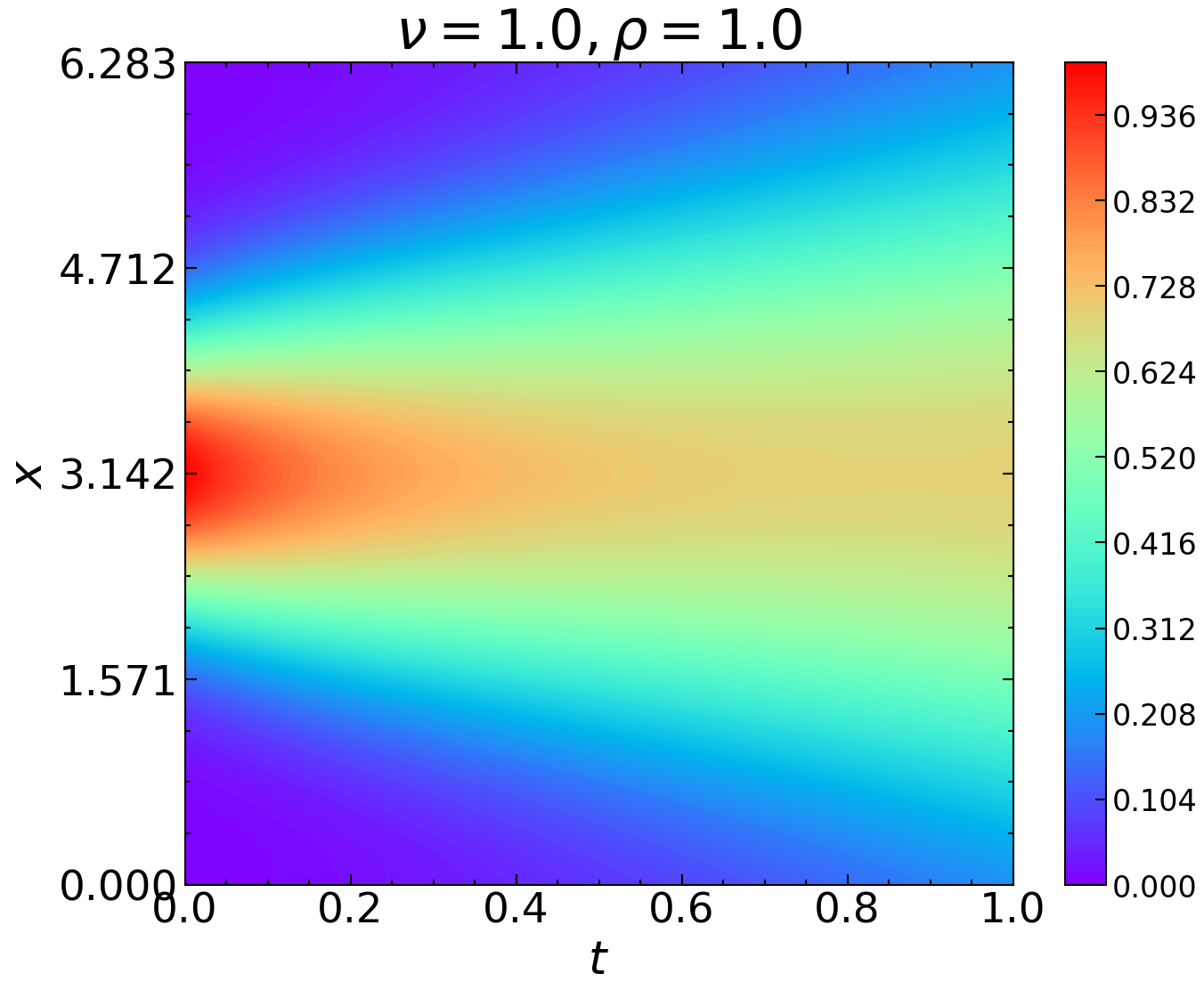}}
  \subfigure[PINN error]{\includegraphics[width=0.24\linewidth]{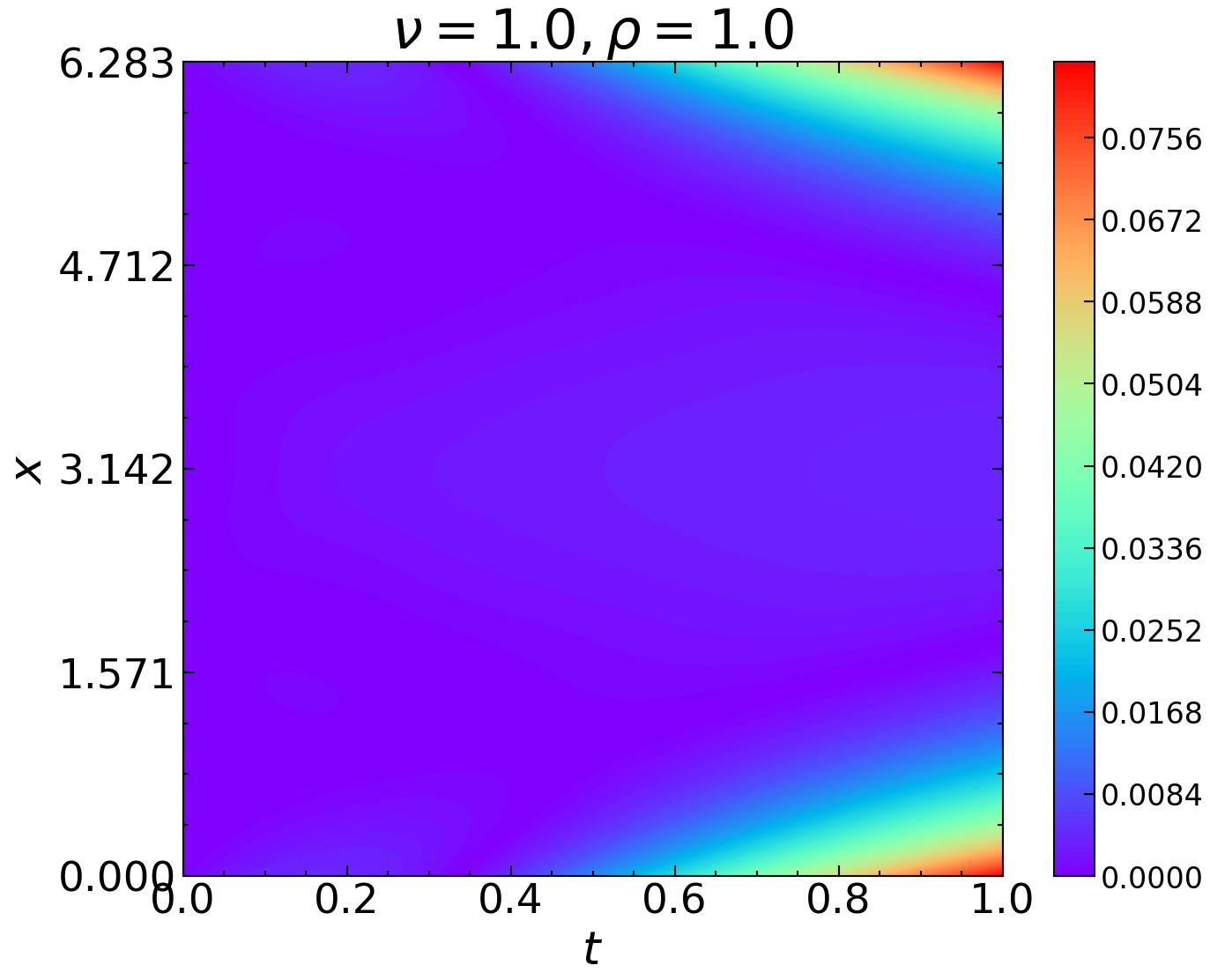}}
  \subfigure[PINN Loss]{\includegraphics[width=0.24\linewidth]{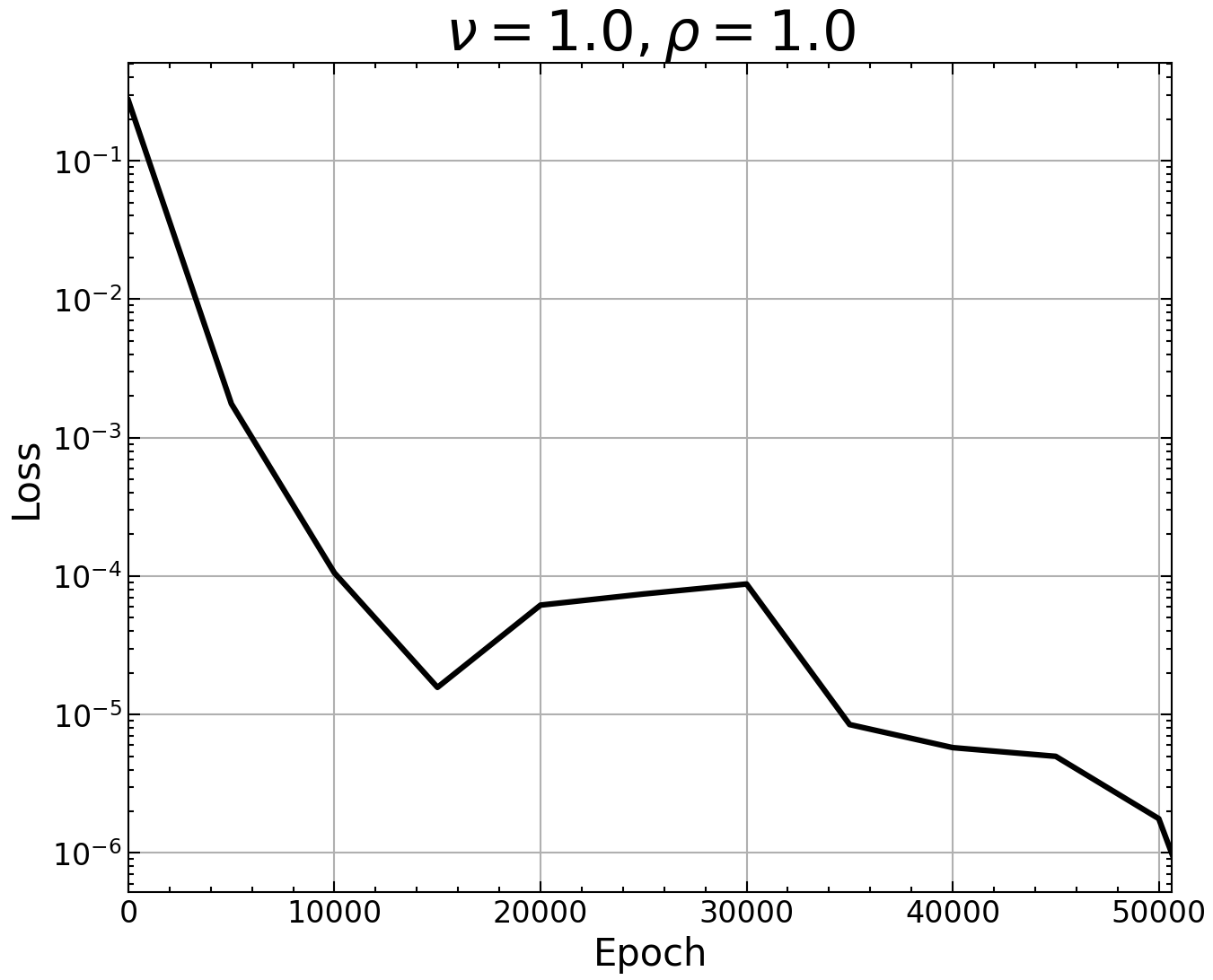}}
  \subfigure[Exact solution]{\includegraphics[width=0.24\linewidth]{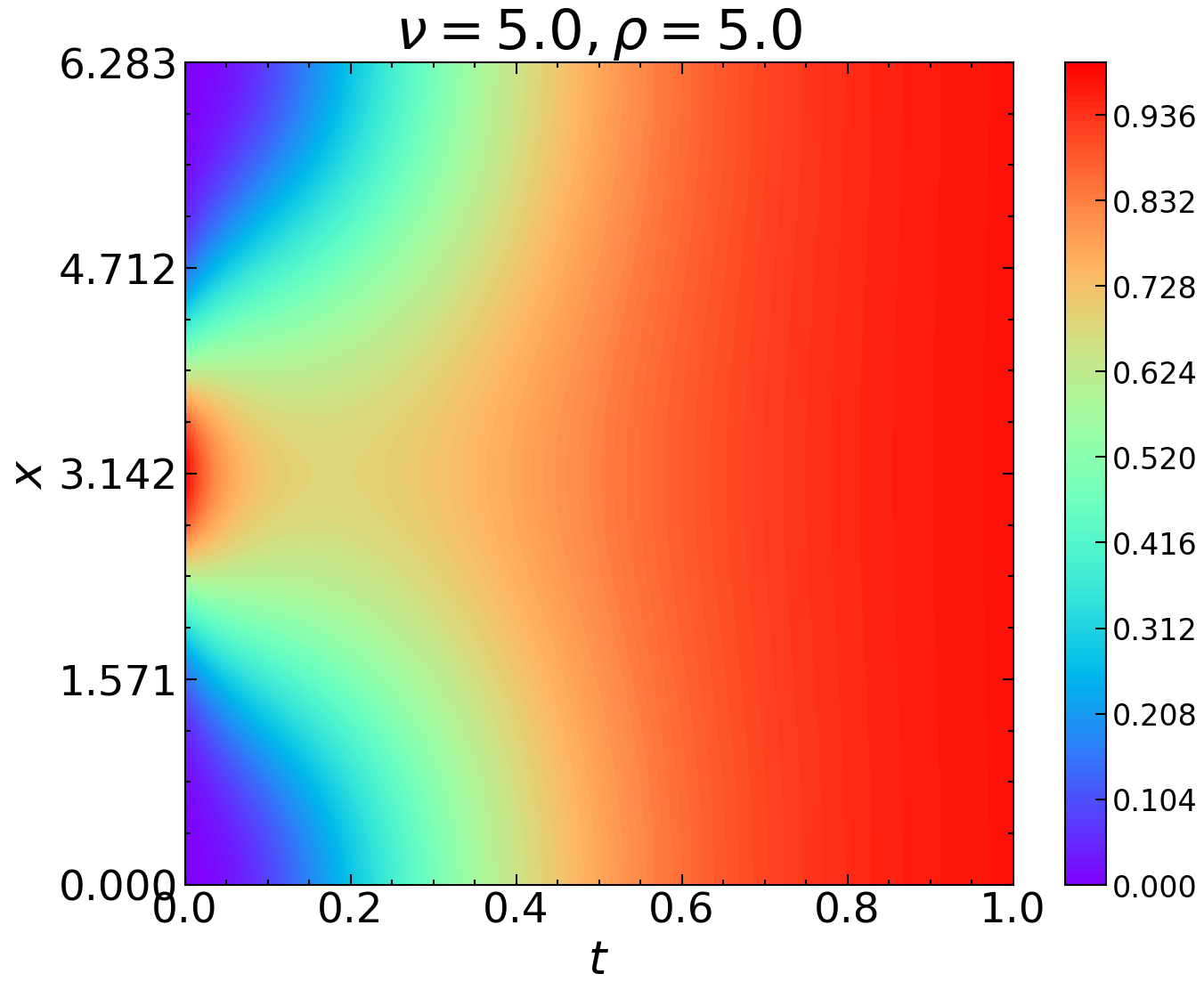}}
  \subfigure[PINN solution]{\includegraphics[width=0.24\linewidth]{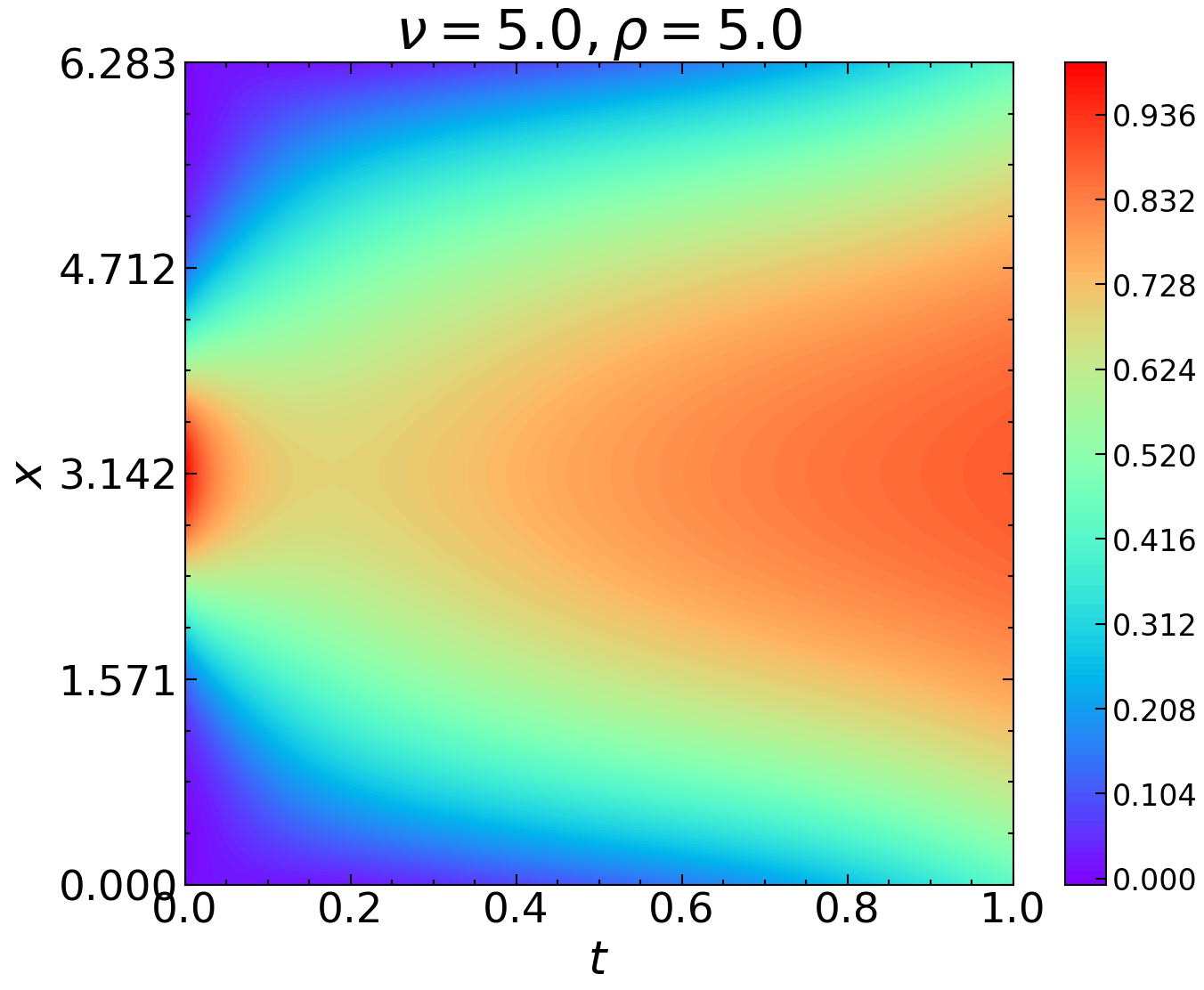}}
  \subfigure[PINN error]{\includegraphics[width=0.24\linewidth]{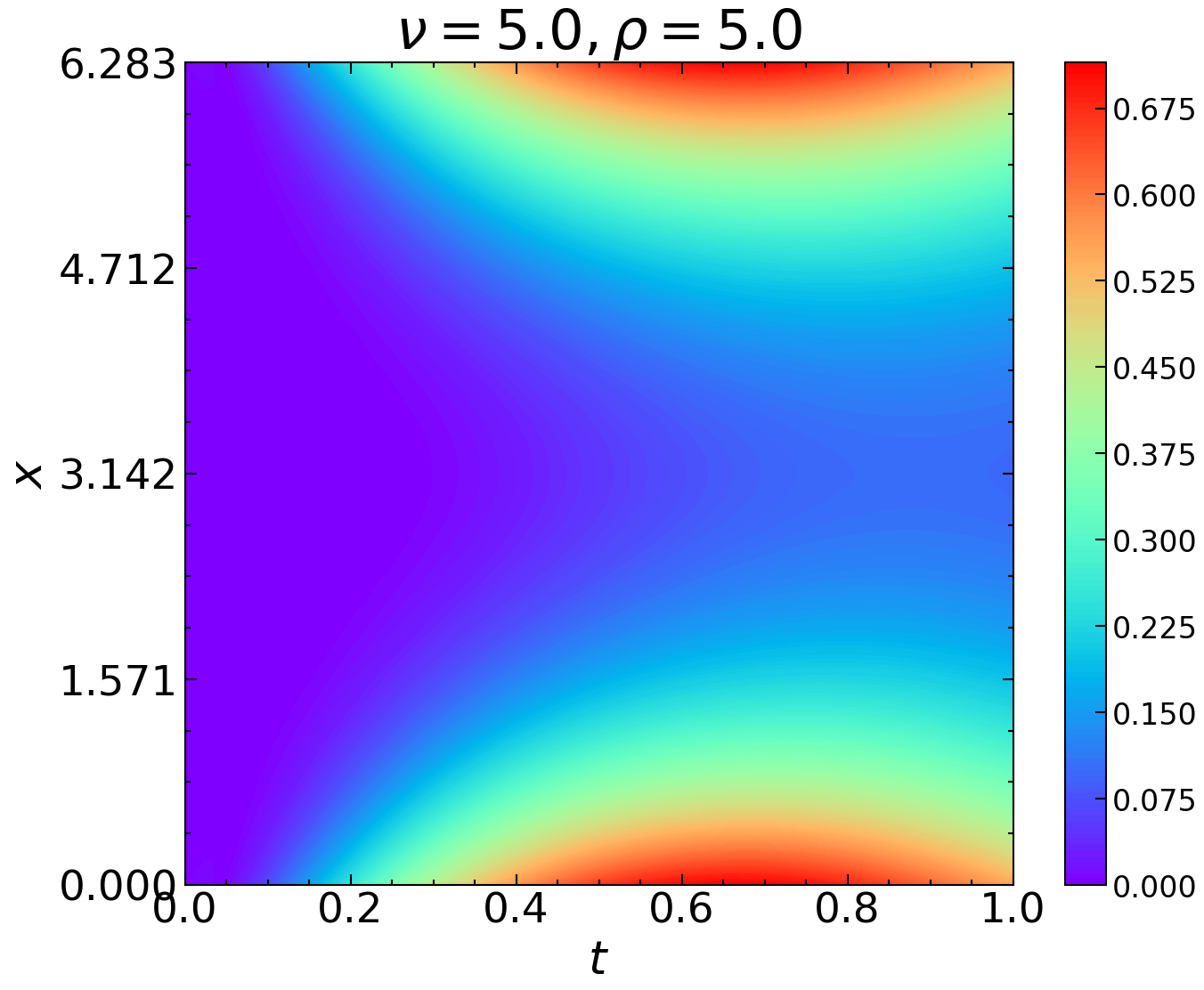}}
  \subfigure[PINN Loss]{\includegraphics[width=0.24\linewidth]{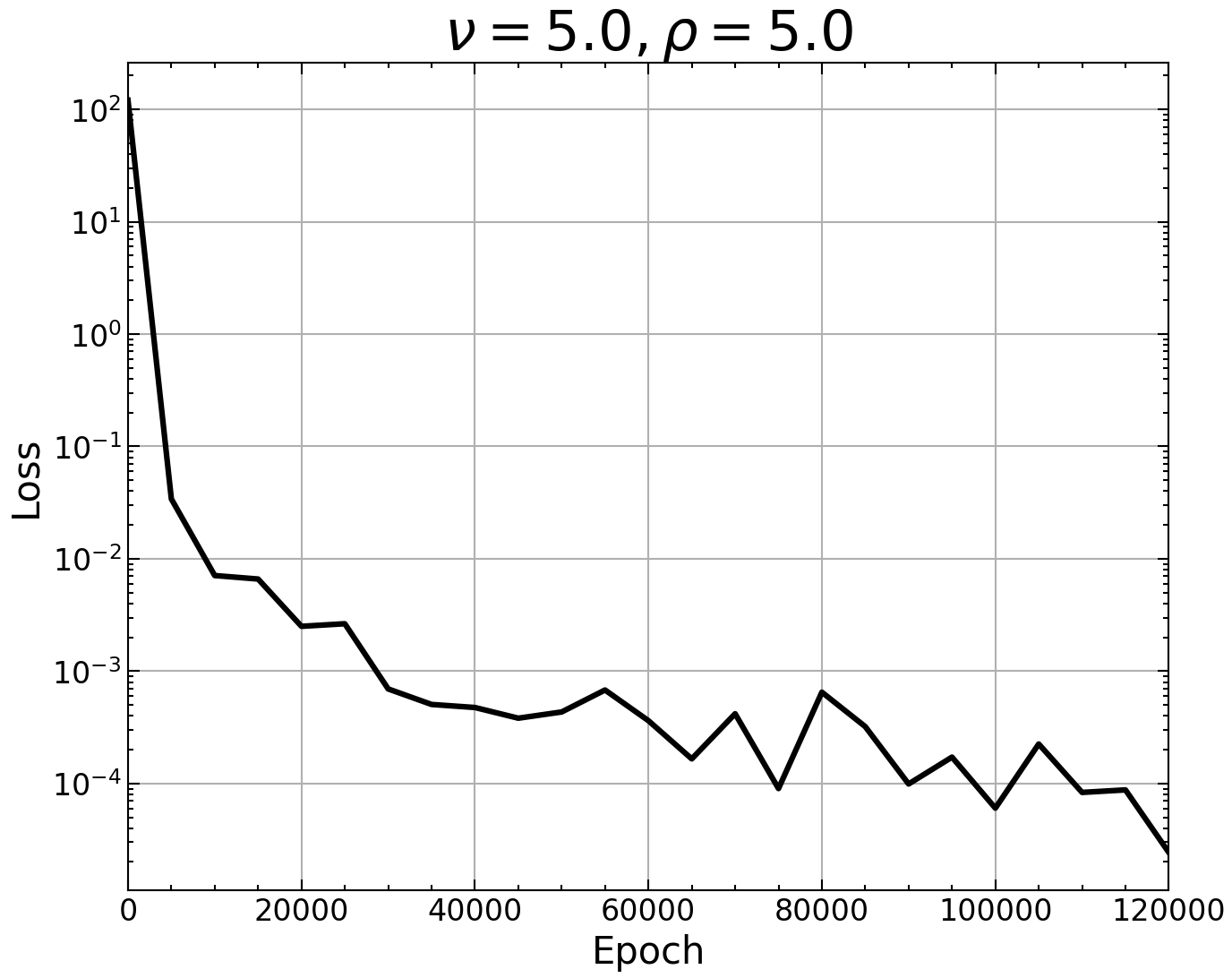}}
  \caption{Results for Section \ref{sssec:results-ppde-rd}: Full PINN for $\nu=1,\rho=1$ (top) and $\nu=5,\rho=5$ (bottom).}
  \label{fig:RD-Full-PINN1}
\end{figure}

\noindent {\bf PINN results:} The PINN snapshot solver employs a fully connected neural network with dimensions [2, 40, 40, 40, 1] with $\lambda \equiv 1$ in $\mathcal{L}_{\rm int}(u)$ of the PINN loss function \eqref{eq:continousloss}. The training samples are randomly selected, but we again employ the optimized activation function from \cite{zhao2023pinnsformer} in \eqref{eq:waveact}. In \Cref{fig:RD-Full-PINN1} we show the accuracy and loss reduction of PINN for the parameter choice $(\nu, \rho) = (1,1)$ (top row) and $(5,5)$ (bottom row). Note that as the parameters $\nu,\rho$ increase, the difficulty of solving the problem significantly increases. 

\begin{figure}[htbp]
  \centering
  \includegraphics[width=0.32\linewidth]{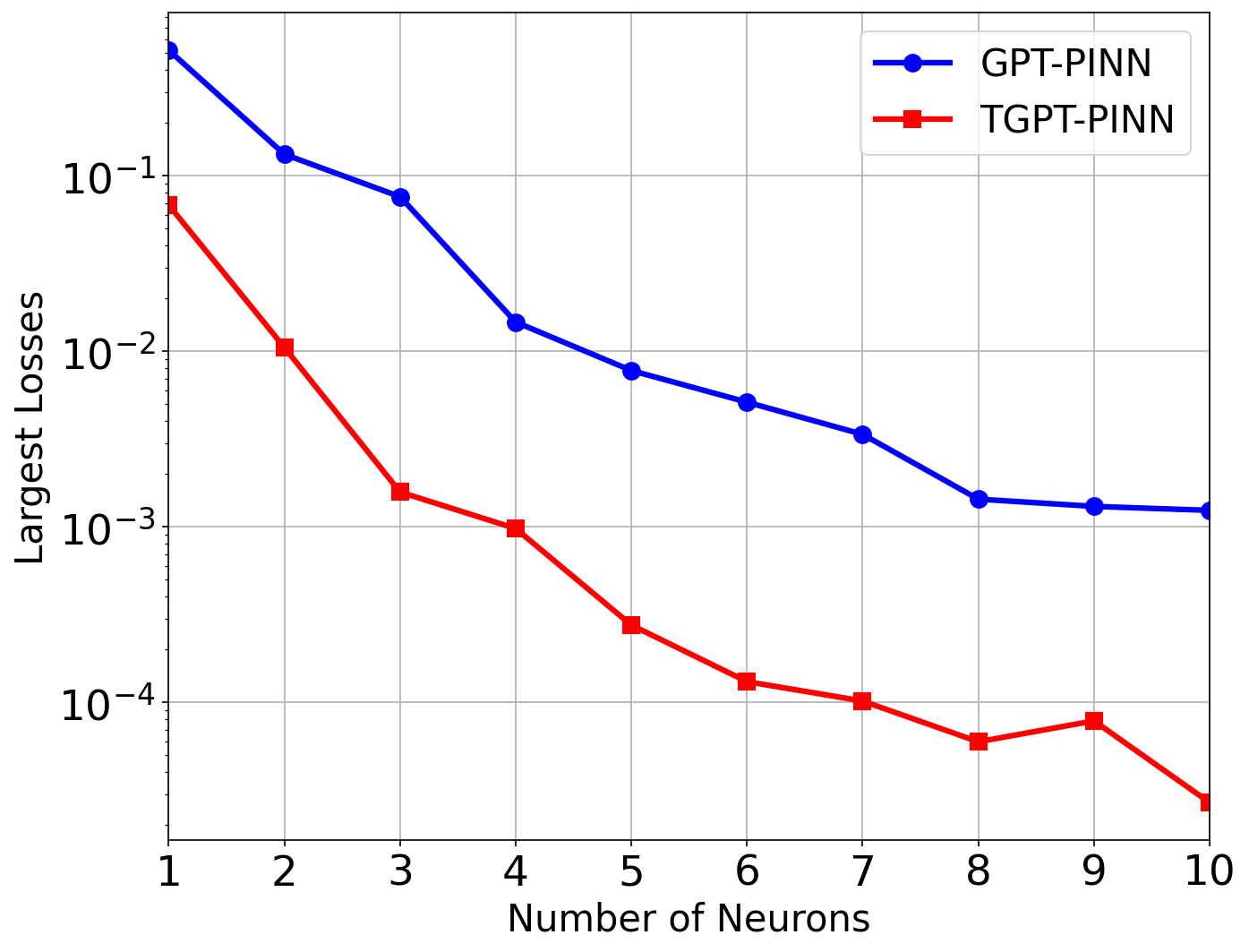}
  \includegraphics[width=0.31\linewidth]{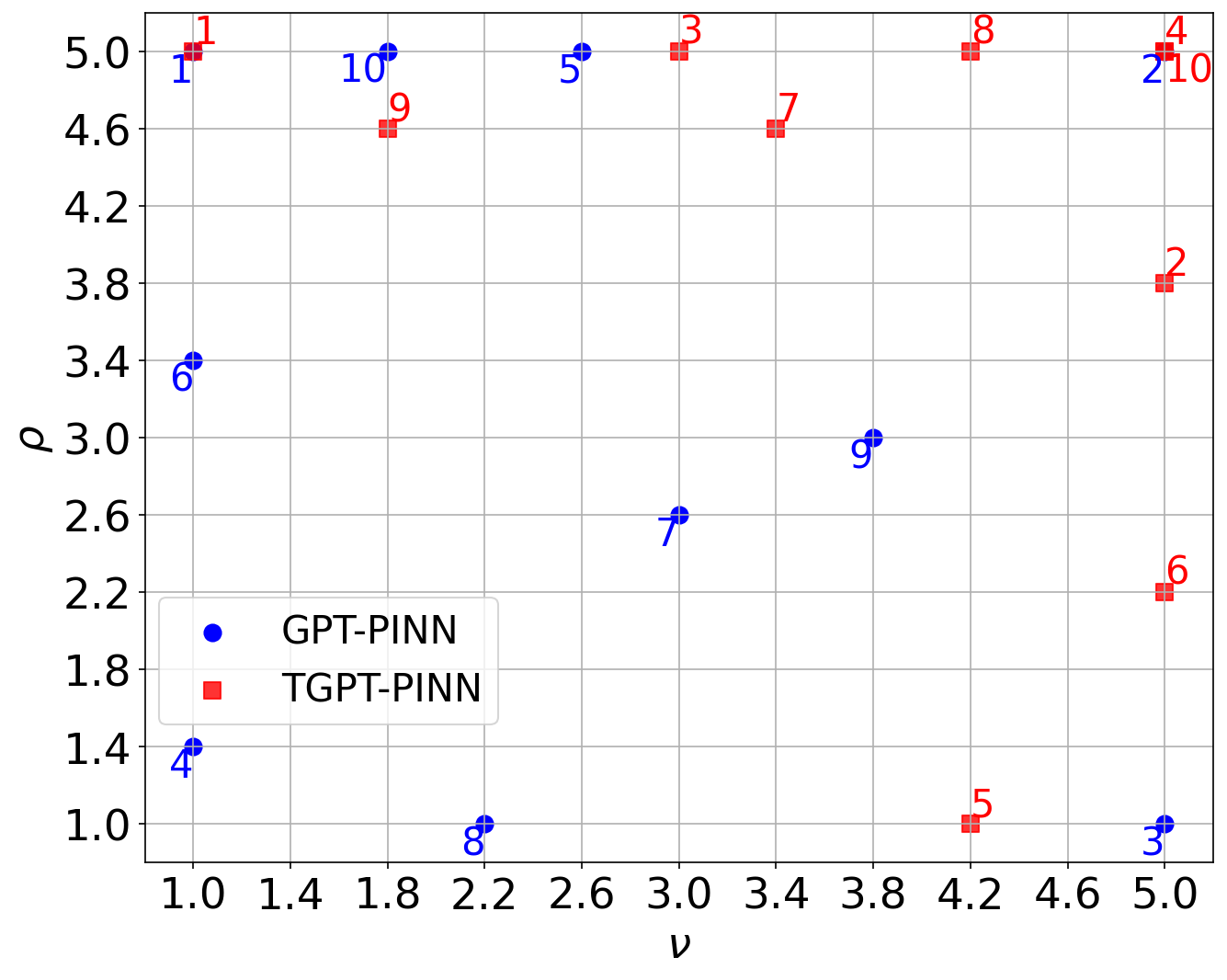}
  \includegraphics[width=0.32\linewidth]{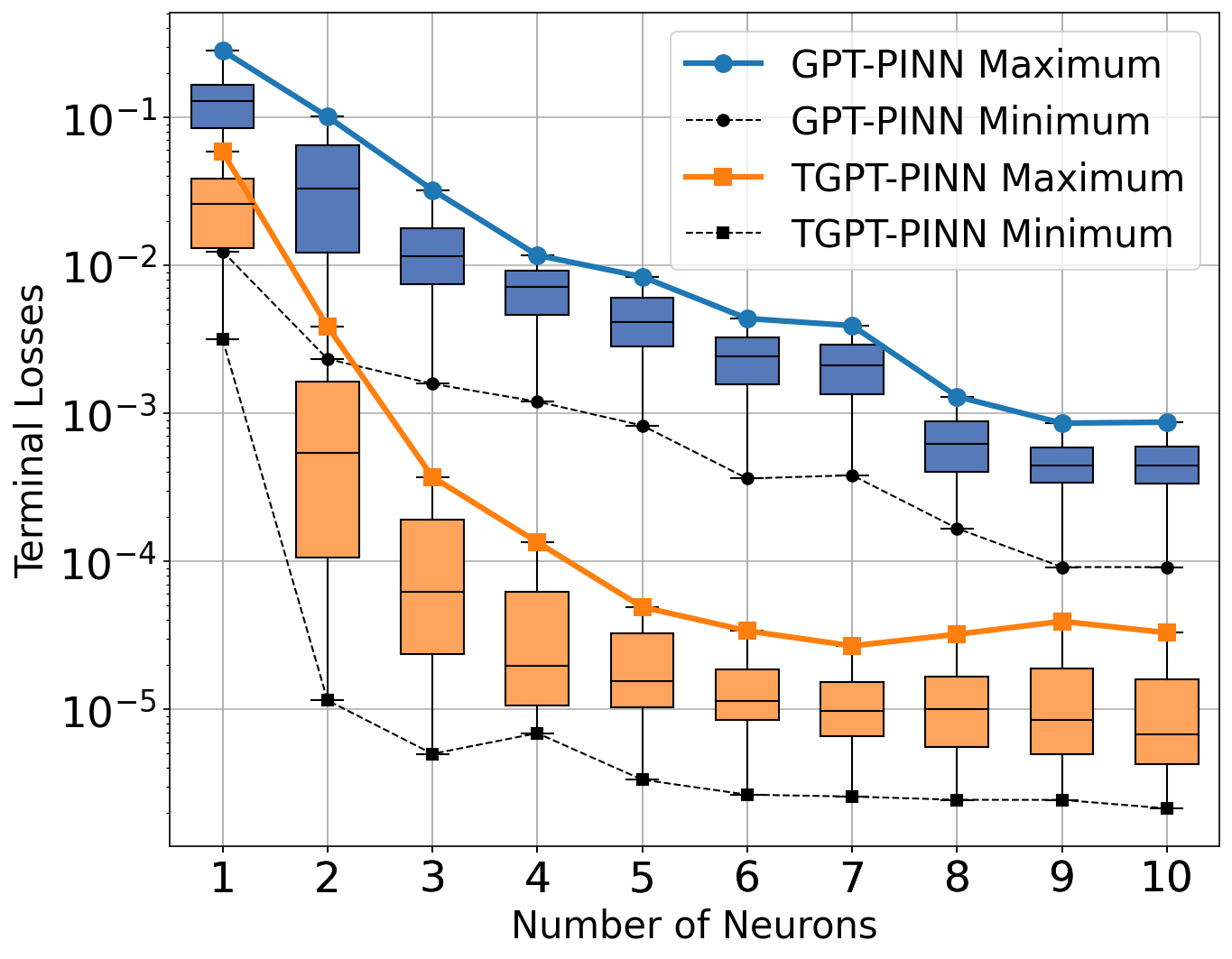}
\put(-360,18){\makebox(0,0){{(a)}}}
\put(-230,18){\makebox(0,0){{(b)}}}
\put(-100,18){\makebox(0,0){{(c)}}}
  \caption{Results for Section \ref{sssec:results-ppde-rd} for (T)GPT-PINN for the reaction-diffusion equation: (a) Worst-case history of convergence during training, (b) parameter domain locations determined by the methods, and (c) box plots of the losses when tested online on unseen locations.}
  \label{fig:RD-GPT-PINN}
\end{figure}

\begin{figure}[htbp]
  \centering
  \subfigure[Exact solution]{\includegraphics[width=0.24\linewidth]{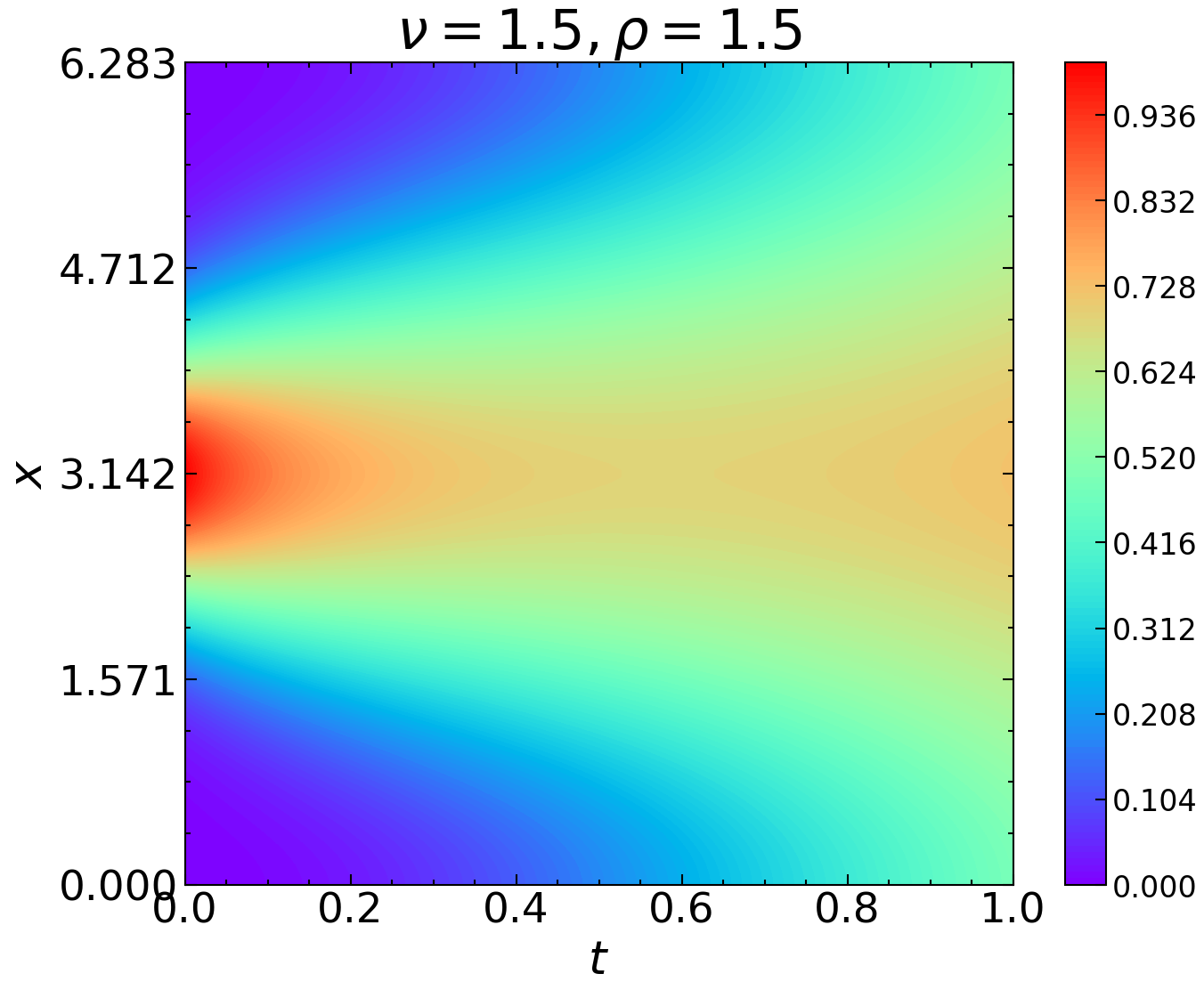}}
  \subfigure[GPT-PINN solution]{\includegraphics[width=0.24\linewidth]{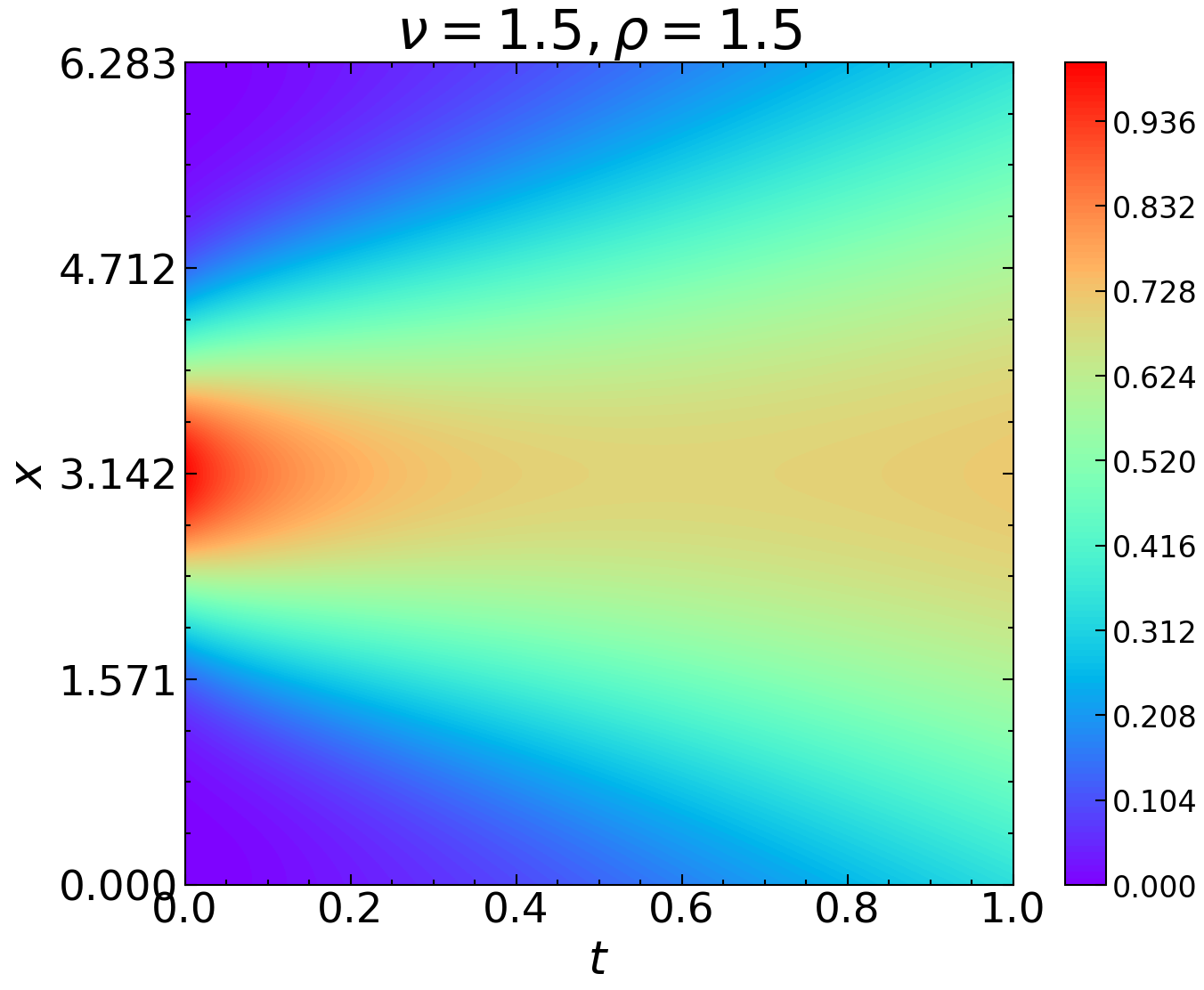}}
  \subfigure[GPT-PINN error]{\includegraphics[width=0.24\linewidth]{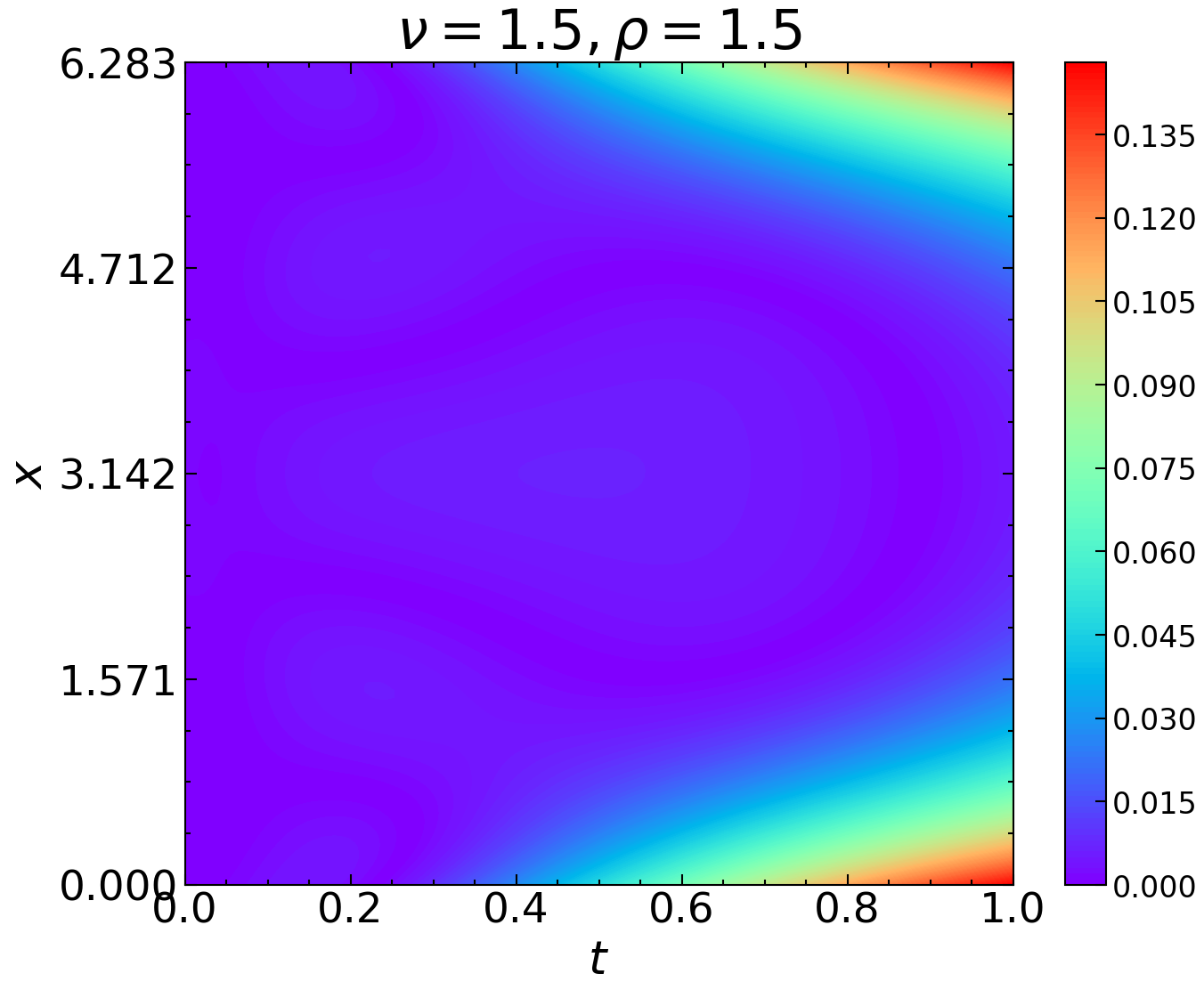}}
  \subfigure[GPT-PINN Loss]{\includegraphics[width=0.237\linewidth]{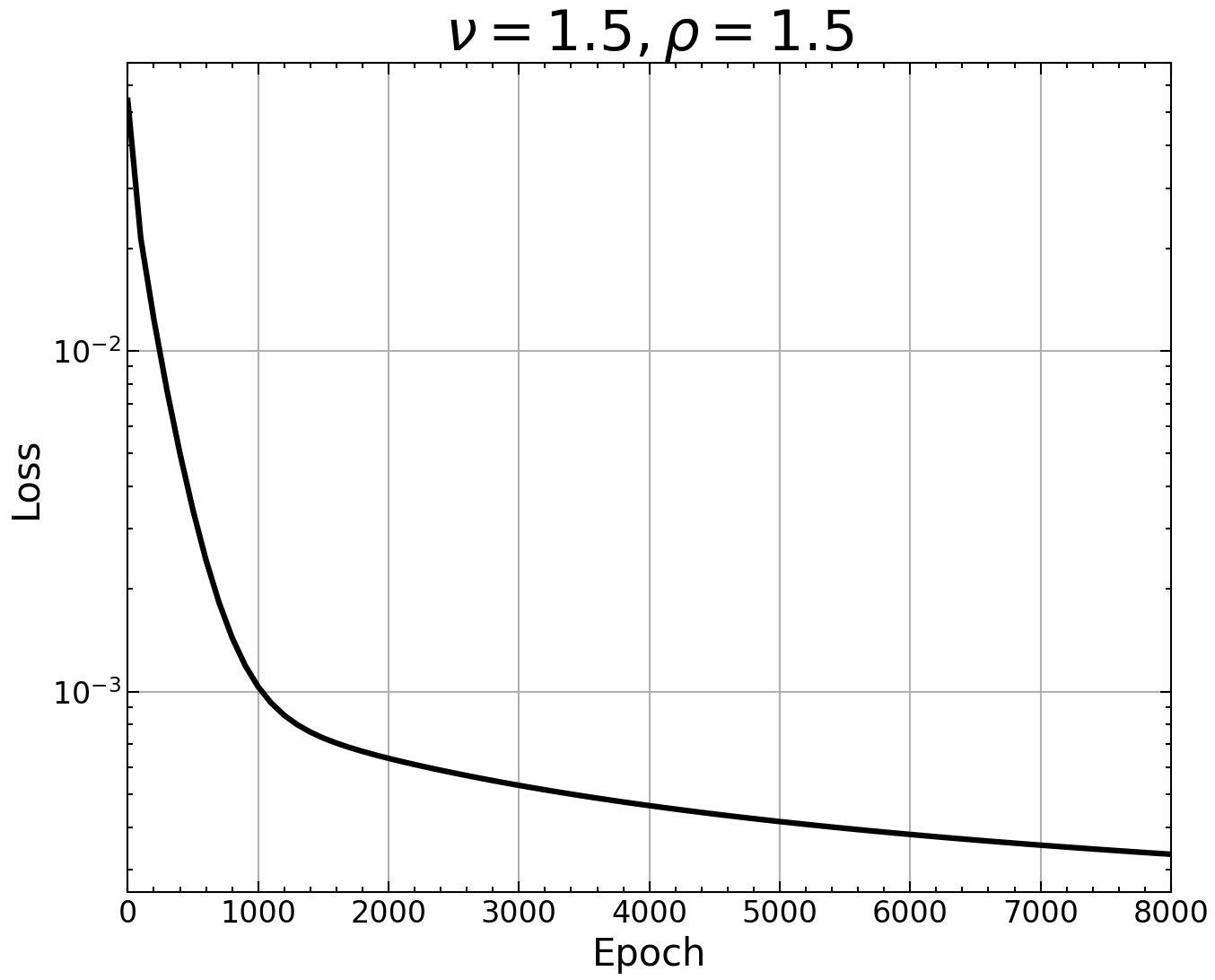}}
  
  \subfigure[Exact solution]{\includegraphics[width=0.24\linewidth]{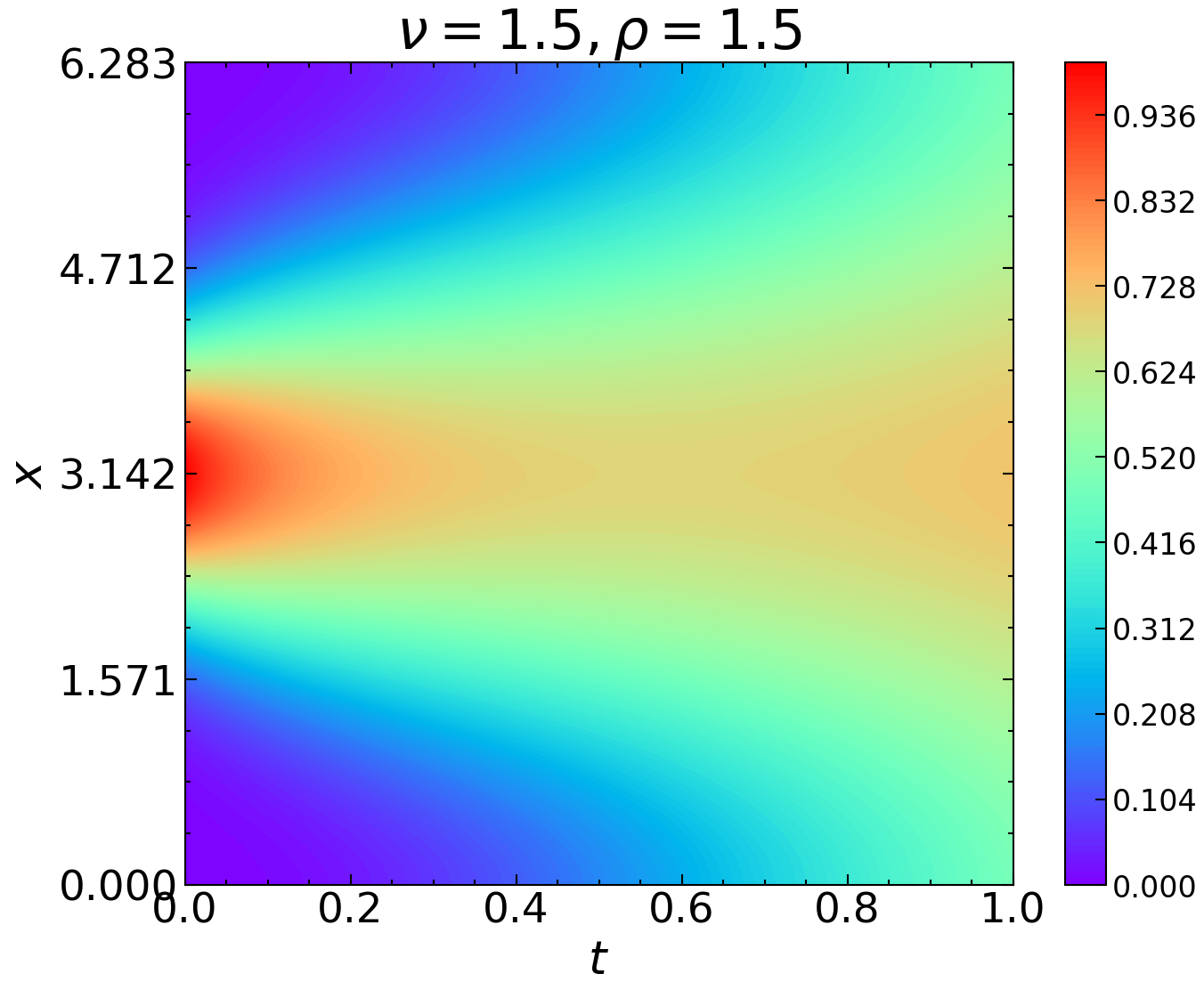}}
  \subfigure[TGPT-PINN solution]{\includegraphics[width=0.24\linewidth]{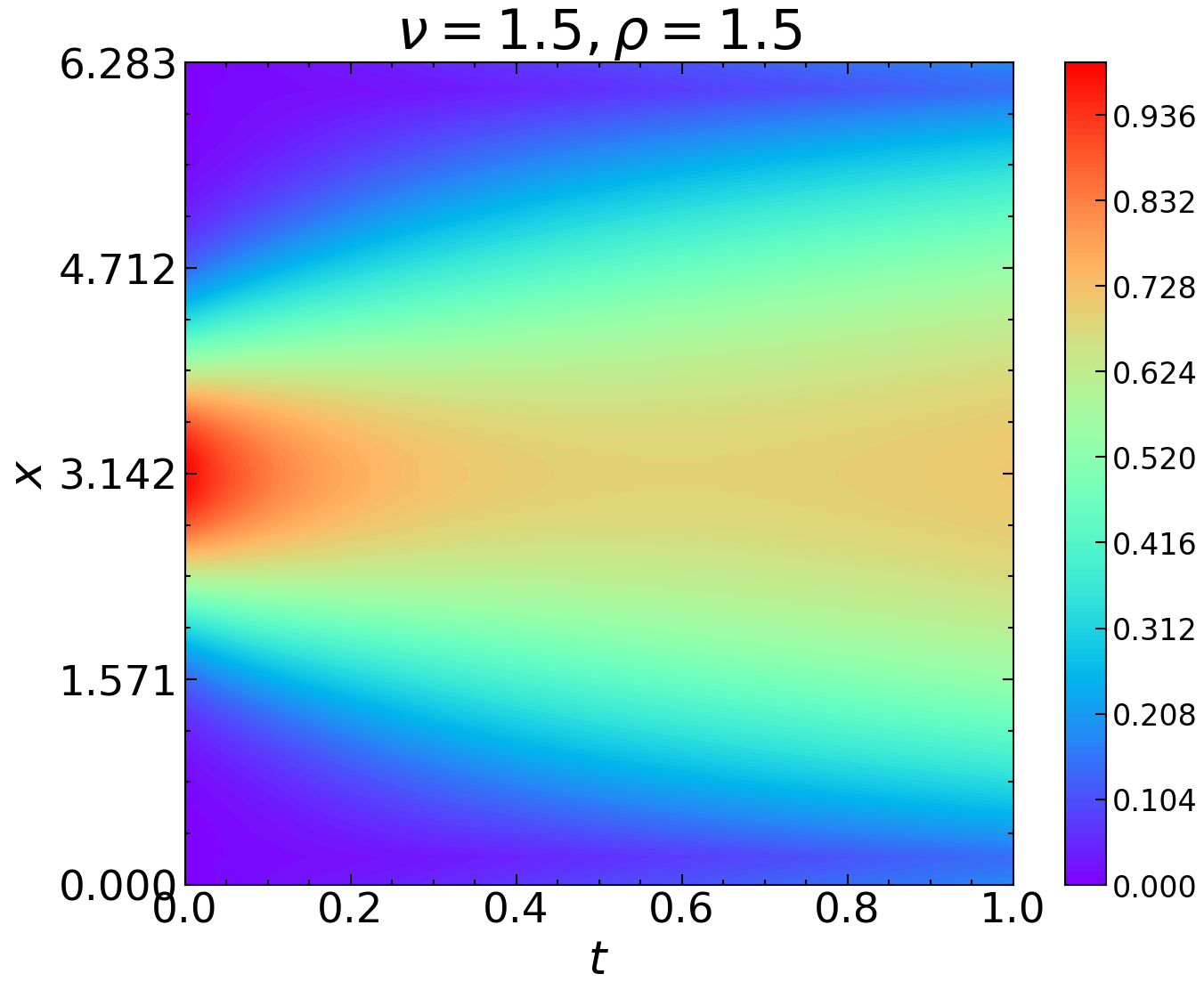}}
  \subfigure[TGPT-PINN error]{\includegraphics[width=0.24\linewidth]{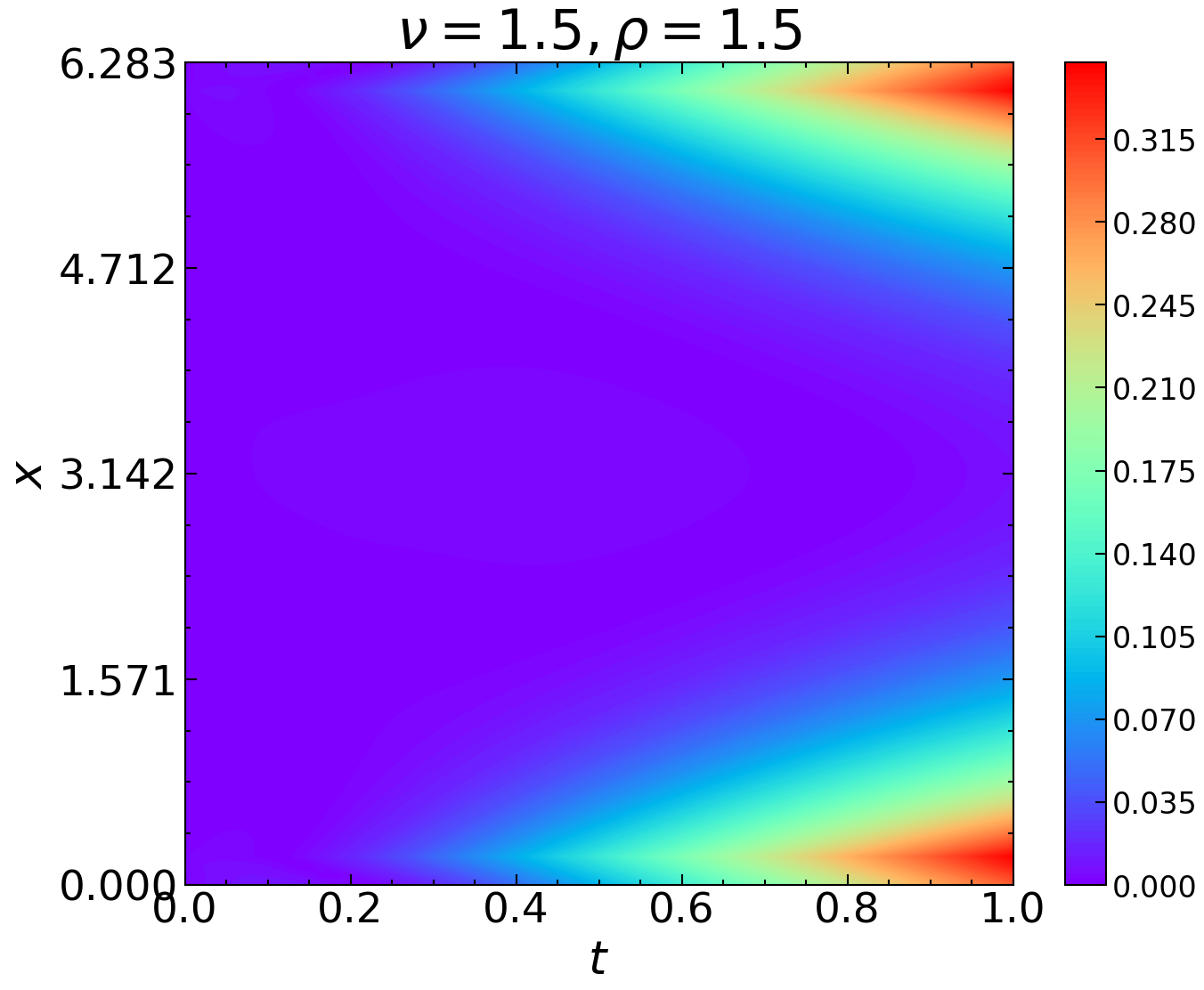}}
  \subfigure[TGPT-PINN Loss]{\includegraphics[width=0.237\linewidth]{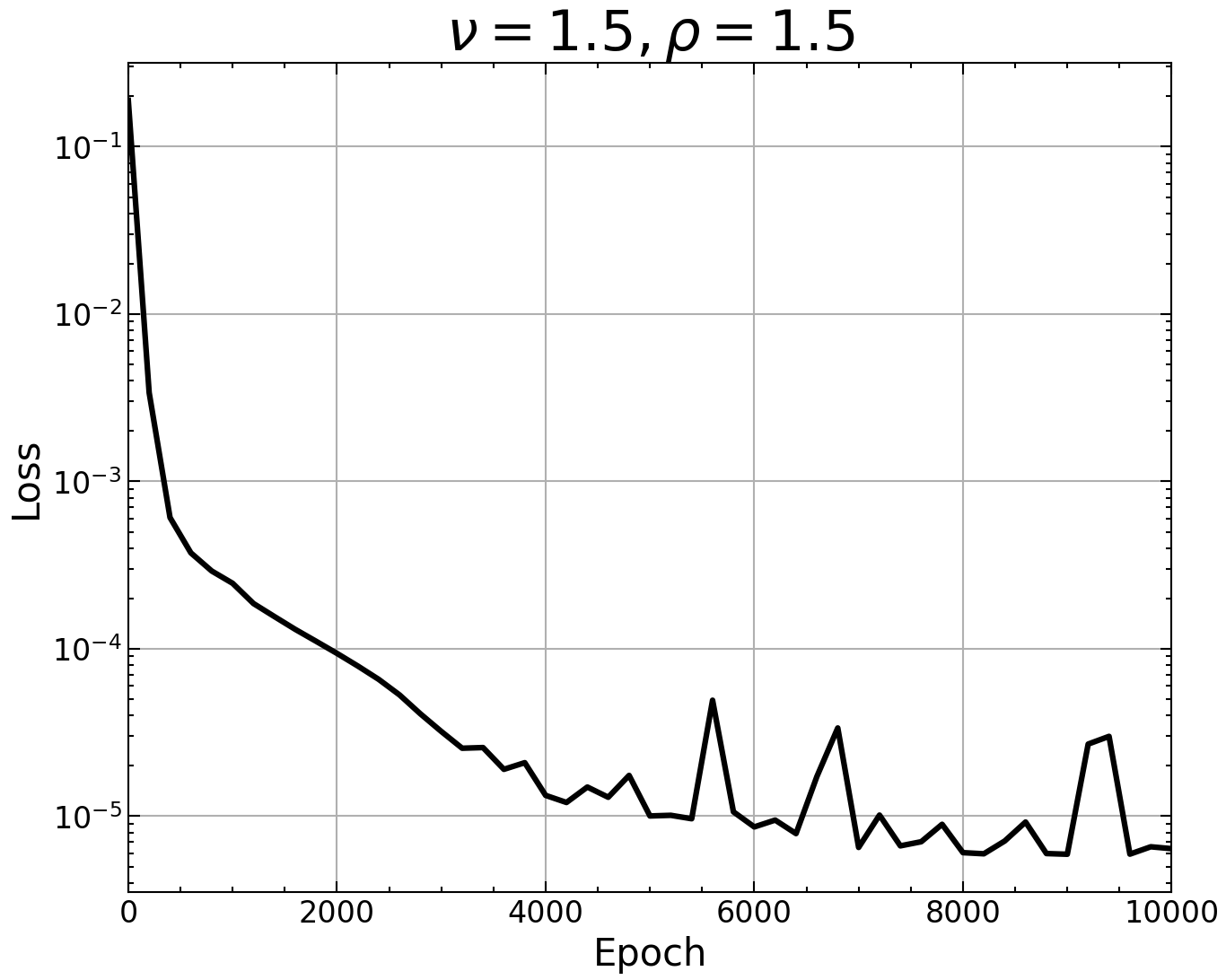}}
  \caption{Results for Section \ref{sssec:results-ppde-rd}: GPT-PINN (top, with 10 neurons) and TGPT-PINN (bottom, with 3 neurons) results for $\rho = 1.5$ and $\nu=1.5$.}
  \label{fig:RD-gpt-pinn1}
\end{figure}

\begin{figure}[htbp]
  \centering
  \subfigure[Exact solution]{\includegraphics[width=0.24\linewidth]{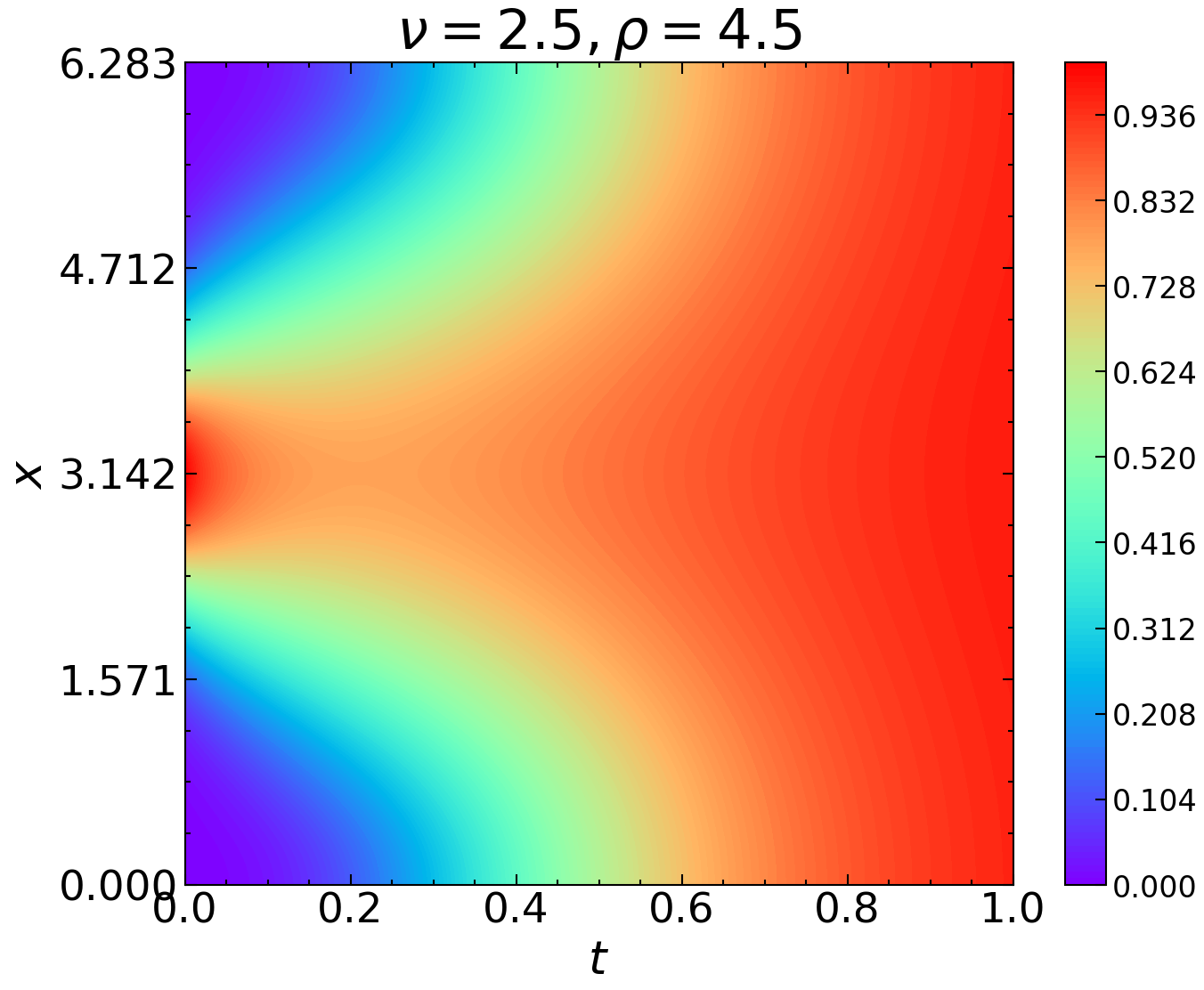}}
  \subfigure[GPT-PINN solution]{\includegraphics[width=0.24\linewidth]{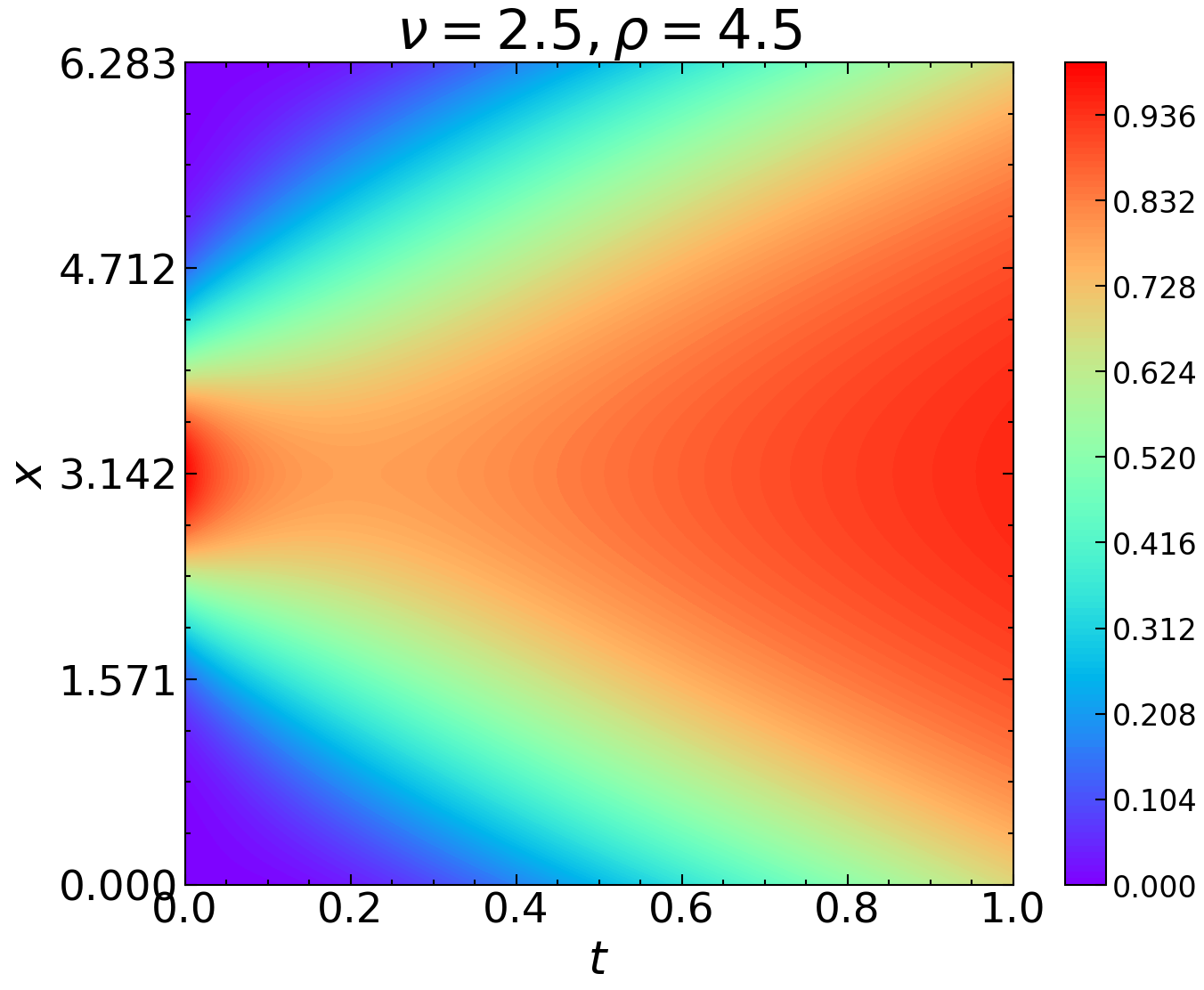}}
  \subfigure[GPT-PINN error]{\includegraphics[width=0.24\linewidth]{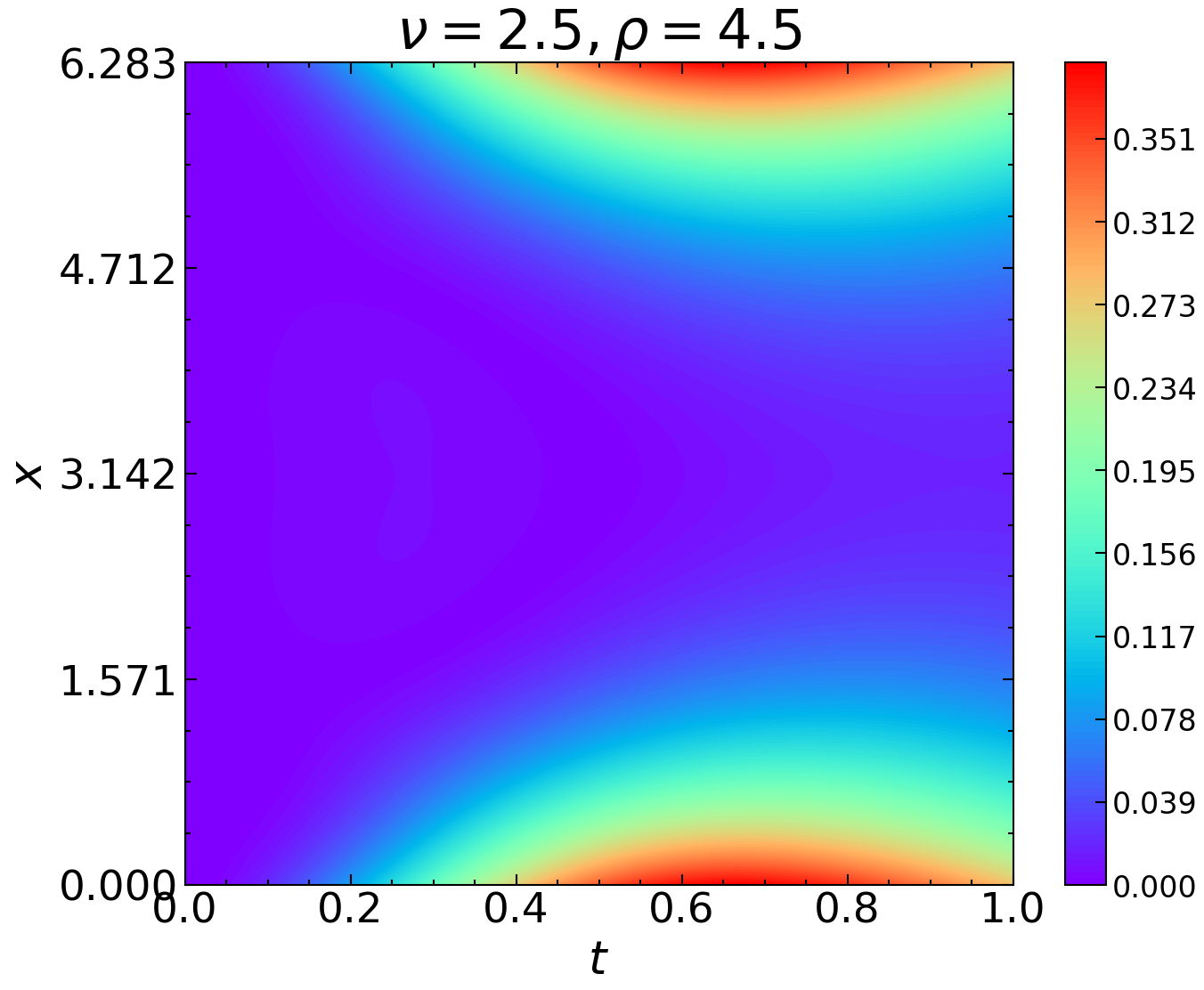}}
  \subfigure[GPT-PINN Loss]{\includegraphics[width=0.237\linewidth]{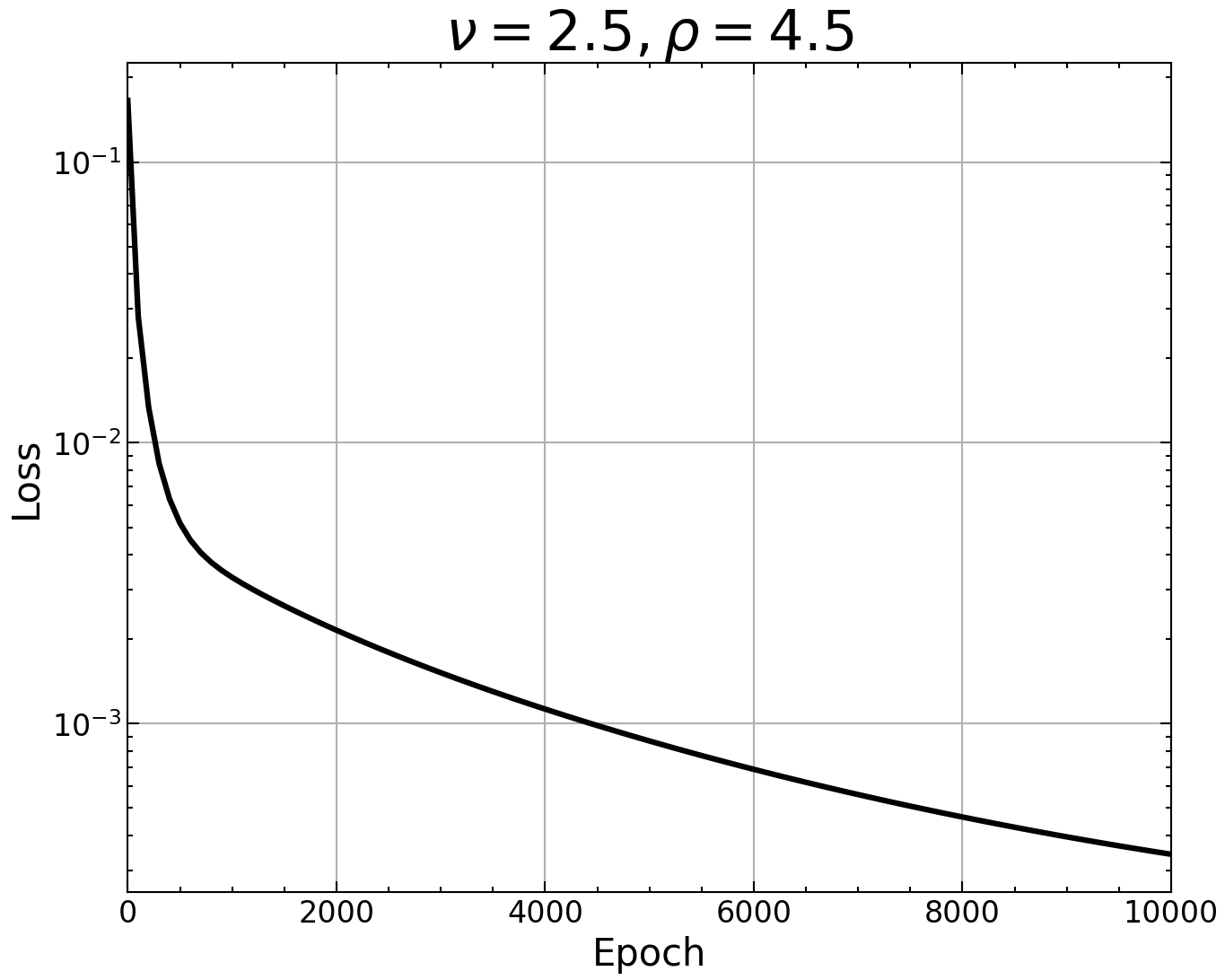}}

  \subfigure[Exact solution]{\includegraphics[width=0.24\linewidth]{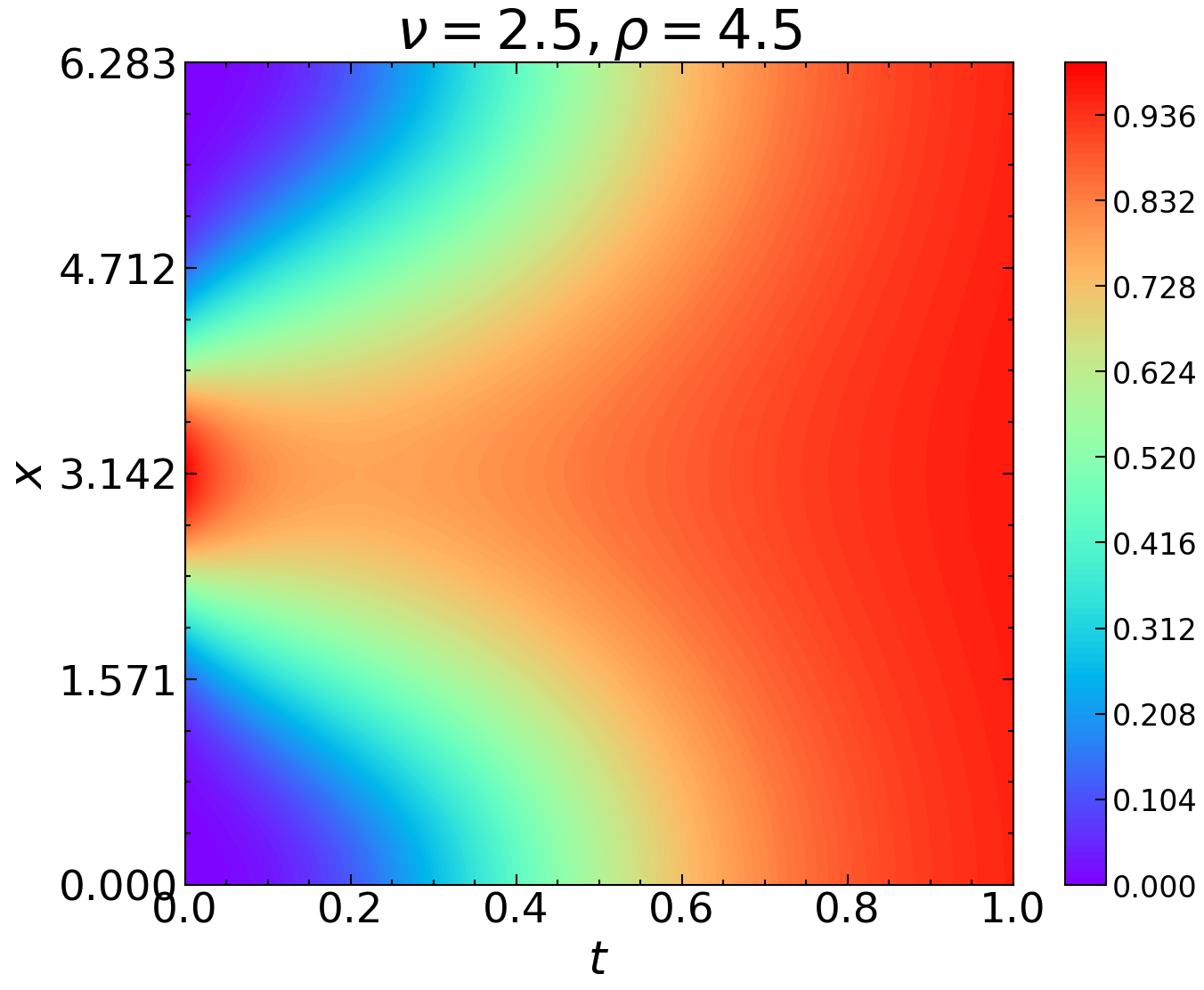}}
  \subfigure[TGPT-PINN solution]{\includegraphics[width=0.24\linewidth]{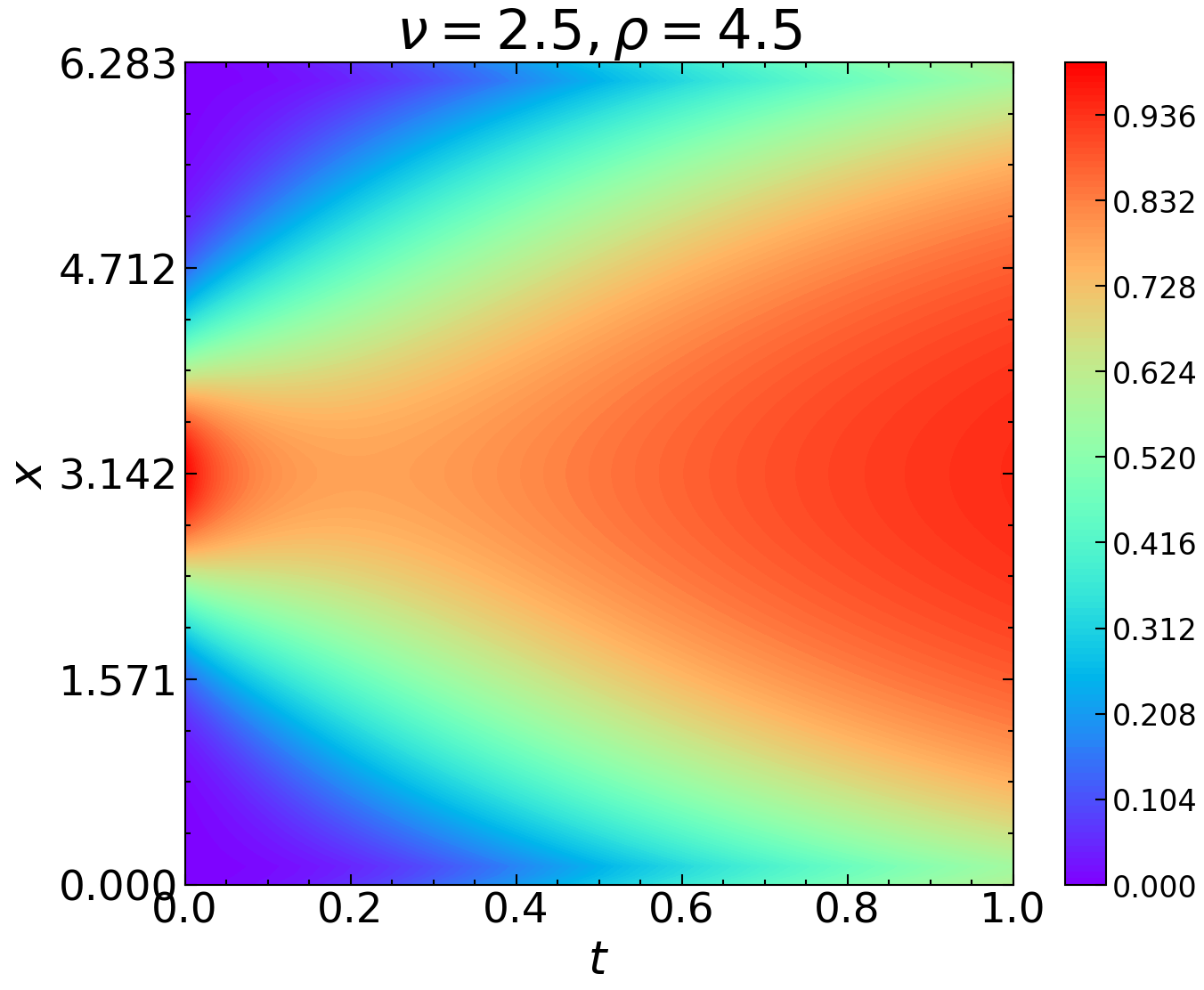}}
  \subfigure[TGPT-PINN error]{\includegraphics[width=0.24\linewidth]{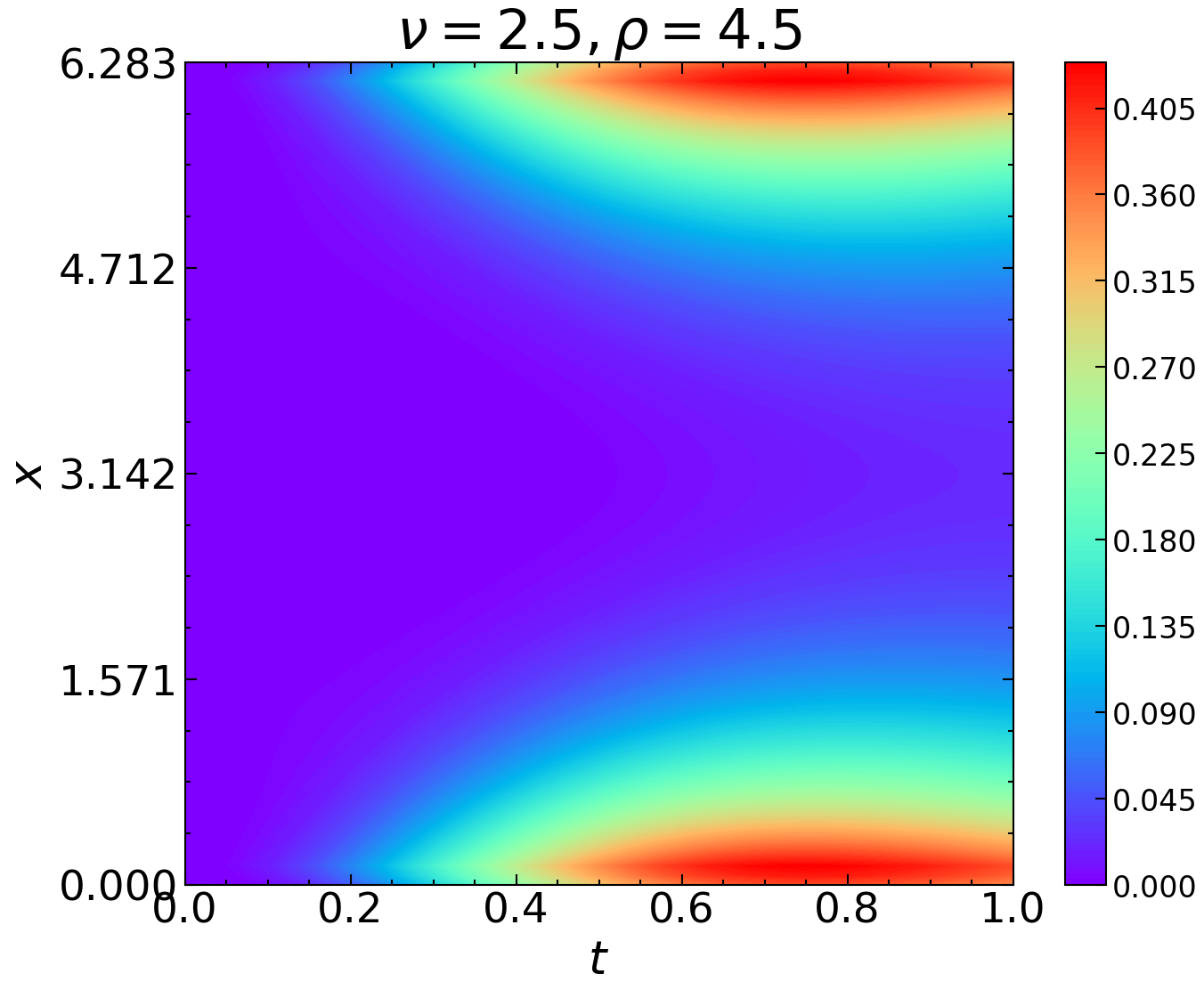}}
  \subfigure[TGPT-PINN Loss]{\includegraphics[width=0.237\linewidth]{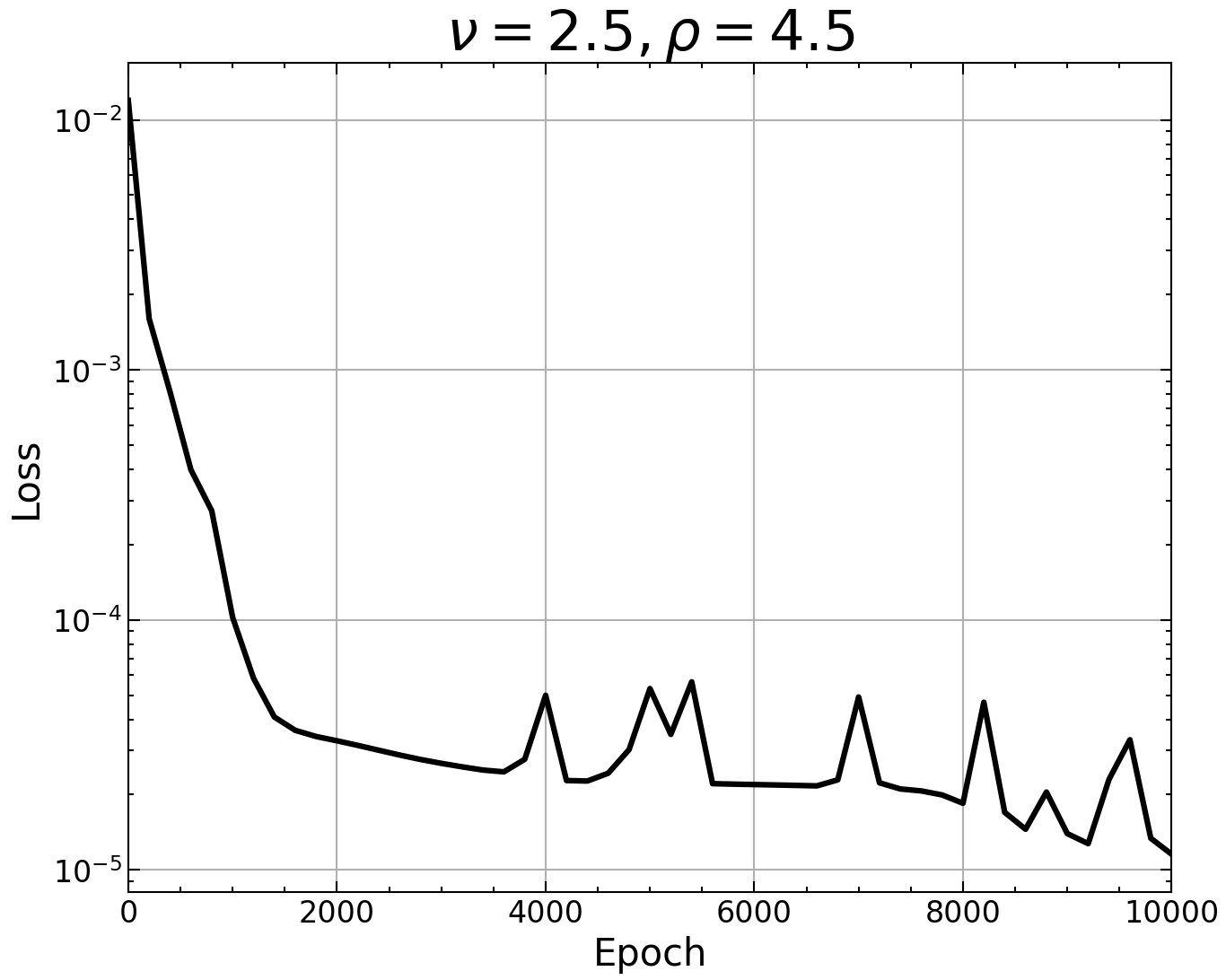}}
  \caption{Results for Section \ref{sssec:results-ppde-rd}: GPT-PINN (top, with 10 neurons) and TGPT-PINN (bottom, with 3 neurons) results for $\rho = 2.5$ and $\nu=4.5$.}
  \label{fig:RD-gpt-pinn4}
\end{figure}

\noindent {\bf GPT/TGPT-PINN results:} We set the parameter domain to be $(\rho,\nu) \in [1, 5] \times [1, 5]$, and the training set to be an $11 \times 11$ equispaced grid on this domain. Figure \ref{fig:RD-GPT-PINN} shows the histories of convergence (left) and the locations of the two hyper-reduced networks (middle) as the number of neurons increases when they are trained offline. On the right is the result when they are tested online at 49 random locaions not seen during training. It is clear that the TGPT-PINN consistently outperforms the GPT-PINN. In fact, with barely $3$ neurons, it already outperforms the GPT-PINN with $10$ neurons. With $2$ neurons, its worst case performance is about the same as the best case for GPT-PINN. We present the solutions obtained at two unseen parameter values for both the GPT-PINN (first row) and TGPT-PINN (second row) in Figure \ref{fig:RD-gpt-pinn1} and Figure \ref{fig:RD-gpt-pinn4}.

\section{Conclusion}\label{sec:conclusion}
We have introduced and investigated TGPT-PINN, a physics-informed nonlinear model reduction framework. By combining the practical efficacy of PINNs-based PDE solutions with model reduction using the GPT-PINN template and introducing new discontinuity-approximating strategies and nonlinear transform layers, the TGPT-PINN can overcome the limitations of linear model reduction in the transport-dominated regime and is effective on a wide range of practical PDE problems.

\end{document}